\documentclass[10pt,reqno,oneside]{amsproc}

\allowdisplaybreaks

\usepackage{amsmath, amsthm, 
amssymb,enumerate, enumitem, bbm,multicol}
\usepackage{marginnote}

\chardef\forshowkeys=0
\chardef\refcheck=0
\chardef\showllabel=0
\chardef\sketches=0

\ifnum\forshowkeys=1

\usepackage[notcite,color]{showkeys}

\fi

\usepackage{times}
\usepackage{graphicx}
\usepackage[usenames,dvipsnames,svgnames,table]{xcolor}
\usepackage[colorlinks=true, pdfstartview=FitV, linkcolor=black, citecolor=black, urlcolor=black]{hyperref}

\usepackage{tikz}

\usepackage{color}
\usepackage{xcolor}

\usepackage{datetime}
\usepackage{fancyhdr,lastpage}
\ifnum\refcheck=1
 \usepackage{refcheck}    %%%
\fi

\chardef\coloryes=1 %%%out
\chardef\isitdraft=0 %%%out

\ifnum\isitdraft=1 %%%outhttps://preview.overleaf.com/public/tfqhbhwccrwh/images/4f4289976560863c3073d36f0bc870df24fe1ade.jpeg
\def\eqref#1{({\ref{#1}})}                %saves writing parenthesis%%%out

\definecolor{mygray}{rgb}{.6, .6, .6}

\definecolor{refkey}{rgb}{.3,0.3,0.3}%%%out

\def\nnewpage{} %\nnewpage does nothing
\else%%%out

\def\startnewsection#1#2{\section{#1}\label{#2}\setcounter{equation}{0}}   %\starts a new section
\def\nnewpage{} %\nnewpage does nothing

\fi%%%out

\begin{document}
	\def\ques{{\colr \underline{??????}\colb}}
	\def\nto#1{{\colC \footnote{\em \colC #1}}}
	\def\fractext#1#2{{#1}/{#2}}
	\def\fracsm#1#2{{\textstyle{\frac{#1}{#2}}}}   %smaller version of frac
	\def\nnonumber{}
	\def\les{\lesssim}
	\def\lec{\lesssim}
	\def\RR{R}
	\def\JJ{J}
	\def\VV{V}
	\def\WW{W}
	\def\pt{\partial_t}

	\def\colr{\color{red}}%%%out
	\def\colg{{}}%%%out
	\def\colu{\color{blue}}%%%out
	\definecolor{coloroooo}{rgb}{0.45,0.0,0}
	\def\cole{\color{coloroooo}}
	\def\colb{\color{black}}%%%out
	\def\colA{{}}%%%out
	\def\colB{{}}%%%out
	\def\colC{{}}%%%out
	\def\colD{{}}%%%out
	\def\colE{{}}%%%out
	\def\colF{{}}%%%out

	%%3
	\ifnum\coloryes=1%%%out
	
	\definecolor{coloraaaa}{rgb}{0.1,0.2,0.8}%%%out
	\definecolor{colorbbbb}{rgb}{0.1,0.7,0.1}%%%out
	\definecolor{colorcccc}{rgb}{0.8,0.3,0.9}%%%out
	\definecolor{colordddd}{rgb}{0.0,.5,0.0}%%%out
	\definecolor{coloreeee}{rgb}{0.8,0.3,0.9}%%%out
	\definecolor{colorffff}{rgb}{0.8,0.9,0.9}%%%out%gray
	\definecolor{colorgggg}{rgb}{0.5,0.0,0.4}%%%out
	\definecolor{coloroooo}{rgb}{0.45,0.0,0}
	
	\def\colb{\color{black}}%%%out
	
	\def\colr{\color{red}}%%%out
	\def\cole{\color{coloroooo}}
	
	\def\colu{\color{blue}}%%%out
	\def\colg{\color{colordddd}}%%%out
	\def\colgray{\color{colorffff}}%%%out
	
	\def\colA{\color{coloraaaa}}%%%out
	\def\colB{\color{colorbbbb}}%%%out
	\def\colC{\color{colorcccc}}%%%out
	\def\colD{\color{colordddd}}%%%out
	\def\colE{\color{coloreeee}}%%%out
	\def\colF{\color{colorffff}}%%%out
	\def\colG{\color{colorgggg}}%%%out

	%%4
	\fi%%%out
	%\fi%%%out
	\ifnum\isitdraft=1%%%out
	\chardef\coloryes=1 %%%out
	\baselineskip=17.6pt%%%out
	%%%%%%%%%%%%%%%%%%%%
	\pagestyle{myheadings}
	%\reversemarginpar
	\def\const{\mathop{\rm const}\nolimits}  %const (for a constant)
	\def\diam{\mathop{\rm diam}\nolimits}    %diameter
	\def\rref#1{{\ref{#1}{\rm \tiny \fbox{\tiny #1}}}}
	\def\theequation{\fbox{\bf \thesection.\arabic{equation}}}
	\def\startnewsection#1#2{\newpage\colg \section{#1}\colb\label{#2}
		\setcounter{equation}{0}
		%\rfoot{\thepage}
		\pagestyle{fancy}
		\lhead{\colb Section~\ref{#2}, #1 }
		\cfoot{}
		%\rfoot{\thepage}
		\rfoot{\thepage\ of \pageref{LastPage}}
		%\lfoot{\colb{\today,~\currenttime}~(kt8, Version~\fbox{\version})}
		%\lfoot{\colb{\today,~\currenttime}~(jkl3), page \thepage\ of \pageref{LastPage}}
		\lfoot{\colb{\today,~\currenttime}~}}
	
	\chead{}
	\rhead{\thepage}
	\def\nnewpage{\newpage}
	\newcounter{startcurrpage}
	\newcounter{currpage}
	\def\llll#1{{\rm\tiny\fbox{#1}}}
	%%%%%%%%%%%%%%%%%%%%
	\def\blackdot{{\color{red}{\hskip-.0truecm\rule[-1mm]{4mm}{4mm}\hskip.2truecm}}\hskip-.3truecm}%%%out
	\def\bluedot{{\colC {\hskip-.0truecm\rule[-1mm]{4mm}{4mm}\hskip.2truecm}}\hskip-.3truecm}%%%out
	\def\purpledot{{\colA{\rule[0mm]{4mm}{4mm}}\colb}}%%%out
	\def\pdot{\purpledot}%%%out
	\else%%%out  
	\baselineskip=12.8pt
	\def\blackdot{{\color{red}{\hskip-.0truecm\rule[-1mm]{4mm}{4mm}\hskip.2truecm}}\hskip-.3truecm}%%%out
	\def\purpledot{{\rule[-3mm]{8mm}{8mm}}}%%%out
	\def\pdot{}
	%%5
	\fi%%%out
        \def\cole{\color{black}}

	%varmac       various macros
\def\and{\text{\ \ \ \ and \ \ \ \ }}      
	\def\textand{\qquad \text{and}\qquad}
	\def\pp{p}
	\def\qq{{\tilde p}}
	\def\KK{K}
		\def\TT{\tilde{T}}
	\def\MM{M}
	\def\ema#1{{#1}}
	\def\emb#1{#1}
	
	\ifnum\isitdraft=1
	\def\llabel#1{\nonumber}
	\else
	\def\llabel#1{\nonumber}
	\fi
	
	\def\tepsilon{\tilde\epsilon}
	\def\epsilonz{\epsilon_0}
	\def\restr{\bigm|}
	\def\into{\int_{\Omega}}
	\def\intu{\int_{\Gamma_1}}
	\def\intl{\int_{\Gamma_0}}
	\def\tpar{\tilde\partial}
	\def\bpar{\,|\nabla_2|}
	\def\barpar{\bar\partial}
        \def\Omegae{\Omega_{\text e}}
        \def\Gammac{\Gamma_{\text c}}
        \def\Gammaf{\Gamma_{\text f}}
        \def\Omegaf{\Omega_{\text f}}
	\def\FF{F}
	\def\gdot{{\color{green}{\hskip-.0truecm\rule[-1mm]{4mm}{4mm}\hskip.2truecm}}\hskip-.3truecm}%%%out
\def\tdot{{\color{green}{\hskip-.0truecm\rule[-.5mm]{5mm}{5mm}\hskip.2truecm}}\hskip-.2truecm}
\def\bdot{{\color{blue}{\hskip-.0truecm\rule[-.5mm]{2mm}{4mm}\hskip.2truecm}}\hskip-.2truecm}
	\def\cydot{{\color{cyan} {\hskip-.0truecm\rule[-1mm]{4mm}{4mm}\hskip.2truecm}}\hskip-.3truecm}%%%out
	\def\rdot{{\color{red} {\hskip-.0truecm\rule[-1mm]{4mm}{4mm}\hskip.2truecm}}\hskip-.3truecm}%%%out
	
	\def\nts#1{{\color{red}\hbox{\bf ~#1~}}} %nts=note to self%%%out
	
	\def\ntsr#1{\vskip.0truecm{\color{red}\hbox{\bf ~#1~}}\vskip0truecm} %nts=note to self%%%out
	
	\def\ntsf#1{\footnote{\hbox{\bf ~#1~}}} %nts=note to self
	\def\ntsf#1{\footnote{\color{red}\hbox{\bf ~#1~}}} %nts=note to self%%%out
	\def\bigline#1{~\\\hskip2truecm~~~~{#1}{#1}{#1}{#1}{#1}{#1}{#1}{#1}{#1}{#1}{#1}{#1}{#1}{#1}{#1}{#1}{#1}{#1}{#1}{#1}{#1}\\}%%%out
	\def\biglineb{\bigline{$\downarrow\,$ $\downarrow\,$}}%%%out
	\def\biglinem{\bigline{---}}%%%out
	\def\biglinee{\bigline{$\uparrow\,$ $\uparrow\,$}}%%%out
	\def\ceil#1{\lceil #1 \rceil}
	\def\gdot{{\color{green}{\hskip-.0truecm\rule[-1mm]{4mm}{4mm}\hskip.2truecm}}\hskip-.3truecm}%%%out
	\def\bluedot{{\color{blue} {\hskip-.0truecm\rule[-1mm]{4mm}{4mm}\hskip.2truecm}}\hskip-.3truecm}%%%out
	\def\rdot{{\color{red} {\hskip-.0truecm\rule[-1mm]{4mm}{4mm}\hskip.2truecm}}\hskip-.3truecm}%%%out
	\def\dbar{\bar{\partial}}
	\newtheorem{Theorem}{Theorem}[section]
	\newtheorem{Corollary}[Theorem]{Corollary}
	\newtheorem{Proposition}[Theorem]{Proposition}
	\newtheorem{Lemma}[Theorem]{Lemma}
	\newtheorem{Remark}[Theorem]{Remark}
	\newtheorem{definition}{Definition}[section]
	\def\theequation{\thesection.\arabic{equation}}
	\def\cmi#1{{\color{red}IK: #1}}
	\def\cmj#1{{\color{red}IK: #1}}
	\def\cml{\rm \colr Linfeng:~} %%%comment from Linfeng
	\def\XX{\mathbf{X}}
	\def\II{\mathcal{I}}

	%varmac (various macros)
	\def\scl{,}
	\def\sqrtg{\sqrt{g}}
	\def\DD{{\mathcal D}}
	\def\OO{\tilde\Omega}
	\def\EE{{\mathcal E}}
	\def\lot{{\rm l.o.t.}}                       %lower order terms
	\def\endproof{\hfill$\Box$\\}
	\def\square{\hfill$\Box$\\}
	\def\inon#1{\ \ \ \ \text{~~~~~~#1}}                %in or on
	\def\comma{ {\rm ,\qquad{}} }            %comma in a formula
	\def\commaone{ {\rm ,\qquad{}} }         %second comma in a formula
	\def\dist{\mathop{\rm dist}\nolimits}    %distance
	\def\ad{\mathop{\rm ad}\nolimits}    %distance
	\def\sgn{\mathop{\rm sgn\,}\nolimits}    %sghttps://preview.overleaf.com/public/tfqhbhwccrwh/images/d61f2116e5636881bf20dd6b89ed619615c006d3.jpegn
	\def\Tr{\mathop{\rm Tr}\nolimits}    %trace
	\def\dive{\mathop{\rm div}\nolimits}    %divergence
	\def\grad{\mathop{\rm grad}\nolimits}    %gradient
	\def\curl{\mathop{\rm curl}\nolimits}    %curl
	\def\det{\mathop{\rm det}\nolimits}    %det
	\def\supp{\mathop{\rm supp}\nolimits}    %support
	\def\re{\mathop{\rm {\mathbb R}e}\nolimits}    %distance
	\def\wb{\bar{\omega}}
	\def\Wb{\bar{W}}
	\def\indeq{\quad{}}                     %indentation in formulas
	\def\indeqtimes{\indeq\indeq\indeq\indeq\times} 
	\def\period{.}                           %period in a formula
	\def\semicolon{\,;}                      %semicolon in a formula
	\def\bfx{\mathbf{x}}
	\newcommand{\cD}{\mathcal{D}}
	\newcommand{\cH}{\mathcal{H}}
	%**end of header
	\newcommand{\imp}{\Rightarrow}
	\newcommand{\tr}{\operatorname{tr}}
	\newcommand{\vol}{\operatorname{vol}}
	\newcommand{\id}{\operatorname{id}}
	\newcommand{\p}{\parallel}
	\newcommand{\norm}[1]{\Vert#1\Vert}
	\newcommand{\abs}[1]{\vert#1\vert}
	\newcommand{\nnorm}[1]{\left\Vert#1\right\Vert}
	\newcommand{\aabs}[1]{\left\vert#1\right\vert}
	
	\ifnum\showllabel=1
	 \def\llabel#1{\marginnote{\color{lightgray}\rm\small(#1)}[-0.0cm]\notag}
	\else
	\def\llabel#1{\notag}
	\fi

\title[Fluid-structure interaction problem]{On the local existence of solutions to the Navier-Stokes-wave system \\with a free interface}

\author[I.~Kukavica]{Igor Kukavica}
\address{Department of Mathematics\\
	University of Southern California\\
	Los Angeles, CA 90089}
\email{kukavica@usc.edu}
%\thanks{XXX}

\author[L.~Li]{Linfeng Li}
\address{Department of Mathematics\\
	University of California Los Angeles\\
	Los Angeles, CA 90095}
\email{lli265@math.ucla.edu}
%\thanks{XXX}

\author[A.~Tuffaha]{Amjad Tuffaha}
\address{Department of Mathematics and Statistics\\
	American University of Sharjah\\
	Sharjah, UAE}
\email{atufaha@aus.edu}
%\thanks{XXX}

\begin{abstract} 
We address a system of equations modeling a compressible fluid interacting with an elastic body in dimension three. 
We prove the local existence and uniqueness of a strong solution when the initial velocity belongs to the space
$H^{2+\epsilon}$ and the initial structure velocity is in $H^{1.5+\epsilon}$, where $\epsilon \in (0,1/2)$.
\end{abstract}

\maketitle
	
\tableofcontents

\startnewsection{Introduction}{sec01}
The objective of this paper is to establish the local-in-time existence of solutions for the free boundary fluid-structure interaction model under low regularity assumptions on the initial data. The model describes the interaction between a viscous compressible fluid and an elastic structure that is immersed in it. Mathematically, the dynamics of the fluid are governed by the compressible Navier-Stokes equations in the velocity and density variables $(u,\rho)$, while the elastic dynamics are described by a second-order elasticity equation (which is replaced by a wave equation for the sake of simplicity) in the vector variables $(w,w_{t})$ representing the displacement and velocity of the structure.

The interaction between the structure and the fluid is mathematically characterized by velocity and stress matching boundary conditions at the moving interface that separates the solid and fluid regions. Since the interface position evolves with time and is unknown a~priori, this is a free-boundary problem. The problem is challenging due to the mismatch between parabolic and hyperbolic regularity, as well as the complexity of the stress-matching condition on the free boundary.

The local-in-time existence and well-posedness results for the fluid-structure interaction model have been extensively studied in the literature. In 2005, the authors of \cite{CS1, CS2} established the local-in-time existence and well-posedness for the incompressible model, using the Lagrangian coordinate system to fix the domain and the Tychonoff fixed point theorem to construct a solution, given an initial fluid velocity $u_0 \in H^5$ and structural velocity $w_1 \in H^3$. Subsequently, \cite{KT1, KT2} obtained a priori estimates for the local existence of solutions using direct estimates for the initial data, namely $u_0 \in H^3$ and $w_1 \in H^{5/2+r}$, where $r \in (0, (\sqrt{2}-1)/2)$. The authors relied on the hidden regularity trace theorem for wave equations, established in \cite{LLT, BL, L1, L2, S, T}, as a key ingredient to obtain their result.
Several works on wave-heat coupled systems on a non-moving domain have contributed to the understanding of the heat-wave interaction phenomena (cf.~\cite{ALT, AT1, AT2, DGHL, BGLT1, BGLT2, KTZ1, KTZ2, KTZ3, LL1}). Recently, Raymond and Vanninathan \cite{RV} obtained a sharp regularity result for the case when the initial domain is a flat channel. They studied the system in the Lagrangian coordinate setting and obtained local-in-time solutions for the 3D model, with the initial velocity $u_0 \in H^{1+\alpha}$ and the initial structural velocity $w_1 \in H^{1/2+\alpha+\beta}$, where $\alpha \in (1/2, 1)$ and $\beta > 0$.
In \cite{BGT}, Boulakia, Guerrero, and Takahashi obtained a unique local-in-time solution for the general domain case, given the initial data $u_0 \in H^2$ and $w_1 \in H^{9/8}$.

The compressible model under consideration was first treated in  \cite{BG1}, where the authors obtained the existence and uniqueness for the initial density $\rho_0$ belonging to $H^3$, the velocity $u_0$ in $H^4$, and the structure displacement and velocity $(w, w_t)$ in $H^3\times H^2$. 
A similar result was later obtained by Kukavica and Tuffaha \cite{KT3} with less regular initial data $(\rho_0, u_0, w_1) \in H^{3/2+r} \times H^3 \times H^{3/2+r}$, where $r\in (0, (\sqrt{2} -1)/2)$.  In \cite{BG2}, the existence of a regular global solution is proved for small initial data.  In a recent work \cite{BG3}, the authors proved the existence of a unique local-in-time strong solution of the interaction problem between a compressible fluid and elastic structure for initial data $(\rho_0, u_0, w_1) \in H^3 \times H^6\times H^3$, where the elastic structure is modeled by the Saint-Venant Kirchhoff system.
For some other works on fluid-structure models, cf.~\cite{AL,B,BS, BuL,BZ1,BZ2,BTZ,DEGL,F,GH,GGCC,GGCCL,IKLT1,IKLT2,KKLTTW, KMT,KOT,LL1,LL2,LT,LTr1,LTr2,MC1,MC2,MC3,SST, Tr}.

In this paper, we provide a natural proof of the existence of a unique local-in-time solution to the system under a low regularity assumptions $u_{0} \in H^{2+\epsilon}$ and $w_{1} \in H^{1.5+\epsilon}$, where $\epsilon \in (0,1/2)$, in the case of the flat initial configuration. 
Our proof relies on a maximal regularity type theorem for the nonhomegeneous linear parabolic problem with Neumann type conditions on the fluid-structure interface, in addition to the hidden regularity theorems (cf.~Lemmas~\ref{L03}--\ref{L13}) for the wave equation.  The time regularity of the solution is obtained using the energy estimates, which, combined with the elliptic regularity, yield the spatial regularity of the solutions.  An essential ingredient of the proof of the main results is a trace inequality 
\begin{align} 
	\begin{split} \Vert u\Vert_{H^{\theta}((-\infty, \infty) \scl L^{2}(\Gammac))} 
		\les 
		\Vert u \Vert_{L^{2}((-\infty, \infty)\scl H^{r}(\Omegaf))}^{1/(2r+1)} \Vert u \Vert_{H^{2\theta r/(2r-1)}((-\infty, \infty)\scl L^{2}(\Omegaf))}^{2r/(2r+1)} 
		+ 
		\Vert u \Vert_{L^{2}((-\infty, \infty)\scl H^{r}(\Omegaf))}, 
	\end{split} 
\llabel{EQ137a} 
\end{align} for functions which are Sobolev in the time variable and square integrable on the boundary (cf.~Lemma~\ref{L06} and \eqref{EQ137} below).  This is used essentially in the proof of the existence for the nonlinear parabolic-wave system, Theorem~\ref{T03}, and in the proof of the main result, Theorem~\ref{T01}. The construction of a unique solution for the fluid-structure problem is obtained via the Banach fixed point theorem. The scheme involves solving the nonlinear parabolic-wave system with the variable coefficients treated as a given forcing perturbations.

One of the essential difficulties in establishing the existence of solutions is that the constants in the inequality are inversely proportional to powers of time $T$, which poses a problem for establishing convergence of a fixed-point scheme for small time. The same issue with the growing constants also arises in the hidden regularity inequalities in Lemmas~\ref{L03}--\ref{L13} for the wave equation. 
We overcome this difficulty by solving a modified system which is posed on the fixed time interval $(0,1]$. 
As opposed to the velocity matching boundary condition \eqref{EQ262} in the original fluid-structure interaction problem, we impose the integrated velocity matching boundary condition \eqref{EQ73} on the unit time interval in the modified system.
These two boundary conditions agree on a small time interval and thus the modified system agrees with the original system when restricted to a small time interval.
In the integrated velocity matching boundary condition \eqref{EQ73}, an important ingredient is the cutoff function in time that depends on a variable time $\tilde{T}$, which is then chosen to be less than a fixed time $T_0$, allowing for contraction estimates on the solution map. 
Another major difficulty is the handling of the normal derivative of the elastic structure on the common boundary, which is estimated by appealing to the hidden trace regularity (see Lemma~\ref{L13}). 
The main issue with proving the fixed-point theorems (for the linear and nonlinear variants) is that time derivatives, which are frequently fractional, fall on the cutoff, showing that the constant dependence on $\TT$ needs to be treated carefully.

Similarly, for the nonlinear system, treated in Section~\ref{sec06}, we also need to modify the definition of the Lagrangian map and the variable coefficient matrix using a cutoff in time function to ensure similar contraction-type estimates on the solution map for the system with given variable coefficients. The solution in each iteration step is used to prescribe new variable coefficients for the next iteration step. The contracting property of the Navier-Stokes-wave system is maintained by taking a sufficiently short time $\tilde{T}$ to ensure closeness of the Jacobian and the inverse matrix of the flow map to their initial states.

Note that the configuration we adopt, \eqref{EQ313} with the periodic boundary conditions in the $y_1$ and $y_2$ directions, is needed only in Lemma~\ref{L13}. In these estimates, Sobolev time norms pose a particular challenge when the cutoff function is involved since they involve singular terms in $\tilde{T}$ that have to be compensated by taking sufficiently high $L^{p}$ norms of time derivatives of~$v$.

The paper is structured as follows. In Section~\ref{sec_setting}, we introduce the fluid-structure model and state our main result. Next, in Section~\ref{sec03}, we present the trace inequality, interpolation, and hidden regularity lemmas. Section~\ref{sec04} provides the maximal regularity for the nonhomogeneous parabolic problem, which is a crucial ingredient in the proof of local existence for the nonlinear parabolic-wave system, discussed in Section~\ref{sec05}. Finally, in Section~\ref{sec06}, we prove our main result, Theorem~\ref{T01}, using the local existence result established in Section~\ref{sec05} and constructing a unique solution via the Banach fixed point theorem.

\startnewsection{The model and main results}{sec_setting}
We consider the fluid-structure problem for a free boundary system involving the motion of an elastic body immersed in a compressible fluid. 
Let $\Omegaf (t)$ and $\Omegae (t)$ be the domains occupied by the fluid and the solid body at time $t$ in $\mathbb{R}^3$, whose common boundary is denoted by~$\Gammac(t)$.
The fluid is modeled by the compressible Navier-Stokes equations,
which in Eulerian coordinates reads
\begin{align}
	&
	\rho_t
	+
	\dive (\rho u)
	=
	0
	\inon{in~$[0,T]\times \Omegaf(t)$}  ,
	\label{EQ21}
	\\
	&
	\rho u_t
	+
	\rho (u \cdot \nabla)u
	-
	\lambda \dive (\nabla u + (\nabla u)^T)
	-
	\mu
	\nabla \dive u
	+
	\nabla p
	=
	0
	\inon{in~$[0,T] \times \Omegaf(t)$}
	,
	\label{EQ22}
\end{align}
where $\rho = \rho (t,x) \in \mathbb{R}_+$ is the density, $u = u(t,x) \in \mathbb{R}^3$ is the velocity, $p = p(\rho (t,x)) \in \mathbb{R}_+$ is the pressure, and $\lambda, \mu>0$ are physical constants.
(We remark that the condition for $\lambda$ and $\mu$ can be relaxed to $\lambda>0$ and $3\lambda + 2\mu >0$.)
The system \eqref{EQ21}--\eqref{EQ22} is defined on $\Omegaf (t)$ which set to $ \Omegaf =\Omegaf (0) $ and evolves in time.
The dynamics of the coupling between the compressible fluid and the elastic body are best described in the Lagrangian coordinates. 
Namely, we introduce the Lagrangian flow map $\eta(t,\cdot) \colon \Omegaf \to \Omegaf (t)$ and
rewrite the system \eqref{EQ21}--\eqref{EQ22} as
\begin{align}
	&
	\RR_t 
	-
	\RR a_{kj} \partial_k v_j
	=
	0
	\inon{in~$[0,T]\times \Omegaf$}
	,
	\label{EQ260}
	\\
	&
	\partial_t v_j 
	- 
	\lambda \RR a_{kl}
	\partial_k (a_{ml} \partial_m v_j + a_{mj} \partial_m v_l)
	-
	\mu \RR a_{kj} \partial_k(a_{mi} \partial_m v_i )
	+
	\RR a_{kj} \partial_k (q (\RR^{-1} ))
	=
	0
	\inon{in~$[0,T]\times \Omegaf$},
	\label{EQ261}
\end{align}
for $j=1,2,3$, where $\RR (t,x)= \rho^{-1} (t, \eta(t,x))$ is the reciprocal of the Lagrangian density, $v(t,x) = u(t, \eta(t,x))$ is the Lagrangian velocity, $a(t,x)= (\nabla \eta (t,x))^{-1}$ is the inverse matrix of the flow map and $q$ is a given function of the density.
The system \eqref{EQ260}--\eqref{EQ261} is expressed in terms of
Lagrangian coordinates and posed in a fixed domain~$\Omegaf$.

On the other hand, the elastic body is modeled by the wave equation
in Lagrangian coordinates, which is posed in a fixed domain $\Omegae$ as
\begin{align}
	w_{tt}
	-
	\Delta w
	=
	0
	\inon{in~$[0,T] \times \Omegae$}
	,
	\label{EQ23}
\end{align}
where $(w, w_t)$ are the displacement and the structure velocity.
The interaction boundary conditions are the velocity and stress matching conditions, which are formulated in Lagrangian coordinates over the fixed common boundary $\Gammac = \Gammac(0)$ as
\begin{align}
	&
	v_j
	= 
	\partial_t w_{j}
	\inon{on~$[0,T] \times \Gammac$},
	\label{EQ262}
	\\&
	\partial_k w_j \nu^k
	=
	\lambda \JJ a_{kl} 
	(a_{ml} \partial_m v_j 
	+
	a_{mj} \partial_m v_l
	) \nu^k
	+
	\mu \JJ a_{kj} a_{mi} \partial_m v_i \nu^k
	-
	\JJ 
	a_{kj} q \nu^k
	\inon{on~$[0,T] \times \Gammac$}
	,
	\label{EQ263}
\end{align}
for $j=1,2,3$, where $J(t,x) = \det (\nabla \eta(t,x))$ is the Jacobian and $\nu$ is 
the unit normal vector to~$\Gammac$, which is outward with respect to~$\Omegae$.
In the present paper, we consider the reference configurations $\Omega = \Omegaf  \cup \Omegae  \cup \Gammac$, $\Omegaf$, and $\Omegae$ given by (see figure 1)
  \begin{align}
	\begin{split}
		&
		\Omega 
		=
		\{y=(y_1,y_2,y_3) \in \mathbb{R}^3 
		:
		(y_1, y_2) \in \mathbb{T}^2, 0<y_3< L_3\}
		,
		\\
		&	
		\Omegaf
		=
		\{y=(y_1,y_2,y_3) \in \mathbb{R}^3 :
		(y_1, y_2) \in \mathbb{T}^2, 0<y_3< L_1  ~\text{~or~}~ L_2<y_3< L_3 \},
		\\&
		\Omegae
		=
		\{y=(y_1,y_2,y_3) \in \mathbb{R}^3 :
		(y_1, y_2) \in \mathbb{T}^2, L_1<y_3< L_2\}
		,
	\end{split}
	\label{EQ313}
  \end{align}
where $0<L_1 <L_2 <L_3$ and $\mathbb{T}^2$ is the two-dimensional
torus with the side~$2\pi$.
Thus, the common boundary is expressed as
\begin{align}
	\Gammac
	=
	\{(y_1, y_2) \in \mathbb{R}^2 :
	(y_1, y_2, y_3) \in \Omega, y_3= L_1 ~\text{~or~}~ y_3= L_2\}
        ,
	\llabel{EQ360}
\end{align}
while the outer boundary is represented by
\begin{align}
	\Gammaf 
	= 
	\{y \in \bar{\Omega} : y_3 = 0\} \cup \{y\in \bar{\Omega} : y_3 = L_3\}
	.
	\llabel{EQ318}
\end{align}

\tikzset{every picture/.style={line width=0.75pt}} %set default line width to 0.75pt        

\tikzset{every picture/.style={line width=0.75pt}} %set default line width to 0.75pt        
\begin{center}

\begin{tikzpicture}[x=0.75pt,y=0.75pt,yscale=-1,xscale=1]

	%Shape: Cube [id:dp8791417391495022] 
	\draw  [fill={rgb, 255:red, 255; green, 255; blue, 255 }  ,fill opacity=1 ] (106,172.8) -- (117.6,161.2) -- (226,161.2) -- (226,207.6) -- (214.4,219.2) -- (106,219.2) -- cycle ; \draw   (226,161.2) -- (214.4,172.8) -- (106,172.8) ; \draw   (214.4,172.8) -- (214.4,219.2) ;
	%Shape: Cube [id:dp8092975371067813] 
	\draw  [fill={rgb, 255:red, 255; green, 255; blue, 255 }  ,fill opacity=1 ] (106,126.4) -- (117.6,114.8) -- (226,114.8) -- (226,161.2) -- (214.4,172.8) -- (106,172.8) -- cycle ; \draw   (226,114.8) -- (214.4,126.4) -- (106,126.4) ; \draw   (214.4,126.4) -- (214.4,172.8) ;
	%Shape: Cube [id:dp011028565419338499] 
	\draw  [fill={rgb, 255:red, 255; green, 255; blue, 255 }  ,fill opacity=1 ] (106,80) -- (117.6,68.4) -- (226,68.4) -- (226,114.8) -- (214.4,126.4) -- (106,126.4) -- cycle ; \draw   (226,68.4) -- (214.4,80) -- (106,80) ; \draw   (214.4,80) -- (214.4,126.4) ;
	%Straight Lines [id:da6609995120401413] 
	\draw    (243,144) -- (353,144) ;
	\draw [shift={(355,144)}, rotate = 180] [color={rgb, 255:red, 0; green, 0; blue, 0 }  ][line width=0.75]    (10.93,-3.29) .. controls (6.95,-1.4) and (3.31,-0.3) .. (0,0) .. controls (3.31,0.3) and (6.95,1.4) .. (10.93,3.29)   ;
	%Shape: Cube [id:dp08854211018077285] 
	\draw  [fill={rgb, 255:red, 255; green, 255; blue, 255 }  ,fill opacity=1 ] (364,79) -- (374.6,68.4) -- (484,68.4) -- (484,207.4) -- (473.4,218) -- (364,218) -- cycle ; \draw   (484,68.4) -- (473.4,79) -- (364,79) ; \draw   (473.4,79) -- (473.4,218) ;
	%Curve Lines [id:da4662873372033345] 
	\draw    (364,117) .. controls (396,89) and (436,150) .. (473.4,117) ;
	%Curve Lines [id:da6966672402256087] 
	\draw    (364,167) .. controls (403,190) and (434,144) .. (473.4,167) ;
	%Curve Lines [id:da056835073149755555] 
	\draw    (473.4,117) .. controls (477,107) and (482,114) .. (484,106) ;
	%Curve Lines [id:da5621186811315967] 
	\draw    (473.4,167) .. controls (481,150) and (477,181) .. (484,160) ;
	%Straight Lines [id:da125335906370873] 
	\draw    (53,230) -- (53,176) ;
	\draw [shift={(53,174)}, rotate = 90] [color={rgb, 255:red, 0; green, 0; blue, 0 }  ][line width=0.75]    (10.93,-3.29) .. controls (6.95,-1.4) and (3.31,-0.3) .. (0,0) .. controls (3.31,0.3) and (6.95,1.4) .. (10.93,3.29)   ;
	
	% Text Node
	\draw (251,123) node [anchor=north west][inner sep=0.75pt]   [align=left] {{\scriptsize under the flow map $\displaystyle \eta $}};
	% Text Node
	\draw (32,176) node [anchor=north west][inner sep=0.75pt]   [align=left] {$\displaystyle y_{3}$};
	% Text Node
	\draw (77,72) node [anchor=north west][inner sep=0.75pt]   [align=left] {$\displaystyle L_{3}$};
	% Text Node
	\draw (77,118) node [anchor=north west][inner sep=0.75pt]   [align=left] {$\displaystyle L_{2}$};
	% Text Node
	\draw (77,165) node [anchor=north west][inner sep=0.75pt]   [align=left] {$\displaystyle L_{1}$};
	% Text Node
	\draw (77,210) node [anchor=north west][inner sep=0.75pt]   [align=left] {$\displaystyle 0$};
	% Text Node
	\draw (148,96) node [anchor=north west][inner sep=0.75pt]   [align=left] {$\displaystyle \Omega _{f}$};
	% Text Node
	\draw (148,188) node [anchor=north west][inner sep=0.75pt]   [align=left] {$\displaystyle \Omega _{f}$};
	% Text Node
	\draw (148,144) node [anchor=north west][inner sep=0.75pt]   [align=left] {$\displaystyle \Omega _{e}$};
	% Text Node
	\draw (404,89) node [anchor=north west][inner sep=0.75pt]   [align=left] {$\displaystyle \Omega _{f}( t)$};
	% Text Node
	\draw (404,186) node [anchor=north west][inner sep=0.75pt]   [align=left] {$\displaystyle \Omega _{f}( t)$};
	% Text Node
	\draw (404,134) node [anchor=north west][inner sep=0.75pt]   [align=left] {$\displaystyle \Omega _{e}( t)$};
	% Text Node
	\draw (160,35) node [anchor=north west][inner sep=0.75pt]   [align=left] {Figure~1. Lagrangian domain to Eulerian domain};
\end{tikzpicture}
\end{center}
To close the system, we impose the homogeneous Dirichlet boundary condition 
\begin{align}
	&
	v = 0
	\inon{on~$[0,T]\times \Gammaf$}
	\label{EQ266}
\end{align}
on the outer boundary $\Gammaf$ and the periodic boundary conditions for $w$, $\rho$, and $u$ on the lateral boundary, i.e.,
  \begin{align}
	w(t, \cdot), \rho(t, \eta(t, \cdot)), u(t, \eta(t, \cdot))
	~~\text{periodic in the $y_1$~and~$y_2$~directions}
	.
	\label{EQ267}
\end{align}
Note that the inverse matrix of the flow map $a$ satisfies the ODE system
\begin{align}
	&
        a_t (t,x)
	=
	-a(t,x) \nabla v(t,x)  a(t,x)
	\inon{in~$[0,T] \times \Omegaf$}
	,
	\label{EQ204}
	\\
	&
	a (0) 
	=
	\mathbb{I}_3
	\inon{in~$\Omegaf$}
	,
	\label{EQ264}
\end{align}
where $\mathbb{I}_3$ is the three-dimensional identity matrix, while the Jacobian satisfies the ODE system
\begin{align}
  \begin{split}
	&
	\JJ_t (t,x)
	=
	\JJ (t,x) a_{kj} (t,x) \partial_k v_j (t,x) 
	\inon{in~$[0,T] \times \Omegaf$}
	,
	\\&
	\JJ(0) = 1
	\inon{in~$\Omegaf$}
	.
      \end{split}
 	\label{EQ210}
     \end{align}
The initial data of the system \eqref{EQ260}--\eqref{EQ23} is given as
\begin{align}
	\begin{split}
		&
	(\RR, v, w, w_t)(0)
	=
	(\RR_0, v_0, w_0, w_1)
	\inon{in~$\Omegaf\times \Omegaf
	\times \Omegae\times \Omegae$}
	,	
	\\&
	(\RR_0, v_0, w_0, w_1)
	~\text{periodic in the $y_1$ and $y_2$ directions}
        ,
	\label{EQ265}
	\end{split}
\end{align}
where $w_0 = 0$.
For $T>0$, we denote 
\begin{align}
	H^{r,s} ((0,T) \times \Omegaf)
	=
	H^r ((0,T), L^2 (\Omegaf)) 
	\cap
	L^2 ((0,T), H^s (\Omegaf))
	,
   \llabel{EQ49}
\end{align}
with the corresponding norm
\begin{align}
	\Vert f\Vert_{H^{r,s} ((0,T) \times \Omegaf)}^2
	=
	\Vert f\Vert_{H^r ((0,T), L^2 (\Omegaf))}^2
	+
	\Vert f \Vert_{L^2 ((0,T), H^s (\Omegaf))}^2,
   \llabel{EQ37}
\end{align}
where $r, s\geq 0$ are constant parameters.
In Sections~\ref{sec04}--\ref{sec06}, we shall work on a modified system with $T=1$ to avoid issue with dependence of constants on small time.
For simplicity of notation, we write
  \begin{equation}
   \Vert f\Vert_{H^{r,s}} = \Vert f\Vert_{H^{r,s} ((0,1)\times \Omegaf)}
   \and
   \Vert f\Vert_{H^{r}_t H_x^{s}} 
   = 
   \Vert f\Vert_{H^{r}  ((0,1), H^{s} (\Omegaf))}
   .
   \llabel{EQ79}
  \end{equation}
It is also convenient to abbreviate
\begin{align}
	K^s= H^{s/2, s} 
	,
   \llabel{EQ47}
\end{align}
where the domain of integration is $(0,1) \times \Omegaf$ unless stated otherwise.
Similarly, for the analogous space of functions defined on the boundary $\Gammac$, we write
  \begin{equation}
   \Vert f\Vert_{H^{r}_t H_x^{s} (\Gammac)} 
   = 
      \Vert f\Vert_{H^{r}  ((0,1), H^{s} (\Gammac))}
   \llabel{EQ87}
  \end{equation}
and abbreviate  
\begin{align}
	K^s_{\Gammac} 
	= 
	H^{s/2, s} (\Gammac),
   \llabel{EQ107}
\end{align}
where the domain of integration is $(0,1)\times \Gammac$ unless stated otherwise.
We emphasize that the time domain of integration in the norms is $(0,1)$ when not indicated.

Our main result states the local-in-time existence of solution to the system \eqref{EQ260}--\eqref{EQ23} with the mixed boundary conditions \eqref{EQ262}--\eqref{EQ267} and the initial data~\eqref{EQ265}.

\begin{Theorem}
	\label{T01}
Let $s\in (2, 2+ \epsilon_0]$ for $\epsilon_0 \in (0,1/2)$.
Assume that $\RR_0 \in H^s (\Omegaf)$, $\RR_0^{-1} \in H^s (\Omegaf)$, $w_1 \in H^{s-1/2} (\Omegae)$,  $v_0 \in H^s (\Omegaf)$,  $v_0 |_{\Gammac} \in
H^{s+1/2} (\Gammac)$, $\partial_3 v_0 |_{\Gammaf} \in
H^{s-1/2} (\Gammaf)$,
and $w_0 = 0$, with the compatibility conditions
\begin{align}
	&
	w_{1j}
	=
	v_{0j}
	\inon{on~$\Gammac$}
	,
	\llabel{EQ270}
	\\&
	v_{0j} = 0
	\inon{on~$\Gammaf$}
	,
	\llabel{EQ271}
	\\&
	\lambda (\partial_k v_{0j} 
	+
	\partial_j v_{0k}
	)\nu^k
	+
	\mu	\partial_i v_{0i} \nu^j 
	- 
	q(\RR_0^{-1}) \nu^j 
	=
	0
	\inon{on~$\Gammac$}
	,
	\llabel{EQ272}
	\\&
	\lambda 
	\partial_k (\partial_k v_{0j}
	+ 
	 \partial_j v_{0k})
	+
	\mu \partial_j \partial_k v_{0k} 
	-
	 \partial_k (q (\RR^{-1}_0 ))
	=
	0
	\inon{on~$\Gammaf$}
	,
	\llabel{EQ273}
\end{align}
for $j=1,2,3$.
Then the system \eqref{EQ260}--\eqref{EQ23} with the coupling conditions \eqref{EQ262}--\eqref{EQ263}, boundary conditions \eqref{EQ266}--\eqref{EQ267}, and the initial data \eqref{EQ265}
admits a unique solution
\begin{align}
  \begin{split}
	&
	v \in 
	 K^{s+1} ((0,T)\times \Omegaf)
	\\&
	\RR \in H^1 ((0,T), H^s (\Omegaf))
	\\&
	w \in C([0,T], H^{s+1/4 -\epsilon_0} (\Omegae))	
	\\&
	w_t \in C([0,T], H^{s-3/4 -\epsilon_0}(\Omegae))	
        ,
  \end{split}
   \llabel{EQ109}
\end{align}
for some constant $T >0$, where the corresponding norms are bounded by a function of the norms of the initial data.
\end{Theorem}

\begin{Remark}
\label{R01}
{\rm
We assume $v_0 \in H^s (\Omegaf)$ for $s\in (2, 2+\epsilon_0]$ where $\epsilon_0 >0$, since the elliptic regularity for $\Vert v\Vert_{L^2_t H^4_x}$ in \eqref{EQ67} requires that $\RR^{-1} \in L^\infty  ((0,T), H^2 (\Omegaf))$. 
From the density equation \eqref{EQ260}, we deduce that the regularity for the initial velocity must be at least in $H^2 (\Omegaf)$,
showing the optimality of the range $s\geq2$. 
It would be interesting
to find whether the statement of the theorem holds for the borderline
case~$s=2$.
\square
}
\end{Remark}

The proof of the theorem is given in Section~\ref{sec06} below.
For simplicity, we present the proof for the pressure law $q(\RR) = \RR$, noting that the case for smooth function $q(R)$ follows completely analogously using the Sobolev and H\"older's inequalities and~\eqref{EQ98}. 
See Remark~\ref{pressure} below for necessary modifications.
\colb

\startnewsection{Space-time trace, interpolation, and hidden regularity inequalities}{sec03}
In this section, we provide several auxiliary results needed in the fixed point arguments.
The first lemma provides an estimate for the trace in a space-time norm and is an essential ingredient when constructing
solutions to the nonlinear parabolic-wave system in Section~\ref{sec05} below.

\cole
\begin{Lemma}
	\label{L06}
	Let $r>1/2$ and $\theta\geq0$.
	If 
	$u\in 
	L^{2}((-\infty,\infty)\scl H^{r}(\Omegaf))
	\cap
	H^{2\theta r/(2r-1)} ((-\infty,\infty), L^{2}(\Omegaf))$, then 
	$
	u \in H^{\theta} ((-\infty,\infty), L^{2}(\Gammac))   
	$, and for all $\epsilon \in(0,1]$,
	we have the inequality
	\begin{align}
		\begin{split}
			\Vert u\Vert_{H^{\theta}((-\infty,\infty)\scl  L^{2}(\Gammac))}
			\leq
			\epsilon
			\Vert u \Vert_{H^{2\theta r/(2r-1)}((-\infty,\infty)\scl  L^{2}(\Omegaf))}
			+ 
			C\epsilon^{1-2r}
			\Vert u \Vert_{L^{2}((-\infty,\infty)\scl  H^{r}(\Omegaf))} 
			,
		\end{split}
		\label{EQ65a}
	\end{align}
	where $C>0$ is a constant.
\end{Lemma}
\colb

The above lemma 
was proven in \cite{G}, where moreover, the interpolation spaces
were identified. Since in this paper, we only use the inequality
\eqref{EQ65a}, which allows a simpler proof, we provide an elementary argument below.

First, however, we point out a consequence when restricting the above
result to a finite time interval.

\cole%%%out
\begin{Corollary}
	\label{C01}
	Let $r>1/2$, $\theta\geq0$, and~$T>0$.
	If 
	$u\in 
	L^{2}((0,T) \scl H^{r}(\Omegaf))
	\cap
	H^{2\theta r/(2r-1)} ((0,T), L^{2}(\Omegaf))$, then 
	$
	u \in H^{\theta} ( (0,T), L^{2}(\Gammac))
	$, and for all $\epsilon \in(0,1]$,
	we have the inequality
	\begin{align}
		\begin{split}
			\Vert u\Vert_{H^{\theta}( (0,T)\scl  L^{2}(\Gammac))}
			\leq
			\epsilon
			\Vert u \Vert_{H^{2\theta r/(2r-1)}( (0,T)\scl  L^{2}(\Omegaf))}
			+ 
			C\epsilon^{1-2r}
			\Vert u \Vert_{L^{2}( (0,T) \scl  H^{r}(\Omegaf))} 
			,
		\end{split}
		\label{EQ65}
	\end{align}
	where $C>0$ is a constant, which depends on $\Omegaf$ and~$T$.
\end{Corollary}
\colb%%%out

The inequality \eqref{EQ65} follows from
Lemma~\ref{L06} using the Sobolev extension operator.
Clearly, the constant is uniform as $T\to\infty$, but may increase to
infinity as~$T\to0$.

\begin{proof}[Proof of Lemma~\ref{L06}]
	%\ntsf{here we provide the proof for the general domain}
	It is sufficient to prove \eqref{EQ65}
	for  $u\in C_{0}^{\infty}({\mathbb R} \times {\mathbb R}^3)$
	with the trace taken on the set
	\begin{equation}
		\Gamma=\bigl\{
		(t,x_1, x_2, x_3)\in{\mathbb R} \times {\mathbb R}^3
		:
		x_3=0
		\bigr\}
		;
		\llabel{EQ322}
	\end{equation}
	the general case is settled by the partition of unity and straightening of the boundary.
	Since it should be clear from the context, we usually do not distinguish in notation between 
	a function and its trace.
	Denoting by $\hat u$ the Fourier transform of $u$ with respect to
	$(t, x_1,x_2,x_3)$, we have
	\begin{equation}
		\Vert u \Vert^{2}_{H^{\theta}((-\infty,\infty)\scl  L^{2}(\Gamma))}
		\les
		\int_{-\infty}^{\infty} \int_{-\infty}^{\infty} \int_{-\infty}^{\infty}
		(1+ \tau^{2})^{\theta}
		\left|\int_{-\infty}^\infty\hat{u}(\xi_{1},\xi_{2}, \xi_3, \tau)\,d\xi_3
		\right|^{2} \,  d \tau \, d \xi_{1} \,  d
		\xi_{2}
		.
		\llabel{EQ226}
	\end{equation}
	Denote by
	\begin{equation}
		\gamma=\frac{2r-1}{2\theta}
		\label{EQ114}
	\end{equation}
	the quotient between the exponents $r$ and $2\theta r/(2r-1)$ in~\eqref{EQ65}. 
	Then, with $\lambda>0$ to be determined below, we have
	\begin{align}
		\begin{split}
			&
			\Vert u \Vert^{2}_{H^{\theta}((-\infty,\infty)\scl  L^{2}(\Gamma))} 
			\les
			\int_{\mathbb{R}^{3}}
			(1+ \tau^{2})^{\theta}
			\left|
			\int_{-\infty}^{\infty} \hat{u}(\xi_{1},\xi_{2}, \xi_{3}, \tau) \, d \xi_{3} 
			\right|^{2}
			\,  d \tau \, d \xi_{1} \,  d \xi_{2}
			\\&\indeq
			\les
			\int_{\mathbb{R}^{3}}
			(1+ \tau^{2})^{\theta}
			\left(
			\int^{\infty}_{-\infty} \frac{(1+(\xi^{2}_{1}+ \xi^{2}_{2})^{\gamma}
				+ \epsilon^{-2} \xi^{2\gamma}_{3}+\tau^{2})^{\lambda/2}}{(1+(\xi^{2}_{1}  + \xi^{2}_{2})^{\gamma} +\epsilon^{-2} \xi^{2\gamma}_{3}+\tau^{2})^{\lambda/2}} |\hat{u}| \, d \xi_{3} \right)^{2} \,  d \tau \, d \xi_{1} \, d   \xi_{2} 
			\\&\indeq
			\les 
			\int_{\mathbb{R}^{3}}
			(1+ \tau^{2})^{\theta}
			\left(\int^{\infty}_{-\infty}
			\bigl(1+(\xi^{2}_{1}+ \xi^{2}_{2})^{\gamma}+ \epsilon^{-2} \xi^{2\gamma}_{3}+\tau^{2}\bigr)^{\lambda}
			|\hat{u}|^{2} \, d \xi_{3}\right) 
			\\&\indeq\indeq\indeq\indeq \times 
			\left( 
			\int_{-\infty}^{\infty} 
			\frac{d\xi_3}
			{(1+(\xi^{2}_{1}+ \xi^{2}_{2})^{\gamma}+ \epsilon^{-2} \xi^{2\gamma}_{3}+\tau^{2})^{\lambda}} 
			\right)\,  d \tau \, d \xi_{1} \, d
			\xi_{2}
			,
		\end{split}
   \llabel{EQ209}
	\end{align}
	where we used the Cauchy-Schwarz inequality in~$\xi_3$.
Using a substitution, we have
	\begin{equation}
		\int_{-\infty}^{\infty}
		\frac{dx}{(A^2+\epsilon^{-2}x^{2\gamma})^{\lambda}}
		\les
		\epsilon^{1/\gamma} A^{1/\gamma-2\lambda}
		\comma A,\epsilon>0
		,
		\label{EQ315}
	\end{equation}
	provided $\lambda$ satisfies
	$
	2 \gamma \lambda >1
	$,
	which is by \eqref{EQ114} equivalent to
	\begin{equation}
		\lambda 
		> 
		\frac{\theta}{2r-1}
		.
		\label{EQ392}
	\end{equation}
	Note that $2\gamma\lambda>1$ implies  $1/\gamma-2\lambda<0$
	for the exponent of $A$ in~\eqref{EQ315}.
	Now we use \eqref{EQ315} for the integral in $\xi_3$
	with 
	$A=(1+(\xi^{2}_{1}+ \xi^{2}_{2})^{\gamma}+ \tau^{2})^{1/2}$, 
	while noting that
	\begin{align}
		\begin{split}
			(1+\tau^2)^{\theta} A^{1/\gamma-2\lambda}
			&
			=
			\frac{
				(1+\tau^2)^{\theta}
			}{
				(1+(\xi_1^{2}+ \xi_2^2)^{\gamma}+ \tau^{2})^{\lambda-1/2\gamma}
			}
			\leq
			(1+\tau^2)^{\theta-\lambda+1/2\gamma}
			\\&
			\leq
			(1+(\xi^{2}_{1}+ \xi^{2}_{2})^{\gamma}+ \epsilon^{-2} \xi^{2\gamma}_{3}+\tau^{2})^{\theta-\lambda+1/2\gamma}
			,
		\end{split}
		\llabel{EQ316}
	\end{align}
	provided $   \lambda-\fractext{1}{2\gamma}\leq \theta$, i.e.,
	\begin{equation}
		\lambda
		\leq 
		\frac{2 r \theta}{2r-1}
		.
		\label{EQ118}
	\end{equation}
	Under the condition \eqref{EQ118}, we thus obtain
	\begin{align}
		\begin{split}
			& \Vert u \Vert^{2}_{H^{\theta}((-\infty,\infty)\scl  L^{2}(\Gammac))}
			\les 
			\epsilon^{1/\gamma}
			\int_{\mathbb{R}\times{\mathbb R}^{3}} 
			\bigl(1+(\xi^{2}_{1}+ \xi^{2}_{2})^{\gamma}
			+ 
			\epsilon^{-2} \xi^{2\gamma}_{3}+\tau^{2}\bigr)^{\theta+1/2\gamma}
			|\hat{u}|^{2} \, d \xi_{3} \, d \xi_{2}\, d    \xi_{1} \, d \tau 
			\\&\indeq\indeq 
			\les 
			\epsilon^{1/\gamma}
			\int_{\mathbb{R}\times{\mathbb R}^{3}} 
			\bigl(1+\epsilon^{-2}(\xi^{2}_{1}+ \xi^{2}_{2}+\xi_{3}^2)^{\gamma}+\tau^{2}\bigr)^{\theta+1/2\gamma}
			|\hat{u}|^{2} \, d \xi_{3} \, d \xi_{2}\, d    \xi_{1} \, d \tau 
			\\&\indeq\indeq 
			\les 
			\epsilon^{-2\theta} 
			\Vert u \Vert^{2}_{L^{2}((-\infty,\infty)\scl H^{\gamma\theta+1/2}(\Omegaf))}
			+  
			\epsilon^{1/\gamma} 
			\Vert u \Vert^{2}_{H^{\theta+1/2\gamma}((-\infty,\infty)\scl L^{2}(\Omegaf))} 
			,
		\end{split}
		\llabel{EQ331}
	\end{align}
	for all $\epsilon\in(0,1]$.
	Using \eqref{EQ114}, we get
	\begin{align}
		\begin{split}
			& \Vert u \Vert^{2}_{H^{\theta}((-\infty,\infty)\scl  L^{2}(\Gammac))}
			\les 
			\epsilon^{-2\theta} 
			\Vert u \Vert^{2}_{L^{2}((-\infty,\infty)\scl H^{r}(\Omegaf))}
			+  
			\epsilon^{2\theta/(2r-1)}
			\Vert u \Vert^{2}_{H^{2\theta r/(2r-1)}((-\infty,\infty)\scl L^{2}(\Omegaf))} 
			,
		\end{split}
		\label{EQ319}
	\end{align}
	for all $\epsilon\in(0,1]$.
	Finally, note that $\lambda=2r\theta/(2r-1)$ satisfies \eqref{EQ392}--\eqref{EQ118} under the condition $r>1/2$.
\end{proof}

Optimizing $\epsilon\in(0,1]$ in \eqref{EQ319} by using
\begin{equation}
	\epsilon
	=
	\left(
	\frac{
		\Vert u \Vert_{L^{2}((-\infty,\infty)\scl H^{r}(\Omegaf))}
	}{
		\Vert u \Vert_{L^{2}((-\infty,\infty)\scl H^{r}(\Omegaf))} 
		+    
		\Vert u \Vert_{H^{2\theta r/(2r-1)}((-\infty,\infty)\scl L^{2}(\Omegaf))} 
	}
	\right)^{(2r-1)/2r \theta}
	,
	\llabel{EQ388}
\end{equation}
we obtain a trace inequality
\begin{align}
	\begin{split}
		\Vert u\Vert_{H^{\theta}((-\infty,\infty)\scl  L^{2}(\Gammac))} 
		\les
		\Vert u \Vert_{L^{2}((-\infty,\infty)\scl H^{r}(\Omegaf))}^{1/(2r+1)}
		\Vert u \Vert_{H^{2\theta r/(2r-1)}((-\infty,\infty)\scl L^{2}(\Omegaf))}^{2r/(2r+1)}
		+
		\Vert u \Vert_{L^{2}((-\infty,\infty)\scl H^{r}(\Omegaf))}
		,
	\end{split}
	\label{EQ137}
\end{align}
which is a more explicit version of~\eqref{EQ65a}.
Note that from \eqref{EQ137}, one may obtain an inequality on the interval
$(0,T)$ with a $T$-dependent constant.

The second lemma provides a space-time interpolation inequality which is needed in several places in Sections~\ref{sec05} and \ref{sec06} below.
\cole
\begin{Lemma}
	\label{L12}
	Let $\alpha,\beta>0$.
	If 
	$u\in H^{\alpha}((-\infty, \infty)\scl L^2(\Omegaf))
	\cap
	L^2((-\infty, \infty), H^{\beta}(\Omegaf))$, then we have that
	$
	u \in H^{\theta} ((-\infty, \infty), H^{\lambda}(\Omegaf))
	$ for all $\theta\in(0,\alpha)$ and $\lambda\in(0,\beta)$ such that
	\begin{equation}
		\frac{\theta}{\alpha} + \frac{\lambda}{\beta} \leq 1
		.
		\llabel{EQ326}
	\end{equation}
	In addition, for all $\epsilon \in (0,1]$,
	we have the inequality
	\begin{align}
		\begin{split}
			\Vert u\Vert_{H^{\theta} ((-\infty, \infty), H^{\lambda}(\Omegaf))}
			\leq
			\epsilon
			\Vert u \Vert_{H^{\alpha}((-\infty, \infty)\scl L^2(\Omegaf))}
			+ 
			C \epsilon^{-\frac{\theta \beta}{\alpha\lambda}}
			\Vert u \Vert_{L^2((-\infty, \infty), H^{\beta}(\Omegaf))} 
			,
		\end{split}
		\llabel{EQ338}
	\end{align}
	where $C>0$ is a constant.
\end{Lemma}
\colb

This statement immediately implies the following one regarding the
same type of interpolation inequality on a finite time interval.

\cole
\begin{Corollary}
	\label{C02}
	Let $\alpha,\beta>0$ and~$T>0$.
	If 
	$u\in H^{\alpha}( (0,T)\scl L^2(\Omegaf))
	\cap
	L^2((0,T), H^{\beta}(\Omegaf))$, then we have that
	$
	u \in H^{\theta} ((0,T), H^{\lambda}(\Omegaf))
	$ for all $\theta\in(0,\alpha)$ and $\lambda\in(0,\beta)$ such that
	\begin{equation}
		\frac{\theta}{\alpha} + \frac{\lambda}{\beta} \leq 1
		.
		\llabel{EQ66}
	\end{equation}
	In addition, for all $\epsilon \in (0,1]$,
	we have the inequality
	\begin{align}
		\begin{split}
			\Vert u\Vert_{H^{\theta} ((0,T), H^{\lambda}(\Omegaf))}
			\leq
			\epsilon
			\Vert u \Vert_{H^{\alpha}((0,T)\scl L^2(\Omegaf))}
			+ 
			C \epsilon^{-\frac{\theta \beta}{\alpha\lambda}}
			\Vert u \Vert_{L^2((0,T), H^{\beta}(\Omegaf))} 
			,
		\end{split}
		\llabel{EQ70}
	\end{align}
	where $C>0$ is a constant depending on $\Omegaf$ and~$T$.
\end{Corollary}
\colb

As above, Corollary~\ref{C02} follows by employing a Sobolev extension
operator in the $t$ variable.

\begin{proof}[Proof of Lemma~\ref{L12}]
	Using a partition of unity, straightening of the boundary, and a Sobolev extension, it is sufficient to prove
	the inequality in the case $\Omegaf={\mathbb R}^3$
	and $u\in C_{0}^{\infty}({\mathbb R}\times{\mathbb R}^3)$. 
	Then, using Parseval's identity and the definition of the
	Sobolev norms,
	we only need to prove
	\begin{align}
		\begin{split}
			(1+|\tau|^{2\theta})(1 + |\xi|^{2 \lambda})
			\leq
			\epsilon (1+|\tau|^{2\alpha})
			+ 
			C \epsilon^{-\frac{\theta\beta}{\alpha\lambda}} 
			(1+|\xi|^{2\beta})
			,
		\end{split}
		\label{EQ328}
	\end{align}
	for $\tau \in {\mathbb R}$ and $\xi\in{\mathbb R}^{n}$, where $\epsilon \in (0,1]$.
	Finally, \eqref{EQ328} follows from the Young's inequality.
\end{proof}

In the last part of this section, we address the regularity for the
wave equation. 
We first recall the hidden regularity result for the wave equation 
\begin{align}
	&
	w_{tt}   - \Delta{w}= 0
	\inon{in~ $ [0,T]\times\Omegae $},
	\label{EQ14}
	\\&
	w = \psi    
	\inon{on~$  [0,T]\times \Gammac $}
	,
	\\&
	w \text{~~periodic~in~the~$y_1$~and~$y_2$~directions}
	,
\end{align}
and the initial data
\begin{equation}
	(w,w_t)(0,\cdot)= (w_0,w_1)   
	\label{EQ358}
\end{equation}
(cf.~\cite{LLT}).

\cole
\begin{Lemma}[\cite{LLT}]
	\label{L03}
	Assume
	that
	$(w_{0},w_{1}) \in H^{\beta}( \Omegae)\times H^{\beta -1}( \Omegae)$,
	where $\beta \geq1$,
	and 
	$$\psi \in C([0,T]\scl  H^{\beta-1/2}(\Gammac))    \cap H^{\beta, \beta}((0,T)\times \Gammac),$$
	with the compatibility conditions $\psi|_{t=0} = w_0|_{\Gammac}$ and $\partial_t \psi|_{t=0} = w_1 |_{\Gammac}$.
	Then there exists a solution $(w, w_{t}) \in C([0,T]\scl   H^{\beta}( \Omegae) \times  H^{\beta-1}( \Omegae))$
	of \eqref{EQ14}--\eqref{EQ358}, which satisfies the estimate
	\begin{align}
		\begin{split}
			&\Vert w \Vert_{C([0,T]\scl   H^{\beta}( \Omegae))} 
			+ 
			\Vert w_{t} \Vert_{C([0,T]\scl   H^{\beta-1}( \Omegae))}
			+ 
			\left\Vert 
			\frac{\partial w}{\partial \nu}
			\right\Vert_{H^{\beta-1, \beta-1}((0,T) \times  \Gammac)} 
			\\&\indeq
			\les
			\Vert w_{0} \Vert_{H^{\beta}( \Omegae)} + \Vert w_{1} \Vert_{H^{\beta-1}( \Omegae)} 
			+
			\Vert \psi \Vert_{H^{\beta, \beta}((0,T) \times \Gammac)}
			,
		\end{split}
		\llabel{EQ339}
	\end{align}
	where the implicit constant depends on $\Omegae$ and~$T$.
\end{Lemma}
\colb

In the final lemma of this section, we recall an essential 
trace regularity result for the wave equation from~\cite{RV}.

\cole
\begin{Lemma}[\cite{RV}]
	\label{L13}
Assume that $(w_0, w_1) \in H^{\beta} (\Omegae) \times H^{\beta+1} (\Omegae)$, where $0<\beta <5/2$, and
	\begin{align}
		\psi \in L^2((0,T), H^{\beta +2} (\Gammac))
		\cap
		H^{\beta/2+1}((0,T), H^{\beta/2+1}(\Gammac)),
   \llabel{EQ117}
	\end{align}
with the compatibility condition $\partial_t \psi|_{t=0} = w_1 |_{\Gammac}$.
Then there exists a solution $w$ of \eqref{EQ14}--\eqref{EQ358} such that
	\begin{align}
		\begin{split}
			\left\Vert \frac{\partial w}{\partial \nu} \right\Vert_{L^2((0,T), H^{\beta+1} (\Gammac))}
			&
			\les
			\Vert w_0 \Vert_{H^{\beta+2} (\Omegae)}
			+
			\Vert w_1 \Vert_{H^{\beta+1} (\Omegae)}
			+
			\Vert \psi \Vert_{L^2((0,T), H^{\beta+2}(\Gammac))}
			\\&\indeq
			+
			\Vert \psi \Vert_{H^{\beta/2 +1}((0,T), H^{\beta/2 +1}(\Gammac))}
			,
			\llabel{EQ03}
		\end{split}
	\end{align}
	where the implicit constant depends on $\Omegae$ and~$T$.
\end{Lemma}
\colb

\startnewsection{The nonhomogeneous parabolic problem}{sec04} 
In this section, we consider the parabolic problem
\begin{align}
	u_t - \lambda \RR \dive (\nabla u + (\nabla u)^T)
	-
	\mu \RR \nabla \dive u 
	&
	= f
	\inon{in~$[0,1] \times\Omegaf$}
	,
	\label{EQ48}
\end{align}
with the nonhomogeneous boundary conditions and the initial data
\begin{align}
	&
	\lambda (\partial_k u_j + \partial_j u_k)
	\nu^k
	+
	\mu \partial_k u_k \nu^j
	=
	h_j
	\inon{on~$[0,1]\times\Gammac$}
	,
	\label{EQ44}
	\\
	&
	u
	=0
	\inon{on~$[0,1]\times \Gammaf$}
	,
	\label{EQ45}
	\\&
	u~~\text{periodic in the $y_1$~and~$y_2$ directions},
	\label{EQ320} 
	\\&
	u(0) 
	= 
	u_0
	\inon{in~$\Omegaf$}
        ,
	\label{EQ46}
\end{align}
for $j=1,2,3$.
To state the maximal regularity for \eqref{EQ48}--\eqref{EQ46}, we consider the homogeneous version
when \eqref{EQ44}--\eqref{EQ46} is replaced by
  \begin{align}
	&
	\lambda (\partial_k u_j  + \partial_j u_k) \nu^k 
	+ 
	\mu  \partial_k u_k \nu^j 
	= 
	0
	\inon{on~$[0,1]\times \Gammac$}
	,
	\label{EQ02}
	\\
	&
	u
	=
	0
	\inon{on~$[0, 1]\times \Gammaf$}
	,
	\label{EQ17}
	\\&
	u~~\text{periodic in the $y_1$~and~$y_2$ directions},
	\label{EQ317} 
	\\&
	u(0)
	= 
	0
	\inon{in~$\Omegaf$}
        ,
	\label{EQ18}
  \end{align}
for   $j=1,2,3$.

\cole
\begin{Lemma}
\label{L01}
Assume
that
$f\in K^{2}((0,1)\times\Omegaf)$ with
$f(0, \cdot) = 0$ on $\Omegaf$
and
\begin{align}
	(\RR, \RR^{-1})
	\in
	(	L^\infty ((0,1), H^2 ( \Omegaf))
	 \cap 
	 H^1  ((0,1), L^\infty (\Omegaf) ))^2
	.
   \label{EQ05}
\end{align}
Then the parabolic problem \eqref{EQ48} with the boundary conditions and the initial data \eqref{EQ02}--\eqref{EQ18} admits a solution 
$u$ satisfying 
\begin{align}
	\begin{split}
	\Vert u \Vert_{K^2 ( (0, 1)\times \Omegaf)}
	\les	
	\Vert f\Vert_{K^0  ((0,1) \times \Omegaf)}
	\label{EQ50}
	\end{split}
\end{align}
and
\begin{align}
	\begin{split}
		\Vert u \Vert_{K^4 ( (0, 1)\times \Omegaf)}
	\les	
	\Vert f\Vert_{K^2  ((0, 1) \times \Omegaf)}
        ,
	\label{EQ59}
	\end{split}
\end{align}
where the implicit constants depend on the norms of $\RR$ and $\RR^{-1}$ in~\eqref{EQ05}.
\end{Lemma}
\colb

\begin{proof}
Analogously to \cite[Theorem~3.2]{LM}, the parabolic problem \eqref{EQ48} admits a solution $u\in K^2((0,1) \times\Omegaf)$ if $f \in K^0((0,1) \times \Omegaf)$ and $u \in K^4((0,1) \times \Omegaf)$ if $f \in K^2((0, 1) \times \Omegaf)$. 
Below, the norm of dependence on time and space are understood as $(0,1)$ and $\Omegaf$, unless stated otherwise.
In the reminder of the proof we shall prove the regularity.
Taking the $L^2$-inner product of \eqref{EQ48} with $u$, we arrive at
\begin{align}
	\frac{1}{2} \frac{d}{dt}
	\int_{\Omegaf} |u|^2
	-
	\lambda \int_{\Omegaf} \RR u_j \partial_k (\partial_k u_j + \partial_j u_k)  
	-
	\mu \int_{\Omegaf}
	\RR u_j \partial_j \partial_k u_k  
	=
	\int_{\Omegaf} f u
	.
	\label{EQ25}
\end{align}
For the second and third terms on the left side of \eqref{EQ25}, we integrate by parts with respect to $\partial_k$ and $\partial_j$ respectively to get
\begin{align}
	\begin{split}
		-\lambda \int_{\Omegaf} \RR u_j \partial_k (\partial_k u_j + \partial_j u_k) 
		&
		=
		\lambda \int_{\Gammac} \RR u_j (\partial_k u_j + \partial_j u_k) \nu^k 
		+
		\lambda \int_{\Omegaf} \RR \partial_k u_j (\partial_k u_j + \partial_{j} u_k )
		\\
		&\indeq
		+
		\lambda \int_{\Omegaf} u_j \partial_k \RR  (\partial_k u_j + \partial_{j} u_k )
		\label{EQ28}
	\end{split}
\end{align}
and
\begin{align}
	\begin{split}
		-
		\mu \int_{\Omegaf}
		\RR u_j \partial_j \partial_k u_k 
		&
		=
		\mu \int_{\Gammac}
		\RR		u_j \partial_k u_k \nu^j
		+
		\mu
		\int_{\Omegaf} \RR \partial_j u_j \partial_k u_k  
		+
		\mu
		\int_{\Omegaf} u_j \partial_j \RR \partial_k u_k 
		.
		\label{EQ29}
	\end{split}
\end{align}
Inserting \eqref{EQ28}--\eqref{EQ29} into \eqref{EQ25} and appealing to \eqref{EQ02}--\eqref{EQ17}, we get
\begin{align}
	\begin{split}
		&
		\frac{1}{2} \frac{d}{dt} 
		\int_{\Omegaf} |u |^2	
		+
		\lambda\int_{\Omegaf} \RR \partial_k u_j (\partial_k u_j+\partial_j u_k) 
		+
		\mu \int_{\Omegaf} \RR  \partial_j u_j \partial_k u_k
		\\
		&\indeq
		=	
		\int_{\Omegaf} f u
		-
		\lambda \int_{\Omegaf} 
		u_j
		\partial_k \RR (\partial_k u_j + \partial_{j} u_k ) 
		-
		\mu
		\int_{\Omegaf} u_j \partial_j \RR \partial_k u_k 
		\\
		&\indeq
		\les
		 \Vert f\Vert_{L^2}^2
		+
		 \Vert u \Vert_{L^2}^2
		+
		\Vert u\Vert_{L^4} 
				\Vert \nabla \RR \Vert_{L^4} 
		\Vert \nabla u\Vert_{L^2}
		,
		\label{EQ30}
	\end{split}
\end{align}
where the last inequality follows from H\"older's and Young's inequalities.
Note that for any $v \in H^1(\Omegaf)$, using the Sobolev and Young's inequalities, we have
\begin{align}
	\Vert v\Vert_{L^4 (\Omegaf)}
	\les
	\Vert v\Vert_{H^1 (\Omegaf)}^{3/4}
	\Vert v\Vert_{L^2 (\Omegaf)}^{1/4}
	\les
	\epsilon \Vert v\Vert_{H^1 (\Omegaf)}
	+
	C_{\epsilon} \Vert v\Vert_{L^2 (\Omegaf)}
	,
	\label{EQ330}
\end{align}
for any $\epsilon \in (0,1]$, where $C_\epsilon>0$ denotes a constant depending on~$\epsilon$.
We integrate \eqref{EQ30} in time from $0$ to $t$ and use
  \begin{equation}
   (\partial_k u_j+\partial_j u_k) \partial_k u_j
   = \frac12  \sum_{j,k=1}^{3}  (\partial_k u_j+\partial_j u_k)^2
   ,
   \label{EQ06}
  \end{equation}
obtaining
\begin{align}
	\begin{split}
		&
		\Vert  u(t) \Vert_{L^2}^2 
		+ 
		\sum_{j,k=1}^3 \int_0^t \int_{\Omegaf} \RR ( \partial_k u_j + \partial_j u_k )^2
		+
		\int_0^t \int_{\Omegaf} \RR |\partial_k u_k|^2
		\\
		 &\indeq
		 \les
		 \Vert f \Vert_{L^2_t L^2_x}^2
		 +
		 \int_0^t \Vert u \Vert_{L^2}^2
		 +
		 \epsilon \int_0^t  \Vert u\Vert_{H^1}^2
		 +
		  C_{\epsilon}  \int_0^t 
		 \Vert u\Vert_{L^2}
		 \Vert u\Vert_{H^1}
		 \\&\indeq
		 \les
		\Vert f \Vert_{L^2_t L^2_x}^2
		 +
		 (\epsilon + C_\epsilon \bar{\epsilon} ) \int_0^t \Vert u \Vert_{H^1}^2
		 +
		 C_{\epsilon,\bar{\epsilon}}  \int_0^t 
		 \Vert u\Vert_{L^2}^2
		,
		\label{EQ32}
	\end{split}
\end{align}
for any $\epsilon, \bar{\epsilon} \in (0,1]$,
where we used \eqref{EQ330} and the Young's inequality.
For the second term on the left, we use
Korn's inequality, which reads
\begin{align}
	\begin{split}
	\int_0^t
	\Vert u\Vert_{H^1}^2	
	\les
	\sum_{j,k=1}^3
	\int_0^t \int_{\Omegaf}
	\RR (\partial_k u_j + \partial_j u_k)^2	
	+
	\int_0^t
	\Vert u\Vert_{L^2}^2
	.
	\label{EQ33}
	\end{split}
\end{align}
From \eqref{EQ32}--\eqref{EQ33} it follows that
\begin{align}
	\begin{split}
		\Vert u(t) \Vert_{L^2}^2 
		+
		\Vert u\Vert_{L^2_t H^1_x}^2
		&
		\les
		\Vert f \Vert_{L^2_t L^2_x}^2
		+
		\int_0^t \Vert u \Vert_{L^2}^2
		,
		\label{EQ35}
	\end{split}
\end{align}	
by choosing suitable $\epsilon, \bar{\epsilon}>0$.
By  Gronwall's inequality, we obtain
\begin{align}
	\begin{split}
	\Vert u(t)\Vert_{L^2}^2
	&
	\les
	\Vert f\Vert_{L^2_t L^2_x}^2
        ,
	\label{EQ34}
	\end{split}
\end{align} 
where we used $e^{C t}\les 1$ for $t\le 1$,
and then, after using \eqref{EQ34} in \eqref{EQ35}, we arrive at
\begin{align}
	\begin{split}
	\Vert u\Vert_{L^2_t H^1_x}^2
	&
	\les
	\Vert f\Vert_{L^2_t L^2_x}^2 
	.
	\label{EQ36}
	\end{split}
\end{align}

Next, we take the $L^2$-inner product of \eqref{EQ48} with $u_t$, obtaining
\begin{align}
	\int_{\Omegaf} |u_t|^2
	-
	\lambda \int_{\Omegaf} \RR u_{tj} \partial_k (\partial_k u_j + \partial_j u_k) 
	-
	\mu \int_{\Omegaf}
	\RR u_{tj} \partial_j \partial_k u_k 
	=
	\int_{\Omegaf} f u_t
	.
	\label{EQ07}
\end{align}
Then, proceeding as in \eqref{EQ28}--\eqref{EQ29}, we get
\begin{align}
	\begin{split}
	-\lambda \int_{\Omegaf} \RR u_{tj} \partial_k
	 (\partial_k u_j + \partial_j u_k) 
	&
	=
	\lambda \int_{\Gammac} \RR u_{tj}
	(\partial_k u_j + \partial_j u_k) \nu^k 
	+
	\lambda \int_{\Omegaf} \RR \partial_k u_{tj}
	(\partial_k u_j + \partial_{j} u_k )
	\\
	&\indeq
	+
	\lambda \int_{\Omegaf} u_{tj} \partial_k \RR 
	(\partial_k u_j + \partial_{j} u_k )
	\label{EQ08}
\end{split}
\end{align}
and
\begin{align}
	\begin{split}
	-
	\mu \int_{\Omegaf}
	\RR u_{tj} \partial_j \partial_k u_k  
	&
	=
	\mu \int_{\Gammac}
	\RR u_{tj} \partial_k u_k \nu^j
	 +
	 \mu
	 \int_{\Omegaf} \RR \partial_j u_{tj} \partial_k u_k 
	 +
	 \mu
	 \int_{\Omegaf} u_{tj} \partial_j \RR \partial_k u_k 
	 .
	 \label{EQ09}
	\end{split}
\end{align}
Inserting \eqref{EQ08}--\eqref{EQ09} into \eqref{EQ07} appealing to
\eqref{EQ02}--\eqref{EQ17},
and using
  \begin{align}
  \begin{split}
   \frac12
   \frac{d}{dt}
   \int_{\Omegaf} \RR 	\partial_k u_j
	(\partial_k u_j +\partial_{j} u_k )
	= 
    \frac12
   \int_{\Omegaf} \RR_t 	\partial_k u_j
	(\partial_k u_j +\partial_{j} u_k )
    +
   \int_{\Omegaf} \RR 	\partial_k u_{tj}
	(\partial_k u_j +\partial_{j} u_k )
   ,
  \end{split}
   \label{EQ41}
  \end{align}
we arrive at
\begin{align}
	\begin{split}
		&
	\int_{\Omegaf} |u_t|^2	
	+
	\frac{\lambda}{2} \frac{d}{dt} \int_{\Omegaf} \RR 	\partial_k u_j
	(\partial_k u_j +\partial_{j} u_k )
	+
	\frac{\mu}{2} \frac{d}{dt}
	 \int_{\Omegaf} \RR \partial_j u_j \partial_k u_k 
	\\
	&\indeq
	=	
	\int_{\Omegaf} f u_t
	+
	\lambda \int_{\Omegaf} u_{tj} \partial_k \RR (\partial_k u_j + \partial_{j} u_k )
	+
	\mu
	\int_{\Omegaf} u_{tj} \partial_j \RR \partial_k u_k  
	\\
	&\indeq\indeq
	+
	\frac{\lambda}{2}	\int_{\Omegaf} \RR_t 	\partial_k u_j
	(\partial_k u_j + \partial_j u_k)
	+
	\frac{\mu}{2} \int_{\Omegaf} 
	\RR_t \partial_j u_j \partial_k u_k 
	\\&\indeq
	\les
	C_\epsilon \Vert f\Vert_{L^2}^2
	+
	\epsilon \Vert u_t \Vert_{L^2}^2
	+
	\Vert \nabla \RR \Vert_{L^4}
	 \Vert \nabla u \Vert_{L^4} \Vert u_t \Vert_{L^2}
	+
	\Vert \RR_t \Vert_{L^\infty} \Vert \nabla u\Vert_{L^2}^2
	,
	\llabel{EQ20}
	\end{split}
\end{align}
for any $\epsilon \in (0,1]$, where we used H\"older's and Young's inequalities.
Integrating in time from $0$ to $t$ and using the Young, Sobolev, and Korn's inequalities with \eqref{EQ330}--\eqref{EQ06}, we get
\begin{align}
	\begin{split}
		&
	\Vert u_t \Vert_{L^2_t L^2_x }^2 
	+ 
	\Vert u(t) \Vert_{H^1}^2
	\\&\indeq
	\les
	 C_\epsilon \Vert f \Vert_{L^2_t L^2_x}^2
	 +
	 \epsilon \Vert u_t\Vert_{L^2_t L^2_x}^2
	 +
	 \Vert u(t) \Vert_{L^2}^2
	 +
	 \int_0^t  (\bar{\epsilon} \Vert u\Vert_{H^2} + C_{\bar{\epsilon}} \Vert u\Vert_{H^1} )
	 \Vert u_t\Vert_{L^2}
	 +
	 \int_0^t \Vert \RR_t \Vert_{L^\infty} \Vert u\Vert_{H^1}^2
	 \\
	 &\indeq
	 \les
	C_\epsilon \Vert f \Vert_{L^2_t L^2_x}^2
	+
	(\epsilon + \bar{\epsilon} + \tilde{\epsilon} C_{\bar{\epsilon}}) \Vert u_t \Vert_{L^2_t L^2_x}^2
	+
	\Vert u(t) \Vert_{L^2}^2
	+
	\bar{\epsilon} \Vert u\Vert_{L^2_t H^2_x}^2
	+
	C_{\bar{\epsilon}, \tilde{\epsilon}} \Vert u\Vert_{L^2_t H^1_x}^2
	\\&\indeq\indeq 
	+
	\int_0^t \Vert \RR_t \Vert_{L^\infty} \Vert u\Vert_{H^1}^2
	 ,
	 \label{EQ24}
	\end{split}
\end{align}
for any $\epsilon, \bar{\epsilon}, \tilde{\epsilon} \in (0,1]$, where we used $\Vert u (0) \Vert_{H^1} = 0$ in the last inequality by~\eqref{EQ18}.
For the space regularity, note that $u$ is the solution of the elliptic problem
\begin{align}
	\begin{split}
		&
		-\lambda \dive(\nabla u + (\nabla u)^T) 
		- 
		\mu  \nabla \dive u 
		=
		-\frac{u_t}{\RR} 
		+ 
		\frac{f}{\RR}
		\inon{in~$[0,T]\times\Omegaf$}
		,
		\label{EQ67}
	\end{split}
\end{align}
with the boundary conditions
\begin{align}
	&
	\lambda (\partial_k u_j  + \partial_j u_k)\nu^k 
	+ 
	\mu  \partial_k u_k \nu^j 
	= 
	0
	\inon{on~$[0,1] \times \Gammac$}
	,
	\label{EQ68}
	\\
	&
	u
	=
	0
	\inon{on~$[0,1] \times\Gammaf$}
	,
	\label{EQ69}
\end{align}
for $j=1,2,3$.
From the elliptic regularity for \eqref{EQ67}--\eqref{EQ69} it follows that
\begin{align}
	\Vert u\Vert_{H^2}
	\les
	\Vert \RR^{-1} u_t \Vert_{L^2}
	+
	\Vert \RR^{-1} f\Vert_{L^2}
	\les
	\Vert u_t \Vert_{L^2}
	+
	\Vert f\Vert_{L^2}
	,
	\label{EQ160}
\end{align}
from where
\begin{align}
	\begin{split}
		\Vert u\Vert_{L^2_t H^2_x}
		&
		\les
		\Vert u_t \Vert_{L^2_t L^2_x}
		+
		\Vert f \Vert_{L^2_t L^2_x}
		.
		\label{EQ158}
	\end{split}
\end{align}
Combining \eqref{EQ34}--\eqref{EQ36}, \eqref{EQ24}, and \eqref{EQ158}, we obtain
\begin{align}
	\begin{split}
	\Vert u_t \Vert_{L^2_t L^2_x }^2 
	+
	\Vert u(t) \Vert_{H^1}^2
	&
	\les
	\Vert f \Vert_{L^2_t L^2_x}^2
	+
	\int_0^t \Vert \RR_t \Vert_{L^\infty} 
	\Vert u \Vert_{H^1}^2
	,
	\label{EQ56}
	\end{split}
\end{align}
by taking suitable $\epsilon, \bar{\epsilon}, \tilde{\epsilon}>0$.
Using Gronwall's inequality, we arrive at
\begin{align}
	\begin{split}
	\Vert u(t)\Vert_{H^1}^2
	&
	\leq
	C
	\Vert f\Vert_{L^2_t L^2_x}^2
	\exp\left(C \int_\tau^t 
	\Vert \RR_t(\tau) \Vert_{L^\infty} 
	\,d\tau \right) 
	\leq
	C
	\Vert f\Vert_{L^2_t L^2_x}^2
	,
	\llabel{EQ57}
	\end{split}
\end{align}
and thus \eqref{EQ56} implies
\begin{align}
	\begin{split}
	\Vert u\Vert_{H^1_t L^2_x}^2
	\les
	\Vert f\Vert_{L^2_t L^2_x}^2
	,
	\label{EQ39}
	\end{split}
\end{align}
where we used $e^{C t}\les 1$ for $t\le 1$.
From \eqref{EQ158} and \eqref{EQ39} it follows that
\begin{align}
	\begin{split}
	\Vert u\Vert_{K^2} 
	&
	\les
	\Vert u\Vert_{L^2_t H^2_x}
	+
	\Vert u\Vert_{H^1_t L^2_x}
	\les
	\Vert f\Vert_{L^2_t L^2_x}
	,
	\end{split}
   \llabel{EQ120}
\end{align}
completing the proof of~\eqref{EQ50}.

Differentiating \eqref{EQ48} in time and taking the $L^2$-inner product with $u_t$, we arrive at
\begin{align}
	\begin{split}
		&
	\frac{1}{2} \frac{d}{dt} \int_{\Omegaf} |u_t |^2
	-
	\lambda \int_{\Omegaf} \RR 	u_{tj}
	\partial_k(\partial_k u_{tj} + \partial_j u_{tk})	
	-
	\mu \int_{\Omegaf} \RR  u_{tj} \partial_j \partial_k u_{tk} 
	\\
	&\indeq
	=
	\int_{\Omegaf} f_t u_t
	+
	\lambda \int_{\Omegaf} \RR_t u_{tj} \partial_k
	(\partial_k u_j + \partial_j u_k) 
	+
	\mu \int_{\Omegaf} \RR_t u_{tj} \partial_{j} \partial_k u_k 
	.
	\label{EQ60}
	\end{split}
\end{align}
We proceed as in \eqref{EQ28}--\eqref{EQ29} to obtain
\begin{align}
	\begin{split}
	-
	\lambda \int_{\Omegaf} \RR u_{tj}	
	\partial_k(\partial_k u_{tj} + \partial_j u_{tk}) 
	&
	=
	\lambda \int_{\Gammac} \RR u_{tj}	
	(\partial_k u_{tj} + \partial_j u_{tk})	\nu^k 
	+
	\lambda\int_{\Omegaf} u_{tj} \partial_k \RR 
	(\partial_k u_{tj} + \partial_j u_{tk})	
	\\&\indeq
	+
	\lambda	\int_{\Omegaf} \RR \partial_k u_{tj}	
	(\partial_k u_{tj} + \partial_j u_{tk})	
	\label{EQ61}
	\end{split}
\end{align}
and
\begin{align}
	\begin{split}
	-
	\mu \int_{\Omegaf} \RR u_{tj} \partial_j \partial_k u_{tk}  
	&
	=
	\mu \int_{\Gammac} \RR u_{tj} \partial_k u_{tk} \nu^j  
	+
	\mu \int_{\Omegaf} u_{tj} \partial_j  \RR \partial_k u_{tk}  
	+
	\mu \int_{\Omegaf}  \RR \partial_j u_{tj}  \partial_k u_{tk}  
	.
	\label{EQ62}
	\end{split}
\end{align}
Inserting \eqref{EQ61}--\eqref{EQ62} into \eqref{EQ60}, we get
\begin{align}
	\begin{split}
	&
	\frac{1}{2} \frac{d}{dt} \int_{\Omegaf} |u_t|^2
	+
	\lambda \int_{\Omegaf} \RR \partial_k u_{tj}
	(\partial_k u_{tj} + \partial_j u_{tk})	
	+
	\mu \int_{\Omegaf} \RR \partial_j u_{tj} \partial_k u_{tk} 
	\\
	&\indeq
	\les
	\Vert f_t\Vert_{L^2}^2 
	+
	\Vert u_{t}\Vert_{L^2}^2
	+
	\Vert \RR_t \Vert_{L^\infty} \Vert u_t \Vert_{L^2} \Vert u\Vert_{H^2}
	+
	\Vert u_t\Vert_{L^4}
		\Vert \nabla \RR\Vert_{L^4}
	 \Vert u_t\Vert_{H^1}
	,
	\end{split}
   \llabel{EQ122}
\end{align}
where we used Young's, H\"older's, and Sobolev inequalities.
Integrating in time from $0$ to $t$ and
using the Young's and Korn's inequalities and \eqref{EQ330}--\eqref{EQ06}, we obtain
\begin{align}
	\begin{split}
		\Vert u_t(t) \Vert_{L^2}^2
		+
		\Vert u_t \Vert_{L^2_t H^1_x}^2
		&
		\les
		\Vert f\Vert_{H_t^1 L_x^2}^2 
		+
		\int_0^t \Vert \RR_t \Vert_{L^\infty}
		\Vert u \Vert_{H^2}^2	
		+
		\int_0^t \Vert \RR_t \Vert_{L^\infty}
		\Vert u_t \Vert_{L^2}^2	
		+
		(\epsilon + \bar{\epsilon} C_\epsilon ) \Vert u_t\Vert_{L^2_t H^1_x}^2
						\\
		&\indeq
		+
		C_{\epsilon, \bar{\epsilon}} \Vert u_t\Vert_{L^2_t L^2_x}^2
		,
		\llabel{EQ103}
	\end{split}
\end{align}
for any $\epsilon, \bar{\epsilon} \in (0,1]$, since $\Vert u_t (0)\Vert_{L^2} = \Vert f(0)\Vert_{L^2} =0$.
From \eqref{EQ160} and \eqref{EQ39} it follows that
\begin{align}
	\begin{split}
		\Vert u_t(t) \Vert_{L^2}^2
		+
		\Vert u_t \Vert_{L^2_t H^1_x}^2
		&
		\les
		\Vert f\Vert_{H_t^1 L_x^2}^2 
		+
		\Vert u_{t}\Vert_{L^2_t L^2_x}^2
		+
		\int_0^t \Vert \RR_t \Vert_{L^\infty}
		\Vert u_t \Vert_{L^2}^2	
		+
		\int_0^t \Vert \RR_t \Vert_{L^\infty}
		\Vert f \Vert_{L^2}^2	
		\\&
		\les
		\Vert f\Vert_{H_t^1 L_x^2}^2 
		+
		\Vert u_{t}\Vert_{L^2_t L^2_x}^2
		+
		\int_0^t \Vert \RR_t \Vert_{L^\infty}
		\Vert u_t \Vert_{L^2}^2	
                ,
		\label{EQ162}
	\end{split}
\end{align}
by taking appropriate $\epsilon, \bar{\epsilon} >0$, where we also used 
\begin{align}
	\Vert f\Vert_{L^\infty_t L^2_x} \les \Vert f\Vert_{H^1_t L^2_x}
	\label{EQ950}
\end{align}
in the last inequality. 
Appealing to Gronwall's inequality, \eqref{EQ162} implies
\begin{align}
	\begin{split}
		\Vert u_t(t) \Vert_{L^2}^2
		\les
		\Vert f\Vert_{K^2}^2	
		,	
		\label{EQ104}
	\end{split}
\end{align}
and then, after using \eqref{EQ104} in \eqref{EQ162}, we arrive at
\begin{align}
	\begin{split}
	\Vert u_t\Vert_{L^2_t H^1_x}^2
	\les
	\Vert f\Vert_{K^2}^2	
	.	
   \llabel{EQ208}
	\end{split}
\end{align}

Differentiating \eqref{EQ48} in time and taking the $L^2$-inner product with $u_{tt}$, we obtain
\begin{align}
	\begin{split}
		&
	\int_{\Omegaf} |u_{tt}|^2
	+
	\frac{\lambda}{2} \frac{d}{dt} \int_{\Omegaf} \RR \partial_k u_{tj} (\partial_k u_{tj} + \partial_j u_{tk}) 
	+
	\frac{\mu}{2} \frac{d}{dt} \int_{\Omegaf} \RR \partial_j u_{tj} \partial_k u_{tk} 
	\\
	&
	=
	\int_{\Omegaf} f_t u_{tt} 
	-
	\lambda \int_{\Omegaf} u_{ttj} \partial_k \RR (\partial_k u_{tj} + \partial_j u_{tk}) 
	-
	\mu \int_{\Omegaf} u_{ttj} \partial_j \RR \partial_k u_{tk} 
	\\&\indeq
	+
	\frac{\lambda}{2} \int_{\Omegaf} \RR_t  \partial_k u_{tj} (\partial_k u_{tj} + \partial_j u_{tk}) 
	+
	\frac{\mu}{2} \int_{\Omegaf} \RR_t \partial_j  u_{tj} \partial_k u_{tk} 
	+
	\lambda \int_{\Omegaf} u_{ttj} \RR_t  \partial_k 
	(\partial_k u_j + \partial_j u_k) 
	\\&\indeq
	+
	\mu \int_{\Omegaf} u_{ttj} \RR_t  \partial_{jk} u_k 
	,
	\end{split}
   \llabel{EQ123}
\end{align}
where we integrated by parts in spatial variables. 
We proceed as in \eqref{EQ60}--\eqref{EQ162} to get
\begin{align}
	\begin{split}
	\Vert u_{tt}\Vert_{L^2_t L^2_x}^2
	+
	\Vert u_t (t)\Vert_{H^1}^2
	&
	\les
	C_{\tilde{\epsilon}}
	\Vert f\Vert_{H^1_t L^2_x}^2
	+
	\Vert u_t (t)\Vert_{L^2}^2
	+
	\epsilon \int_0^t
	\Vert u_t\Vert_{H^2} 
	\Vert u_{tt} \Vert_{L^2}
	+
	C_{\tilde{\epsilon}}
	\int_0^t  \Vert \RR_t \Vert_{L^\infty}^2
	\Vert f \Vert_{L^2}^2 
	\\
	&\indeq\indeq
	+ 
	(\bar{\epsilon} C_\epsilon + \tilde{\epsilon} + \tilde{\epsilon})
	\Vert u_{tt} \Vert_{L^2_t L^2_x}^2
	+
	C_{\epsilon, \bar{\epsilon}}  \Vert u_t\Vert_{L^2_t H^1_x}^2
	+
	C_{\tilde{\epsilon}}
	\int_0^t (1+\Vert \RR_t\Vert_{L^\infty}^2  ) \Vert u_t\Vert_{H^1}^2
	 ,
	 \label{EQ63}
	\end{split}
\end{align}
for any $\epsilon, \bar{\epsilon}, \tilde{\epsilon} \in (0,1]$, where we used the Young's, H\"older, Sobolev, and Korn's inequalities. 
Note that $u_t$ is the solution of the elliptic problem
\begin{align}
	\begin{split}
		&
		-\lambda \dive(\nabla u_t + (\nabla u_t)^T) 
		- 
		\mu  \nabla \dive u_t 
		=
		-\RR^{-1} u_{tt}
		+
		\RR^{-2} u_t \RR_t
		+ 
		\RR^{-1} f_t
		-
		\RR^{-2} \RR_t f
		\inon{in~$[0,1]\times\Omegaf$}
		,
	\end{split}
   \llabel{EQ126}
\end{align}
with the boundary conditions
\begin{align}
	&
	\lambda (\partial_k u_{tj } + \partial_j u_{tk})\nu^k 
	+ 
	\mu  \partial_k u_{tk} \nu^j 
	= 
	0
	\inon{in~$[0,1] \times \Gammac$},
	\\
	&
	u_{tj}
	=
	0
	\inon{in~$[0,1] \times \Gammaf$},
   \llabel{EQ129}
\end{align}
for $j=1,2,3$.
The elliptic regularity implies that
\begin{align}
	\begin{split}
	\Vert u_t\Vert_{H^2}
	&
	\les
	\Vert u_{tt}\Vert_{L^2}
	+
	\Vert u_t \RR_t \Vert_{L^2}
	+
	\Vert f_t \Vert_{L^2}
	+
	\Vert \RR_t f \Vert_{L^2}
	\\&
	\les	
	\Vert u_{tt}\Vert_{L^2}
	+
	\Vert u_t \Vert_{L^2} \Vert \RR_t \Vert_{L^\infty}
	+
	\Vert f_t \Vert_{L^2}
	+
	\Vert \RR_t \Vert_{L^\infty}
	\Vert f\Vert_{L^2}
	,
	\label{EQ163}
	\end{split}
\end{align}
where we used H\"older's inequality.
From \eqref{EQ104}--\eqref{EQ163}, we obtain
\begin{align}
	\begin{split}
		\Vert u_{tt}\Vert_{L^2_t L^2_x}^2
		+
		\Vert u_t (t)\Vert_{H^1}^2
		&
		\les
		\Vert f\Vert_{K^2}^2
		+
		\Vert u_t (t)\Vert_{L^2}^2
		+
		\int_0^t
		( 1
		+
		 \Vert \RR_t\Vert_{L^\infty}^2
		 )
		 \Vert u_t\Vert_{H^1}^2
		+
		\int_0^t  \Vert \RR_t \Vert_{L^\infty}^2
		\Vert f \Vert_{L^2}^2 
		\\&
		\les
		\Vert f\Vert_{K^2}^2
		+
		\int_0^t
		( 1
		+
		\Vert \RR_t\Vert_{L^\infty}^2
		)
		\Vert u_t\Vert_{H^1}^2
		,
		\llabel{EQ182}
	\end{split}
\end{align}
by taking $\epsilon, \bar{\epsilon}, \tilde{\epsilon}>0$ sufficiently small, where we used~\eqref{EQ950}.
Appealing to Gronwall's inequality, we arrive at
\begin{align}
	\begin{split}
	\Vert u_t (t)\Vert_{H^1}^2
	&
	\les
	\Vert f\Vert_{K^2}^2		
	 ,
	 \llabel{EQ64}
	\end{split}
\end{align}
whence
\begin{align}
	\begin{split}
	\Vert u_{tt} \Vert_{L^2_t L^2_x}^2
	\les
	\Vert f\Vert_{K^2}^2		
        .
	\label{EQ105}
	\end{split}
\end{align}
From the $H^4$ regularity of the elliptic problem \eqref{EQ67}--\eqref{EQ69} and \eqref{EQ163} it follows that
\begin{align}
	\begin{split}
	\Vert u\Vert_{H^4}
	&
	\les
	\Vert \RR^{-1} u_t \Vert_{H^2}
	+
	\Vert \RR^{-1} f \Vert_{H^2}
	\\&
	\les
	\Vert u_{tt} \Vert_{L^2}
	+
	\Vert \RR_t \Vert_{L^\infty} \Vert u_t\Vert_{L^2}
	+
	\Vert \RR_t \Vert_{L^\infty} \Vert f\Vert_{L^2}
	+
	\Vert f_t\Vert_{L^2}
	+
	\Vert f\Vert_{H^2}
	,
	\label{EQ164}
		\end{split}
\end{align}
since $H^2$ is an algebra.
We combine \eqref{EQ104} and \eqref{EQ105}--\eqref{EQ164} to get
\begin{align}
	\begin{split}
	\Vert u\Vert_{K^4}
	&
	=
	\Vert u\Vert_{L^2_t H^4_x}
	+
	\Vert u\Vert_{H^2_t L^2_x}
	\\
	&
	\les
	\Vert u_{tt}\Vert_{L^2_t L^2_x}
	+
	\Vert \RR_t \Vert_{L^2_t L^\infty_x} 
	\Vert u_t\Vert_{L^\infty_t L^2_x}
	+
	\Vert \RR_t \Vert_{L^2_t L^\infty_x} 
	\Vert f\Vert_{H^1_t L^2_x}
	+
	\Vert f\Vert_{K^2}
	\les
	\Vert f\Vert_{K^2}	
	,
	\end{split}
   \llabel{EQ130}
\end{align}
completing the proof of~\eqref{EQ59}.
\end{proof}

The following lemma provides a maximal regularity for the parabolic system~\eqref{EQ48}--\eqref{EQ46}.

\begin{Lemma}
	\label{L02}
Let $s\in (2, 2+\epsilon_0]$, where $\epsilon_0 \in (0,1/2)$
is arbitrary.
Assume the compatibility conditions
	\begin{align}
		&
		h_j(0)	
		=
		\lambda (\partial_k u_{0j} 
		+ 
		\partial_j u_{0k}) \nu^k
		+
		\mu \partial_k u_{0k} \nu^j 
		\inon{on~$\Gammac$}
		,
		\label{EQ42}
		\\
		&
		u_{0j}
		=
		0
		\inon{on~$\Gammaf$}
		,
		\label{EQ43}
	\end{align}
for $j=1,2,3$.
Suppose that
\begin{align}
	\begin{split}
		&
	(\RR, \RR^{-1})
	\in
	(L^\infty ((0,1), H^2(\Omegaf)) \cap H^1 ((0,1), L^\infty (\Omegaf)))^2
	\label{EQ332}
	\end{split}
\end{align}
and
\begin{align}
	\begin{split}
	(u_0|_{\Gammac}, \partial_3 u_0 |_{\Gammaf})
	\in H^{s+1/2} (\Gammac)
	\times 
	H^{s-1/2} (\Gammaf)
	\end{split}
\end{align}
with the nonhomogeneous terms satisfying
\begin{align}
		&
	( h, f, f(0))
	\in
	K^{s-1/2} ((0,1)\times \Gammac)
	\times 
	K^{s-1} ((0,1)\times \Omegaf)
	\times
	H^{s-2} (\Omegaf)
	.
	\label{EQ951}
\end{align}
Then the system \eqref{EQ48}--\eqref{EQ46} admits a solution $u$ satisfying
\begin{align}
	\begin{split}
	\Vert u\Vert_{K^{s+1} ( (0,1) \times \Omegaf )}
	&
	\les	
					\Vert h\Vert_{K^{s-1/2}_{\Gammac}}
	+
	\Vert u_0 |_{\Gammac} \Vert_{H^{s+1/2} (\Gammac)}
	+
	\Vert \partial_3 u_0 |_{\Gammaf} \Vert_{H^{s-1/2} (\Gammaf )}
	\\&\indeq
	+
	\Vert u_0\Vert_{H^s}
	+
	\Vert f  \Vert_{K^{s-1} ((0,1)\times \Omegaf)}
	+
	\Vert f (0) \Vert_{H^{s-2} (\Omegaf)}
	,
	\label{EQ92}
	\end{split}
\end{align}
where the implicit constant depends on the norms of $\RR$ and $\RR^{-1}$ in~\eqref{EQ332}.
\end{Lemma}

\begin{proof}
	In order to apply a lifting result in \cite{LM}, we consider the boundary conditions
\begin{align}
		&
	v
	=
	u_0 |_{\Gammac}
	\inon{on~$[0,1]\times\Gammac$}
	,
	\label{EQ1010}
	\\&
	\lambda (\partial_k v_j + \partial_j v_k)
	\nu^k
	+
	\mu \partial_k v_k \nu^j
	=
	h_j
	\inon{on~$[0,1]\times\Gammac$}
	,
	\label{EQ1043}
	\\&
	\partial_k \partial_m v_j \nu^k \nu^m
	=
	0
	\inon{on~$[0,1]\times\Gammac$}
	,
	\label{EQ1044}\\&
	v
	=0
	\inon{on~$[0,1]\times \Gammaf$}
	,
	\label{EQ1045}
	\\&
	\partial_k v_j \nu^k
	=\partial_k u_{0j} \nu^k
	\inon{on~$[0,1]\times \Gammaf$}
	,
	\label{EQ1048}
	\\&
	\partial_m \partial_k v_j \nu^k \nu^m
	=0
	\inon{on~$[0,1]\times \Gammaf$}
	,
	\label{EQ1047}
	\\&
	v~~\text{periodic in the $y_1$~and~$y_2$ directions},
	\label{EQ10320} 
\end{align}
for $j=1,2,3$,
and the initial data
\begin{align}
	&
	v(0) 
	= 
	u_0 
	\inon{in~$\Omegaf$}
	,
	\label{EQ1046}
	\\&
	\partial_t v(0) 
	= 
		\lambda \RR_0 \dive (\nabla u_0  + (\nabla u_0)^T)
	+
	\mu \RR_0 \nabla \dive u_0
	+ f (0)
	\inon{in~$\Omegaf$}
	.
	\label{EQ1049}
\end{align}
Below, the norm of dependence on time and space are understood as $(0,1)$ and $\Omegaf$, unless stated otherwise.
From \cite[Theorem~2.3]{LM} and the compatibility conditions \eqref{EQ42}--\eqref{EQ43} and since $s>1/2$ it follows that there exists $v \in K^{s+1} ((0,1)\times \Omegaf)$ satisfying the boundary conditions and initial conditions
 \eqref{EQ1010}--\eqref{EQ1049} with
\begin{align}
	\begin{split}
	\Vert v \Vert_{K^{s+1}}
	\les
	\Vert h\Vert_{K^{s-1/2}_{\Gammac}}
		+
	\Vert u_0  |_{\Gammac} \Vert_{K^{s+1/2}_{\Gammac}}
		+
	\left\Vert \frac{\partial u_0}{\partial \nu}  \right\Vert_{K^{s-1/2}_{\Gammaf}}
	+
	\Vert u_0\Vert_{H^s}
	+
		\Vert 	 \RR_0 D^2 u_0 \Vert_{H^{s-2}}
	+
	\Vert f (0) \Vert_{H^{s-2}}
	,
	\llabel{EQ83}
	\end{split}
\end{align}
from where
\begin{align}
	\begin{split}
	\Vert v \Vert_{K^{s+1}}
	&
	\les
	\Vert h\Vert_{K^{s-1/2}_{\Gammac}}
	+
	\Vert u_0 |_{\Gammac} \Vert_{H^{s+1/2} (\Gammac)}
	+
	\Vert   \partial_3 u_0|_{\Gammaf} \Vert_{H^{s-1/2} (\Gammaf)}
		+
	\Vert u_0\Vert_{H^s}
	+
	\Vert f (0) \Vert_{H^{s-2}}
	.
	\label{EQ83}
	\end{split}
\end{align}
\colb
Now we consider the homogeneous parabolic problem
\begin{align}
	\begin{split}
	w_t
	-
	\lambda \RR \dive(\nabla w + (\nabla w)^T)
	-
	\mu \RR \nabla \dive w 
	&
	= 
	F
	\inon{in~$[0,1]\times\Omegaf$}
        ,
	\label{EQ268}
	\end{split}
\end{align}
with the homogeneous boundary conditions and the initial data
\begin{align}
	&
	\lambda (\partial_k w_j + \partial_j w_k)
	 \nu^k
	+
	\mu \partial_k w_k \nu^j
	=
	0
	\inon{on~$[0,1]\times\Gammac$}
	,
	\\
	&
	w
	=0
	\inon{on~$[0,1]\times\Gammaf$}
	,
	\\
	&
	w~~\text{periodic in the $y_1$~and~$y_2$ directions},
	\\
	&
	w(0, \cdot)
	=
	0
	\inon{in~$\Omegaf$}
	,
	\label{EQ04}
\end{align}
for $j=1,2,3$,
where
\begin{align}
	\begin{split}
		F=
		v_t
		-f 
		-
		\lambda \RR \dive(\nabla v + (\nabla v)^T)
		-
		\mu \RR \nabla \dive v
		\inon{in~$[0,1]\times\Omegaf$}.
		\label{EQ85}
	\end{split}
\end{align}
Note that \eqref{EQ1049} implies that
  \begin{equation}
   F(0, \cdot)=0
   	\inon{in~$\Omegaf$}
   .
   \label{EQ221}
  \end{equation}
By \eqref{EQ332}, \eqref{EQ221}, and Lemma~\ref{L01}, there exists a solution $w$ to the system \eqref{EQ268}--\eqref{EQ85} satisfying
\begin{align}
	\begin{split}
	\Vert w\Vert_{K^2}
	\les
	\Vert F\Vert_{K^0}
	\label{EQ81}
\end{split}
\end{align}
and
\begin{align}
	\Vert w\Vert_{K^4}
	&
	\les
	\Vert F\Vert_{K^2}
	,
	\label{EQ82}
\end{align}
where the implicit constants depend on the norms of $\RR$ and $\RR^{-1}$ in~\eqref{EQ332}.
From \cite[Theorem~6.2]{LM} and \eqref{EQ81}--\eqref{EQ82} it follows that
\begin{align}
	\begin{split}
	\Vert w\Vert_{K^{s+1}}
	&
	\les
	\Vert F\Vert_{K^{s-1}}
	,
	\label{EQ84}
	\end{split}
\end{align}
since $s \notin 1/2 + \mathbb{Z}$ and $s/2 \notin \mathbb{Z}$.
From \eqref{EQ85}, we get
\begin{align}
	\begin{split}
	\Vert F\Vert_{K^{s-1}}
	\les
	\Vert f\Vert_{K^{s-1}}
	+
	\Vert v_t\Vert_{K^{s-1}}
	+
	\Vert R D_x^{2} v\Vert_{K^{s-1}}
	.
	\label{EQ86}
	\end{split}
\end{align}
For the second term on the right side of \eqref{EQ86}, we obtain
\begin{align}
	\begin{split}
	\Vert v_t\Vert_{K^{s-1}}
	&	
	\les
	\Vert v_t\Vert_{L^2_t H^{s-1}_x}	
	+
	\Vert v_t\Vert_{H^{(s-1)/2}_t L^2_x}	
	\les
	\Vert v\Vert_{K^{s+1}}
	,
	\llabel{EQ300}
	\end{split}
\end{align}
where we used Corollary~\ref{C02}.
To treat the last term on the right side of \eqref{EQ86}, we claim that
\begin{align}
	\Vert AB\Vert_{H^{(s-1)/2}_t L^2_x}
	\les
	\Vert A\Vert_{H^{1}_t L^\infty_x} 
	\Vert B\Vert_{H^{(s-1)/2}_t L^2_x}
	+
	\Vert A\Vert_{L^{\infty}_t L^\infty_x} 
	\Vert B\Vert_{H^{(s-1)/2}_t L^2_x}
	\label{EQ57}
\end{align}
on the domain $(0,1)\times\Omegaf$.
Using extensions, we may assume that the domain is actually $\mathbb{R}\times\mathbb{R}^{3}$.
From the H\"older inequality it follows that
\begin{align*}
	\begin{split}
		\Vert AB\Vert_{H^{(s-1)/2}_t L^2_x}
		&
	\les
	\Vert A\Vert_{W^{(s-1)/2, 4}_t L^\infty_x} 
	\Vert B\Vert_{L^{4}_t L^2_x}
	+
	\Vert A\Vert_{L^{\infty}_t L^\infty_x} 
	\Vert B\Vert_{H^{(s-1)/2}_t L^2_x}.
	\\&
	\les
		\Vert A\Vert_{W^{3/4, 4}_t L^\infty_x} 
	\Vert B\Vert_{L^{4}_t L^2_x}
	+
	\Vert A\Vert_{L^{\infty}_t L^\infty_x} 
	\Vert B\Vert_{H^{(s-1)/2}_t L^2_x}.
	\end{split}
\end{align*}
since $2<s<5/2$.
The claim \eqref{EQ57} is thus completed by appealing to the Sobolev inequality.
For the last term on the right side of \eqref{EQ86}, we use the H\"older's inequality, yielding
\begin{align}
	\begin{split}
	\Vert \RR D_x^2 v\Vert_{L^2_t H^{s-1}_x}
	&
	\les
	\Vert \RR  \Vert_{L^\infty_t H^2_x}
	\Vert D_x^{2} v \Vert_{L^2_t H^{s-1}_x}
	\les
	 \Vert v \Vert_{K^{s+1}}
	 \llabel{EQ88}
	\end{split}
\end{align}
and
\begin{align}
	\begin{split}
	\Vert \RR D_x^2 v\Vert_{H^{(s-1)/2}_t L^2_x}
	&
	\les
\Vert \RR \Vert_{H^1_t L^\infty_x}
\Vert D_x^2 v\Vert_{H^{(s-1)/2}_t L^2_x}
	+
	\Vert \RR \Vert_{L^\infty_t L^\infty_x}
	\Vert D_x^2 v\Vert_{H^{(s-1)/2}_t L^2_x}
	\les
	\Vert v\Vert_{K^{s+1}}
	,
	\label{EQ89}
	\end{split}
\end{align}
where we appealed to \eqref{EQ57} and Corollary~\ref{C02}.
Note that from \eqref{EQ268}--\eqref{EQ85}, we infer that the difference $u=v-w$ is a solution of the system~\eqref{EQ48}--\eqref{EQ46}.
From \eqref{EQ83}, \eqref{EQ84}--\eqref{EQ86}, and~\eqref{EQ89} it follows that
\begin{align}
	\begin{split}
		\Vert u\Vert_{K^{s+1}}
		&
		\les
		\Vert w\Vert_{K^{s+1}}
		+
		\Vert v\Vert_{K^{s+1}}	
		\\&
		\les
				\Vert h\Vert_{K^{s-1/2}_{\Gammac}}
		+
		\Vert u_0 |_{\Gammac} \Vert_{H^{s+1/2} (\Gammac)}
		+
		\Vert \partial_3 u_0 |_{\Gammaf} \Vert_{H^{s-1/2} (\Gammaf )}
		+
		\Vert u_0\Vert_{H^s}
		+
		\Vert f  \Vert_{K^{s-1}}
		+
		\Vert f (0) \Vert_{H^{s-2}}
		,
		\llabel{EQ90}
	\end{split}
\end{align}
concluding the proof of~\eqref{EQ92}.
\end{proof}

\startnewsection{Solution to a parabolic-wave system}{sec05}
In this section, we consider the coupled parabolic-wave system
\begin{align}
	&
	v_t - \lambda \RR \dive (\nabla v + (\nabla v)^T)
	-\mu \RR\nabla \dive v 
	+
	\RR \nabla ( \RR^{-1} )
	=
	 f
	 \inon{in~$[0, T]\times \Omegaf$}
	 ,
	 \label{EQ71a}
	 \\
	 &
	 \RR_t - \RR \dive v
	 =
	 0
	 \inon{in~$[0, T]\times \Omegaf$}
	 ,
	 \label{EQ93a}
	\\
	&
        w_{tt} - \Delta w 
	= 0
	\inon{in~$[0, T]\times\Omegae$}
	,
	\label{EQ72a}
\end{align}
with the boundary conditions
\begin{align}
	&
	v
	= 
	w_t
	\inon{on~$[0,T]\times\Gammac$}
	,
	\label{EQ73a}
	\\
	&
	\lambda (\partial_k v_j + \partial_j v_k) \nu^k
	+
	\mu \partial_k v_k \nu^j
	=
	\partial_k w_j \nu^k
	+
	\RR^{-1} \nu^j
	+
	h_j
	\inon{on~$[0, T]\times \Gammac$}
	,
	\label{EQ75a}
	\\&
	v, w~~\text{periodic in the $y_1$~and~$y_2$ directions},
	\label{EQ321a}
	\\&
	v
	=
	0
	\inon{on~$[0, T]\times\Gammaf$}
	,
	\label{EQ76a}
\end{align}
for $j =1,2,3$, 
and the initial data
\begin{align}
	\begin{split}
		&
	(v, \RR, w, w_t)(0) 
	= 
	(v_0, \RR_0, w_0, w_1)
	\inon{in~$\Omegaf \times \Omegaf 
	\times \Omegae \times \Omegae$}
	,
	\\
	&
	(v_0, \RR_0, w_0, w_1) ~\text{periodic in the $y_1$~and~$y_2$ directions}
     ,
     \\
     &
     w_0
     =
     0.
	\label{EQ77}
	\end{split}
\end{align}

In order to avoid issues of dependence of constants for small time, we introduce a cutoff function in time and work on the unit time interval~$(0,1)$.
Let $\TT \in (0,1/4)$, and let
$\phi_{\TT} (t)$ be a smooth cutoff function valued in $[0,1]$ such that
\begin{equation}
	\phi_{\TT} (t)=
	\begin{cases}
		1 \inon{on} ~~~[0,\TT],\\
		0 \inon{on} ~~~[2\TT, 1],
		\label{EQ991}
	\end{cases}
\end{equation}
and $\Vert \phi_{\TT}' \Vert_{L^\infty (0,1)} \les 1 /\TT$. 
The following lemma provides a necessary estimate for the cutoff function.
\begin{Lemma}
	\label{L14}
	We have
	$\Vert \phi_{\TT} \Vert_{ H^{(s-2)/2}_t}
	\les
	1$.
\end{Lemma}

\begin{proof}[Proof of Lemma~\ref{L14}]
By the Sobolev interpolation inequality, we have
 \begin{align*}
	\Vert \phi_{\TT} \Vert_{ H^{(s-2)/2}_t}
	&
	\les
	\Vert \phi_{\TT} \Vert_{H^1_t}^{(s-2)/2}
	\Vert \phi_{\TT} \Vert_{L^2_t}^{(4-s)/2}
	\les
	(1 + \TT^{-1/2})^{(s-2)/2}
	\TT^{(4-s)/4}
	\les
	\TT^{(3-s)/2}\les 1,
\end{align*}
since $s<3.$
\end{proof}
\colb

To obtain the existence of solutions and avoid issues with the dependence of constants for small time, we replace
\eqref{EQ71a}--\eqref{EQ72a} and \eqref{EQ73a}--\eqref{EQ76a} with
\begin{align}
	&
	v_t - \lambda \RR \dive (\nabla v + (\nabla v)^T)
	-\mu \RR\nabla \dive v 
	+
	\RR \nabla ( \RR^{-1} )
	=
	 f
	 \inon{in~$[0, 1]\times \Omegaf$}
	 ,
	 \label{EQ71}
	 \\
	 &
	 \RR_t - \phi_{\tilde T}\RR \dive v
	 =
	 0
	 \inon{in~$[0, 1]\times \Omegaf$}
	 ,
	 \label{EQ93}
	\\
	&
	w_{tt} - \Delta w 
	= 0
	\inon{in~$[0, 1]\times\Omegae$}
	,
	\label{EQ72}
\end{align}
with the boundary conditions
  \begin{align}
	&w(t,x) 
	= 
	\int_0^t \phi_{\TT} (\tau) v(\tau,x) \,d\tau
	+
	\left(t
        	- 
         	\int_0^t \phi_{\TT}(\tau) \,d\tau
        \right)
	v_0 (x)
	\inon{on~$[0, 1]\times \Gammac$}
	,
	\label{EQ73}
	\\&
	\lambda (\partial_k v_j + \partial_j v_k) \nu^k
	+
	\mu \partial_k v_k \nu^j
	=
	\partial_k w_j \nu^k
	+
	\RR^{-1} \nu^j
	+
	h_j
	\inon{on~$[0, 1]\times \Gammac$}
	,
	\label{EQ75}
	\\&
	v, w~~\text{periodic in the $y_1$~and~$y_2$ directions},
	\label{EQ321}
	\\&
	v
	=
	0
	\inon{on~$[0, 1]\times\Gammaf$}
	,
	\label{EQ76}
\end{align}
for $j =1,2,3$, where $\phi_{\TT} (t)$ is as in~\eqref{EQ991}.
Note that from \eqref{EQ73} it follows that
\begin{align}
	w_t(t, x)
	=
	\phi_{\TT} (t) (v(t, x) - v_0 (x))
	+
	v_0 (x)
	\inon{on~$[0, 1]\times \Gammac$}
	,
	\label{EQ007}
\end{align}
and thus the boundary condition \eqref{EQ007} agrees with \eqref{EQ73a} on the time interval $[0,\TT]$,
and the solutions of \eqref{EQ71}--\eqref{EQ76} agree with the
solution of \eqref{EQ71a}--\eqref{EQ76a} on the time interval~$[0,\tilde T]$, with the same initial and boundary conditions~\eqref{EQ77}.

To provide the maximal regularity for the system~\eqref{EQ71}--\eqref{EQ76}, we state the following necessary a~priori density estimates.
\cole
\begin{Lemma}
	\label{L04}
Let $s\in (2, 2 + \epsilon_0]$, where $\epsilon_0 \in (0,1/2)$ is arbitrary.
Consider the ODE system
\begin{align}
		&
		\RR_t - \RR \phi_{\TT} \dive v 
		= 
		0
		\inon{in~$[0, 1] \times \Omegaf$} 
		,
		\label{EQ146}
		\\
		&
		\RR(0)
		=
		\RR_0
		\inon{on~$\Omegaf$}
		.
		\label{EQ147}
\end{align}
Assume that
$
	(\RR_0, \RR_0^{-1}, v_0)
	\in
	H^s (\Omegaf)
	\times
	H^s (\Omegaf)
	\times H^s (\Omegaf)
$
and $\Vert v\Vert_{K^{s+1} ((0,1)\times \Omegaf)} \leq M$, where $M\geq 1$.
Let $\delta \in (0,1/2)$.
Then for a sufficiently small constant $\TT>0$, depending on $M$ and $\delta$, we have
\begin{enumerate}[label=(\roman*)]
	\item $\Vert \RR \Vert_{L^\infty_t L^\infty_x} +\Vert \RR^{-1} \Vert_{L^\infty_t L^\infty_x}+ \Vert \RR \Vert_{L^\infty_t H^s_x} 
	+
	\Vert \RR^{-1} \Vert_{L^\infty_t H^s_x} 
	 \les 1$	 
	 ,
\item $\Vert R^{-1}\Vert_{H^1_t H^{3/2+\delta}_x} +\Vert R\Vert_{H^1_t H^{3/2+\delta}_x} \les 1$,
	\item $ \Vert \RR\Vert_{H^1_t H^s_x} \les M$,
\end{enumerate}
where the norm of dependence is $(0,1) \times \Omegaf$.
\end{Lemma}
\colb

We emphasize that the implicit constants in the above inequalities (i)--(iii) are independent of~$M$ and~$\delta$.

\begin{proof}[Proof of Lemma~\ref{L04}]
(i) The solution of the ODE system \eqref{EQ146}--\eqref{EQ147} reads
\begin{align}
	\begin{split}
	\RR(t,x)
	=
	\RR_0 (x) e^{\int_0^t 
	\phi_{\TT} (\tau) \dive v(\tau) \,d\tau}
	\inon{in~$[0, 1]\times \Omegaf$}
       .
	\label{EQ148}
	\end{split}
\end{align}
Let $\tilde T\in (0,T]$ be a small time to be determined below.
From H\"older's and Sobolev inequalities it follows that
\begin{align}
	\begin{split}
	\Vert \RR \Vert_{L^\infty_t L^\infty_x} 
	\les
	\Vert \RR_0 \Vert_{H^s} 
	e^{\int_0^{2\TT}  \Vert \phi_{\TT} (\tau)\dive v(\tau)  \Vert_{L^\infty} \,d\tau }
	\les 
	C^{\TT^{1/2} M }
	\les
	1
	\llabel{EQ150}
\end{split}
\end{align}
and
\begin{align}
	\begin{split}
		\Vert \RR^{-1} \Vert_{L^\infty_t L^\infty_x} 
		\les
		\Vert \RR_0^{-1} \Vert_{H^s} 
		e^{ \int_0^{2\TT}  
			\Vert \phi_{\TT} (\tau) 
			\dive v(\tau)  \Vert_{L^\infty} \,d\tau }
		\les 
		C^{\TT^{1/2} M}
		\les
		1
		,
		\llabel{EQ151}
	\end{split}
\end{align}
for some sufficiently small~$\TT >0$.
Similarly, we have
\begin{align}
	\begin{split}
		\Vert \RR \Vert_{L^\infty_t H^s_x}
		\les
		\Vert \RR_0 \Vert_{H^s}	
		\Vert e^{\int_0^{2\TT}
		\phi_{\TT}(\tau) \dive v(\tau) \,d\tau} \Vert_{L^\infty_t H^s_x}
		\les
		1
		\llabel{EQ187}
	\end{split}
\end{align}
and
\begin{align}
\begin{split}
	\Vert \RR^{-1} \Vert_{L^\infty_t H^s_x}
		\les
	\Vert \RR_0^{-1} \Vert_{H^s}	
	\Vert e^{\int_0^{2\TT}
		\phi_{\TT}(\tau) \dive v(\tau) \,d\tau} 
	\Vert_{L^\infty_t H^s_x}
	\les
	1
	.
	\llabel{EQ186}
\end{split}	
\end{align}

(ii) From \eqref{EQ146}, we use H\"older's and the Sobolev inequalities to get
\begin{align}
	\begin{split}
		\Vert (\RR^{-1})_t \Vert_{L^2_t H^{3/2+\delta}_x}
		\les
		\Vert \RR^{-2} \RR_t \Vert_{L^2_t H^{3/2+\delta}_x}
		\les	
		\Vert \dive v \Vert_{L^2  ((0,2\TT), H^{3/2+\delta} (\Omegaf))}
				\les
	\Vert v\Vert_{L^2  ((0,2\TT), H^{5/2+\delta} (\Omegaf))}
	.
	\label{EQ205}
	\end{split}
\end{align}
Recall that for any $0<r<r'$ and $f \in H^{r'}$, we have the Sobolev interpolation inequality
\begin{align}
	\begin{split}
		\Vert f\Vert_{H^r }
		\les
		\epsilon \Vert f\Vert_{H^{r'} }
		+
		\epsilon^{r/(r-r')} \Vert f\Vert_{L^2}
		,
		\label{EQ203}
	\end{split}
\end{align}
for any $\epsilon\in (0,1]$.
From \eqref{EQ205}--\eqref{EQ203} it follows that
\begin{align}
	\begin{split}
	\Vert (R^{-1})_t \Vert_{L^2_t H^{3/2+\delta}_x}
	&
	\les
	\epsilon \Vert v\Vert_{L^2  ((0,2\TT), H^{s+1} (\Omegaf))}
	+
	C_\epsilon
	\Vert v\Vert_{L^2  ((0,2\TT), L^{2} (\Omegaf))}
	\les
	\epsilon M
	+C_\epsilon \TT^{1/2} \Vert v\Vert_{L^\infty_t L^2_x}
		\\&
		\les
		\epsilon M
	+C_\epsilon \TT^{1/2}\Vert v\Vert_{H^{(s+1)/2}_t L^2_x}
	\les
	(	\epsilon 
		+
	C_\epsilon	
	\TT^{1/2})
	 M
	,
	\label{EQ230}
	\end{split}
\end{align}
since~$s>2$.
Taking $\epsilon=1/M$ in \eqref{EQ230}, we arrive at
\begin{align}
	\begin{split}
		\Vert (R^{-1})_t \Vert_{L^2_t H^{3/2+\delta}_x}
		\les
		1,
	\end{split}
   \llabel{EQ131}
\end{align}
for some sufficiently small~$\TT>0$.
Similarly, we have
\begin{align}
	\begin{split}
		\Vert R_t \Vert_{L^2_t H^{3/2+\delta}_x}
		\les
		1.
	\end{split}
   \llabel{EQ132}
\end{align}
Thus, we conclude the proof of (ii) by combining (i).

(iii) From \eqref{EQ146} and H\"older's inequality it follows that
\begin{align}
	\begin{split}
		\Vert \RR_t \Vert_{L^2_t H^s_x}
		\les 
		\Vert \RR  \phi_{\TT} \Vert_{L^\infty_t H^s_x}
		\Vert \dive v\Vert_{L^2_t H^s_x }	
		\les
		\Vert v\Vert_{L^2_t H^{s+1}_x}
		\les
		M
		.
	\end{split}
   \llabel{EQ133}
\end{align}
Therefore, we conclude the proof of (iii).
\end{proof}

The following lemma provides necessary estimates for the structure displacement and velocity on the boundary.
\begin{Lemma}
\label{Lomega}
Let $s\in (2, 2 + \epsilon_0]$, where $\epsilon_0 \in (0,1/2)$ is arbitrary. 
Assume that $\Vert v\Vert_{K^{s+1} ((0,1)\times \Omegaf)} \leq M$ for some~$M\geq 1$. 
Suppose that $v$ and $w$ satisfy \eqref{EQ73} and \eqref{EQ321} with the initial data satisfying
$
	(v_0, w_0, w_1)
	\in
	H^s(\Omegaf) \times H^{s+ 1/2} (\Omegae)
	\times
	H^{s-1/2} (\Omegae)
$
and
$v_0 |_{\Gammac}
	\in
	H^{s+1/2} (\Gammac)$.
Then we have
\begin{enumerate}[label=(\roman*)]
	\item $\Vert w \Vert_{L^2_t H_x^{s+1/2} (  \Gammac)}
	\les 
	\TT^{1/2} M
	+1$,
	\item $
	\Vert w_t\Vert_{H^{s/2-3/4}_t H^{s/2+1/4}_x ( \Gammac)}
	+
	\Vert w \Vert_{H^{s/2-3/4}_t H_x^{s/2+1/4} (\Gammac)}
	\les 
	(\epsilon 
	+
	\tilde{\epsilon} C_\epsilon 
	+
	C_{\tilde{\epsilon}, \epsilon}
	\TT^{1/2})
	M
	+
	C_\epsilon$,
	\item $\Vert w\Vert_{H^{s/2+3/4}_t L^2_x (\Gammac)}
	 \les 
	(\epsilon
	+	\tilde{\epsilon}
	C_\epsilon  
	+
	C_{\epsilon, \tilde{\epsilon}} \TT^{1/2}) M
	+ C_\epsilon$,
\end{enumerate}
for any $\epsilon, \tilde{\epsilon} \in (0,1]$, where the implicit constants depend on the initial data. 
\end{Lemma}

Here and below, when not indicated, the time and space domains are
understood to be $(0, 1)$ and~$\Omegaf$, respectively.

\begin{proof}[Proof of Lemma~\ref{Lomega}]
(i) Using \eqref{EQ73} we get
\begin{align}
	\begin{split}
		\Vert w \Vert_{L^2_t H_x^{s+1/2}( \Gammac)}
		&
		\les
		\left(
		\int_{0}^{1}
		\left\Vert \int_0^{t}
		\phi_{\TT}  v \,d\tau
		\right\Vert^2_{
			H^{s+1/2} (\Gammac)} dt
		\right)^{1/2}
		+
		\left(
		\int_{0}^1
		\left\Vert 
		\left(t- \int_0^t
		\phi_{\TT}
		\right)
		v_0
		\right\Vert^2_{
			H^{s+1/2} (\Gammac)} dt
		\right)^{1/2}
		\\&
		\les
		\TT^{1/2} \Vert v\Vert_{L^2_t H^{s+1}_x}
		+
		1
		\les
		\TT^{1/2} M
		+
		1
		,
		\label{EQ350}
	\end{split}
\end{align}
since $v_0|_{\Gammac} \in H^{s+1/2} (\Gammac)$,
where we also used that for every $t\in[0,1]$ we have
\begin{equation}
	\left(
	\int_0^{t}
	\phi_{\TT}
	\Vert v\Vert_{H^{s+1/2} (\Gammac)}
	\,d\tau
	\right)^2
	=
	\left(
	\int_0^{2\tilde T}
	\phi_{\TT}
	\Vert v\Vert_{H^{s+1/2} (\Gammac)}
	\,d\tau
	\right)^2
	\lec
	\tilde T
	\int_0^{2\tilde T}
	\Vert v\Vert_{H^{s+1/2} (\Gammac)}^2
	\,d\tau
	.
	\llabel{EQ51}
\end{equation}

(ii) We use the Sobolev interpolation and Young inequalities to write
\begin{align}
	\begin{split}
		\Vert w_t\Vert_{H^{s/2-3/4}_t H^{s/2+1/4}_x (\Gammac)}
		&\les	
		\Vert w_t\Vert_{H^{1}_t H^{s-3/2}_x (\Gammac)}^{s/2-3/4}
		\Vert w_t\Vert_{L^2_t H_x^{(-4s^2+16s-7)/2(7-2s)} (\Gammac)}^{7/4-s/2}
		\\&
		\les
		\epsilon
		\Vert w_t\Vert_{H^{1}_t H^{s-3/2}_x (\Gammac)}
		+
		C_\epsilon
		\Vert w_t\Vert_{L^2_t H_x^{(-4s^2+16s-7)/2(7-2s)} (\Gammac)}
		:=\mathcal{I}_1
		+
		\mathcal{I}_2
		,
		\label{EQ599}
	\end{split}
\end{align}
for any $\epsilon \in (0,1]$.

Note that the implicit constant in the first inequality is independent of $\TT$ since the interpolation is applied on a fixed domain $(0,1)\times \Gammac$.
For the term $\mathcal{I}_1$, we use \eqref{EQ007}, the trace inequality, and the Leibniz rule, to obtain
\begin{align}
	\begin{split}
		\mathcal{I}_1
		&
		\les
		\epsilon \Vert \phi_{\TT}' (v-v_0)\Vert_{L^2_t H^{s-1}_x}
		+
		\epsilon \Vert  v' \Vert_{L^2_t H^{s-1}_x}
		+
		\epsilon \Vert \phi_{\TT} (v - v_0) +v_0 \Vert_{L^2_t H^{s-1}_x}
		=:
		\mathcal{I}_{11}
		+
		\mathcal{I}_{12}
		+
		\mathcal{I}_{13}
		.
		\label{EQ351}
	\end{split}
\end{align}
The term $\II_{11}$ is estimated using the Sobolev and H\"older inequalities as
\begin{align}
	\begin{split}
		\II_{11}
		&
		\les
		\epsilon
		\TT^{-1}
		\Vert v - v_0\Vert_{L^2  ((0,2\TT), H^{s-1} (\Omegaf))}
		\les
		\epsilon \TT^{-1}
		\left\Vert  \int_0^t v_t \right \Vert_{L^2 ((0,2\TT), H^{s-1} (\Omegaf)) }
		\les
		\epsilon 
		\Vert v'\Vert_{L^{2}_t H^{s-1}_x }
		\les
		\epsilon M
		,
		\label{EQ588}
	\end{split}
\end{align}
where we used
\begin{align}
	\begin{split}
		&
		\left\Vert  \int_0^t v_t \right\Vert_{L^2 ((0,2\TT), H^{s-1} (\Omegaf)) }^2
		\leq
		\int_{0}^{2\tilde T}
		\left(
		\int_{0}^{2\tilde T}
		\Vert v_t(s)\Vert_{H^{s-1}}
		\,ds
		\right)^2
		\,dt
		\\&\indeq
		= 2\tilde T
		\left(
		\int_{0}^{2\tilde T}
		\Vert v_t(s)\Vert_{H^{s-1}}
		\,ds
		\right)^2
		\lec
		\tilde T^2
		\int_{0}^{2\tilde T}
		\Vert v_t(s)\Vert_{H^{s-1}}^2
		\,ds
	\end{split}
	\label{EQ74}
\end{align}
in the second inequality
and
Corollary~\ref{C02} in the last.
Next, the terms $\II_{12}$ and $\II_{13}$ are estimated as
\begin{align}
	\begin{split}
		\II_{12}
		&
		\les
		\epsilon M
	\end{split}
\end{align}
and
\begin{align}
	\begin{split}
		\II_{13}
		=
		\epsilon \Vert \phi_{\TT} v  + (1-\phi_{\tilde T} )v_0 \Vert_{L^2_t H^{s-1}_x}
		\lec
		\epsilon 
		\Vert v\Vert_{L^2_t H^{s-1}_x}
		+
		\Vert (1-\phi_{\TT} ) v_0\Vert_{L^2_t H^{s-1}_x}
		\les
		\epsilon M
		+1.
		\label{EQ366}
	\end{split}
\end{align}
For the term $\II_2$, we use \eqref{EQ007}, \eqref{EQ203}, and the trace inequality to get
\begin{align}
	\begin{split}
		\II_2
		&
		\les
		C_\epsilon
		\Vert v\Vert_{L^2  ((0,2\TT), H^{(-4s^2+16s-7)/(14-4s)+1/2} (\Omegaf))}
		+
		C_\epsilon
		\Vert v_0 (1-\phi_{\TT}) \Vert_{L^2  ((0,2\TT), H^{(-4s^2+16s-7)/(14-4s)+1/2} (\Omegaf))}
		\\&
		\les
		C_\epsilon
		\Vert v\Vert_{L^2  ((0,2\TT), H^s (\Omegaf))}
		+
		C_\epsilon
		\Vert v_0  \Vert_{L^2  ((0,2\TT), H^{s} (\Omegaf))}
		,
	\end{split}
	\llabel{EQ153}
\end{align}
where the last inequality follows from the identity
$
(-4s^2+16s-7)/(14-4s) +1/2
=s
$.
Using the Sobolev interpolation inequality, we get
\begin{align}
	\begin{split}
		\II_2
		\les
		\tilde{\epsilon} 	C_\epsilon \Vert v\Vert_{L^2  ((0,2\TT), H^{s+1} (\Omegaf))}
		+
		C_{\tilde{\epsilon}, \epsilon} \Vert v\Vert_{L^2  ((0,2\TT), L^2 (\Omegaf))}
		+
		C_\epsilon
		\les
		\tilde{\epsilon} C_\epsilon M
		+
		C_{\tilde{\epsilon}, \epsilon}
		\TT^{1/2}
		M
		+
		C_\epsilon
		,
		\label{EQ590}	
	\end{split}
\end{align}
for any $\tilde{\epsilon} \in (0,1]$.
Combining \eqref{EQ599}--\eqref{EQ590}, we arrive at
\begin{align}
	\begin{split}
		\Vert w_t
		\Vert_{H^{s/2-3/4}_t H^{s/2+3/4}_x	
			(\Gammac)}
		\les	
		&
		\epsilon M
		+
		\tilde{\epsilon} C_\epsilon M
		+
		C_{\tilde{\epsilon}, \epsilon}
		\TT^{1/2}
		M
		+
		C_\epsilon
		.
		\label{EQ365}
	\end{split}
\end{align}

For the second term on the left side of (ii), we proceed as in \eqref{EQ599}, obtaining
\begin{align}
	\begin{split}
		\Vert w\Vert_{H^{s/2-3/4}_t H^{s/2+1/4}_x (\Gammac)}
		&\les
		\Vert w\Vert_{H^1_t H^{1/2}_x (\Gammac)}^{s/2-3/4}	
		\Vert w\Vert_{L^2_t H^{(2s+5)/2(7-2s)}_x (\Gammac)}^{7/4-s/2}
		\\&
		\les
		\Vert w\Vert_{H^{1}_t H^{1/2}_x (\Gammac)}
		+
		\Vert w\Vert_{L^2_t H_x^{(2s+5)/2(7-2s)} (\Gammac)}
		\\&
		\lec
		\Vert 	\phi_{\TT}  v  \Vert_{L^2_t H^1_x}
		+
		\Vert (1-\phi_{\TT})
		v_0\Vert_{L^2_t H^1_x}
		+
		\Vert w\Vert_{L^2_t H^{(2s+5)/(14-4s)}_x (\Gammac)}
		,
	\end{split}
	\llabel{EQ154}
\end{align}
since $1/2\leq(2s+5)/(14-4s)$.
Note that $(2s+5)/(14-4s)<s+1/2$ for $2<s<5/2$.
Thus, using \eqref{EQ203} and \eqref{EQ350}, we obtain
\begin{align}
	\Vert w\Vert_{H^{s/2-3/4}_t H^{s/2+1/4}_x (\Gammac)}
	\les
	(\epsilon+ C_\epsilon \TT^{1/2} )
	M
	+1,
	\label{EQ009}
\end{align}
for any $\epsilon \in (0,1]$.

(iii) First, we write
\begin{align}
	\begin{split}
		\Vert w\Vert_{H^{s/2+3/4}_t L^2_x (\Gammac)}
		&
		\les
		\Vert w_t\Vert_{H^{s/2-1/4}_t L^2_x (\Gammac)}	
		+
		\Vert w\Vert_{L^2_t L^2_x (\Gammac)}
		\\&
		\les
		\Vert w_t\Vert_{H^{1}_t L^2_x (\Gammac)}^{s/2-1/4}	
		\Vert  w_t\Vert_{L^2_t L^2_x (\Gammac)}^{5/4-s/2}	
		+
		\Vert w\Vert_{L^2_t H^{s+1/2}_x (\Gammac)}
		\\&
		\les
		\epsilon \Vert w_{tt}\Vert_{ L^2_t L^2_x (\Gammac)}
		+
		C_\epsilon \Vert w_t\Vert_{L^2_t L^2_x (\Gammac)}
		+\TT^{1/2} M +1
		,
		\label{EQ970}
	\end{split}
\end{align}
for any $\epsilon \in (0,1]$,
where the last inequality follows from~\eqref{EQ350}.
Note that the implicit constant in the second inequality is independent of $\TT$ since the interpolation is performed on $(0,1)\times \Gammac$.
From \eqref{EQ007} it follows that
\begin{align}
	w_{tt}(t) 
	= 
	\phi_{\TT}'(t)
	(v(t)-v_0)
	+
	\phi_{\TT} (t)
	v_t(t, x) 
	\inon{on~$[0, 1]\times \Gammac$}
	.
	\label{EQ960}
\end{align}
For the first term on the far right side of \eqref{EQ970}, we use \eqref{EQ960} and obtain
\begin{align}
	\begin{split}
		\epsilon
		\Vert w_{tt} \Vert_{L^2_t L^2_x (\Gammac)}	
		&
		\les
		\epsilon \Vert \phi_{\TT}' (v-v_0) \Vert_{L^2_t L^2_x (\Gammac)}	
		+
		\epsilon
		\Vert \phi_{\TT} v_t \Vert_{L^2_t L^2_x (\Gammac)}	
		\\&
		\les
		\epsilon \TT^{-1}
		\Vert v-v_0 \Vert_{L^2  ( (0,2\TT), H^1 (\Omegaf))}
		+
		\epsilon
		\Vert v_t\Vert_{L^2_t H^1_x}
		\les
		\epsilon M
		,
	\end{split}
	\label{EQ78}
\end{align}
where the last inequality follows from \eqref{EQ588} and Corollary~\ref{C02}.
For the second term on the far right side of \eqref{EQ970},
we use \eqref{EQ007} to arrive at
\begin{align}
	\begin{split}
		C_\epsilon
		\Vert w_t\Vert_{L^2_t L^2_x (\Gammac)}
		&
		\les	
		C_\epsilon
		\Vert \phi_{\TT} v\Vert_{L^2_t L^2_x (\Gammac)}
		+
		C_\epsilon
		\Vert v_0 (1-\phi_{\TT}) \Vert_{L^2_t L^2_x (\Gammac)}
		\les	
		(	\tilde{\epsilon}
		C_\epsilon  
		+
		C_{\epsilon, \tilde{\epsilon}} \TT^{1/2}) M
		+ C_\epsilon
		,
		\label{EQ971}
	\end{split}
\end{align}
for any $\tilde{\epsilon} \in (0,1]$, where we used the trace inequality and~\eqref{EQ203}.
The proof of (iii) is concluded by combining \eqref{EQ970} and \eqref{EQ78}--\eqref{EQ971}.

\end{proof}
\colb

The following theorem provides the local existence for the parabolic-wave system~\eqref{EQ71}--\eqref{EQ76}.

\begin{Theorem}
\label{T03}
Let $s\in (2, 2+\epsilon_0]$, where $\epsilon_0 \in (0, 1/2)$.
Assume the compatibility conditions
\begin{align}
	\begin{split}
		&
	w_{1j}
	=
	v_{0j}
	\inon{on~$\Gammac$}
	,
	\\&
	v_{0j} 
	= 
	0
	\inon{on~$\Gammaf$}
	,
	\\&
	\lambda (\partial_k v_{0j} 
	+
	\partial_j v_{0k}
	)\nu^k
	+
	\mu	\partial_i v_{0i} \nu^j - \RR_0^{-1} \nu^j 
	-
	\partial_k w_{0j} \nu^k
	=
	h_j (0)
	\inon{on~$\Gammac$}
	,
	\\&
	\lambda R_0
	\partial_k (\partial_k v_{0j}
	+ 
	\partial_j v_{0k})
	+
	\mu R_0\partial_j \partial_k v_{0k} 
	-
	R_0
	\partial_j (\RR^{-1}_0 )
	=
	-
	f_j (0)
	\inon{on~$\Gammaf$}
	,
	\label{EQ140}	
	\end{split}
\end{align}
for $j=1,2,3$.
Suppose that the initial data satisfy
\begin{align}
	&
	(v_0, w_0, w_1, \RR_0^{-1}, \RR_0, f(0) )
	\in
	H^s(\Omegaf) \times H^{s+ 1/2} (\Omegae)
	\times
	H^{s-1/2} (\Omegae)
	\times
	H^s(\Omegaf)
	\times
	H^s(\Omegaf)
		\times 
	H^{s-2} (\Omegaf)
   \llabel{EQ143}
\end{align}
and
\begin{align}
	&
	( v_0 |_{\Gammac}, \partial_3 v_0 |_{\Gammaf})
	\in
	H^{s+1/2} (\Gammac)
	\times 
	H^{s-1/2} (\Gammaf)
	\llabel{EQ143}
\end{align}
with the nonhomogeneous terms satisfying
\begin{align}
	(f, h) 
	\in
	K^{s-1} ((0, 1) \times \Omegaf) 
	\times 
	K^{s- 1/2} ((0, 1) \times \Gammac )
	.
   \llabel{EQ145}
\end{align}
Then there exists a unique solution 
\begin{align}
	\begin{split}
	(v, \RR, w, w_t)
	&
	\in
	K^{s+1} ((0,\TT) \times \Omegaf)
	\times
	H^1 ((0, \TT), H^s(\Omegaf))
	\times
	C([0,\TT], H^{s+ 1/4-\epsilon_0}(\Omegae))
		\\&\indeqtimes
	C([0,\TT], H^{s- 3/4-\epsilon_0}(\Omegae))
        ,
	\end{split}
   \llabel{EQ149}
\end{align}
to the system \eqref{EQ71}--\eqref{EQ76},
where $\TT>0$ is a constant and the corresponding norms are bounded by a function of the initial data and the nonhomogeneous terms.
\end{Theorem}

Let
\begin{align}
	\begin{split}
	\mathcal{Z}
	&
	= 
	\{v \in K^{s+1} ((0,1) \times \Omegaf):
	v(0) = v_0 ~~\text{in}~~\Omegaf, 
	v=0 ~~\text{on}~~ [0,1] \times\Gammaf,
		\\&\indeq\indeq\indeq
	v~~\text{periodic in the $y_1$~and~$y_2$ directions},
	~~\text{and}~~ \Vert v\Vert_{K^{s+1} ((0,1) \times \Omegaf)} \leq M \}
	,
	\label{EQ188}
	\end{split}
\end{align}
where $M\geq 1$ is a constant to be determined below.
For $v\in \mathcal{Z}$, define $\RR$ by~\eqref{EQ148}.
Next, we solve the wave equation \eqref{EQ72} for $w$ with the boundary condition
\eqref{EQ73} and the initial data
$(w, w_t) (0) 
	= 
	(w_0, w_1)$
in~$\Omegae$.
With $(\RR, w)$ constructed this way, we define a mapping
\begin{align}
	\Lambda \colon v (\in \mathcal{Z} )\mapsto \bar{v}
	,
   \llabel{EQ152}
\end{align}
where $\bar{v}$ is the solution of the nonhomogeneous parabolic problem
\begin{align}
	\begin{split}
		&
		\bar{v}_t
		- 
		\lambda \RR \dive 
		(\nabla \bar{v}
		+ 
		(\nabla \bar{v})^T
		)
		-
		\mu \RR \nabla \dive \bar{v}
		= 
		f
		-
		\RR \nabla \RR^{-1}
		\inon{in~$[0, 1]\times \Omegaf$}
		,
		\label{EQ124}
	\end{split}
\end{align}
with the boundary conditions and the initial data
\begin{align}
	\begin{split}
&
\lambda (\partial_k \bar{v}_j 
+ 
\partial_j \bar{v}_k)
\nu^k
+
\mu \partial_k \bar{v}_k  \nu^j
=
\partial_k w_j \nu^k
+
\RR^{-1} \nu^j
+
h_j
\inon{on~$[0, 1]\times\Gammac$}
,
\\
&
\bar{v}
=
0
\inon{on~$[0,1]\times\Gammaf$}
,
\\&
\bar{v}~~\text{periodic in the $y_1$~and~$y_2$ directions},
\\&
\bar{v}(0) = v_0
\inon{in~$\Omegaf$}
,
\label{EQ125}
	\end{split}
\end{align}
for $j=1,2,3$. 
We shall prove below that $\Lambda$ is a contraction mapping and then use the Banach fixed-point theorem.

\subsection{Uniform boundedness of the iterative sequence}
\label{sec051}
In this section, we show that the mapping $\Lambda$ is well-defined from $\mathcal{Z}$ to $\mathcal{Z}$, for some sufficiently large constant~$M\geq 1$.
Let $\TT \in (0, 1/4)$ be a constant.
We emphasize that the implicit constants in this section  below depend on the initial data but are independent of $M$ and~$\TT$.
Denote the right side of \eqref{EQ125}$_1$ by~$\tilde{h}_j$. One may easily verify that
\begin{align}
	\tilde{h}_j (0)
	=
		\lambda (\partial_k v_{0j} 
	+
	\partial_j v_{0k}
	)\nu^k
	+
	\mu	\partial_i v_{0i} \nu^j
	\inon{on~$\Gammac$}
	\label{EQ1112}
\end{align}
by~\eqref{EQ140}$_3$.
From \eqref{EQ140}$_2$, \eqref{EQ1112}, Lemma~\ref{L02}, and Lemma~\ref{L04}, it follows that 
\begin{align*}
	\begin{split}
	\Vert \bar{v} \Vert_{K^{s+1} }
	&
	\les	
	\Vert \tilde{h}
	\Vert_{K^{s-1/2}_{\Gammac}}
	+
	\Vert v_0|_{\Gammac} \Vert_{H^{s+1/2} (\Gammac)}
	+
	\Vert \partial_3 v_0|_{\Gammaf} \Vert_{H^{s-1/2} (\Gammaf)}
	+
	\Vert v_0\Vert_{H^s}
	+
	\Vert f\Vert_{K^{s-1} }
		\\&\indeq
	+
	\Vert \RR^{-1} \nabla \RR \Vert_{K^{s-1} }
	+
	\Vert f(0)\Vert_{H^{s-2}}
	+
	\Vert \RR_0 \nabla \RR_0^{-1}\Vert_{H^{s-2}}
	\end{split}
\end{align*}
from where
\begin{align}
	\begin{split}
		\Vert \bar{v} \Vert_{K^{s+1} }
		&
		\les
		\left\Vert \frac{\partial w}{\partial \nu} \right\Vert_{K^{s-1/2}_{\Gammac}}
		+
		\Vert \RR^{-1} \nabla \RR \Vert_{K^{s-1} }
			+
		\Vert \RR^{-1} \Vert_{K^{s-1/2}_{\Gammac}}
				+
		\Vert h\Vert_{K^{s-1/2}_{\Gammac}}
		+
		\Vert f\Vert_{K^{s-1} }
		+
		\Vert v_0\Vert_{H^s}
		\\&\indeq
		+
		\Vert v_0 |_{\Gammac} \Vert_{H^{s+1/2} (\Gammac)}
				+
		\Vert \partial_3 v_0 |_{\Gammaf} \Vert_{H^{s-1/2} (\Gammaf)}
		+
		\Vert f(0)\Vert_{H^{s-2}}
		+
		\Vert \RR_0^{-1} \nabla \RR_0\Vert_{H^{s-2}}
		.
		\label{EQ111}
	\end{split}
\end{align}
Here and below, when not indicated, the time and space domains are
understood to be $(0, 1)$ and~$\Omegaf$, respectively.

For the space component of the first term on the right side of \eqref{EQ111}, we appeal to Lemma~\ref{L13} to obtain
\begin{align}
	\begin{split}
		\left\Vert \frac{\partial w}{\partial \nu} \right\Vert_{L^2_t H_x^{s-1/2}( \Gammac)}	
		&
		\les
		\Vert w_0\Vert_{H^{s+1/2} (\Omegae)}
		+
		\Vert w_1\Vert_{H^{s-1/2} (\Omegae)}
		+
		\Vert w \Vert_{L_t^2 H_x^{s+1/2} (\Gammac)}
		+
		\Vert w \Vert_{H^{s/2+1/4}_t H_x^{s/2+1/4} (\Gammac)}
		\\&
		\les
		\Vert w \Vert_{L^2_t H_x^{s+1/2} (  \Gammac)}
		+
		\Vert w_t\Vert_{H^{s/2-3/4}_t H^{s/2+1/4}_x ( \Gammac)}
		+
		\Vert w \Vert_{H^{s/2-3/4}_t H_x^{s/2+1/4} (\Gammac)}
		+1
	.
		\label{EQ115}
	\end{split}
\end{align}
From \eqref{EQ115} and Lemma~\ref{Lomega} it follows that
\begin{align}
	\begin{split}
	\left\Vert \frac{\partial w}{\partial \nu} \right\Vert_{L^2_t H_x^{s-1/2}( \Gammac)}
	\les	
	&
	(\epsilon 
	+
	\tilde{\epsilon} C_\epsilon 
	+
	C_{\tilde{\epsilon}, \epsilon}
	\TT^{1/2})
	M
	+
	C_\epsilon,
	\label{EQ370}
	\end{split}
\end{align}
for any $\epsilon, \tilde{\epsilon}\in (0,1]$.
For the time component of the first term on the right side of \eqref{EQ111}, we use Lemma~\ref{L03} to get
\begin{align}
	\begin{split}
		\left\Vert \frac{\partial w}{\partial \nu}\right\Vert_{H^{s/2-1/4}_t L_x^2(\Gammac)}	
		&
		\les
		\Vert w_0\Vert_{H^{s/2+3/4} (\Omegae)}
		+
		\Vert w_1\Vert_{H^{s/2-1/4} (\Omegae)}
		+
		\Vert w\Vert_{ L^2_t H^{s/2+3/4}_x (\Gammac)}				
		\\&\indeq
		+
		\Vert w\Vert_{H^{s/2+3/4}_t L^2_x (\Gammac)}
		.
		\label{EQ353}
	\end{split}
\end{align}
For the third term on the right side of \eqref{EQ353}, we appeal to \eqref{EQ350} to get
\begin{align}
	\Vert w\Vert_{ L^2_t H^{s/2+3/4}_x (\Gammac)}
	\les
	\TT^{1/2} M +1
	,
	\label{EQ972}
\end{align}
since $s/2+3/4 \leq s+1/2$.
Applying Lemma~\ref{Lomega} and \eqref{EQ972} in \eqref{EQ353}, we get
\begin{align}
	\left\Vert \frac{\partial w}{\partial \nu}\right\Vert_{H^{s/2-1/4}_t L_x^2(\Gammac)}	
	\les
	(\epsilon 
	+ 
	\tilde{\epsilon}C_\epsilon
	+
	C_{\epsilon, \tilde{\epsilon}}
	\TT^{1/2}) M		
	+
	C_\epsilon
	,
	\label{EQ973}
\end{align}
for any $\epsilon, \tilde{\epsilon}\in (0,1]$.

For the space component of the second term on the right side of \eqref{EQ111}, we use the H\"older's and the Sobolev inequalities to obtain
\begin{align}
	\begin{split}
	\Vert 
	\RR^{-1} \nabla \RR\Vert_{L^2_t H^{s-1}_x}
	&
	\les
	\Vert \RR^{-1} \Vert_{L^\infty_t H^s_x} 
	\Vert
	\nabla \RR \Vert_{L^2_t  H^{s-1}_x}
	\les
	1
	,
	\label{EQ214}
	\end{split}
\end{align}
where we appealed to Lemma~\ref{L04}.
For the time component, we use H\"older's and the Sobolev inequalities
with Lemma~\ref{L04} to obtain
\begin{align}
	\begin{split}
	\Vert 
	\RR^{-1} \nabla \RR\Vert_{H^{(s-1)/2}_t L^2_x}
	&\les
	\Vert 
	\RR^{-1} \nabla \RR
	\Vert_{H^1_t L^2_x}
	\les
	\Vert \RR_t \nabla R\Vert_{L^2_t L^2_x}
	+
	\Vert \nabla \RR_t\Vert_{L^2_t L^2_x}
	+
	\Vert \nabla R\Vert_{L^2_t L^2_x}
	\\&
	\les
\Vert \RR_t \Vert_{L^2_t L^6_x}
	\Vert \nabla R\Vert_{L^\infty_t L^3_x}
	+
	\Vert \RR_t\Vert_{L^2_t H^1_x}
	+
	\Vert \RR\Vert_{L^2_t H^s_x}
	\les
	1
	.
	\label{EQ216}
	\end{split}
\end{align}

For the space component of the third term on the right side of \eqref{EQ111}, we use the trace inequality to obtain
\begin{align}
	\begin{split}
	\Vert \RR^{-1} \Vert_{L^2_t H_x^{s-1/2}(\Gammac)}
	\les		
	\Vert \RR^{-1} \Vert_{L^2_t H^s_x}
	\les
	1
	,
	\label{EQ189}
	\end{split}
\end{align}
where the last inequality follows from Lemma~\ref{L04}.
For the time component, we proceed analogously to \eqref{EQ216}, obtaining
\begin{align}
	\begin{split}
	\Vert \RR^{-1} \Vert_{H^{s/2-1/4}_t L^2(\Gammac)}
	&
	\les
	\Vert \RR^{-1} \Vert_{H^1_t H^1_x}	
	\les
	1
         ,
	\label{EQ190}
	\end{split}
\end{align}
since~$s\leq 5/2$.

For the last term on the right side of \eqref{EQ111}, we proceed analogously as in \eqref{EQ214}, obtaining
\begin{align}
	\begin{split}
		\Vert \RR_0^{-1} \nabla \RR_0\Vert_{H^{s-2}}
		\les
		\Vert \RR_0^{-1} \Vert_{H^s} 
		\Vert\nabla \RR_0\Vert_{H^{s-2}}	
		\les
		1.	
		\label{EQ1113}
	\end{split}
\end{align}

Combining \eqref{EQ111}, \eqref{EQ370}, and \eqref{EQ973}--\eqref{EQ1113}, we arrive at
\begin{align}
	\begin{split}
		\Vert \bar{v} \Vert_{K^{s+1}}
	&\les
	(\epsilon 
	+ 
	\tilde{\epsilon}C_\epsilon
	+
	C_{\epsilon, \tilde{\epsilon}}
	\TT^{1/2}) M		
	+
	C_\epsilon
	,
	\llabel{EQ119}
	\end{split}
\end{align}
for any $\epsilon, \tilde{\epsilon} \in (0,1]$.
Taking appropriate $\epsilon$, $\tilde{\epsilon}$, and $\TT>0$
(first $\epsilon$ sufficiently small, then $\tilde \epsilon$ sufficiently small depending on $\epsilon$, and then $\TT$ sufficiently small, depending on $\epsilon$ and~$\tilde \epsilon$), we get
\begin{align}
	\begin{split}
	\Vert \bar{v}\Vert_{K^{s+1}}
	\leq
	M
	,
	\label{EQ245}
	\end{split}
\end{align}
by allowing $M\geq 1$ sufficiently large.

Thus, we have shown that the mapping $\Lambda \colon v\mapsto \bar{v}$ is well-defined from $\mathcal{Z}$ to $\mathcal{Z}$ and satisfies \eqref{EQ245} for some $M\geq 1$, which depends on the size of the initial data and nonhomogeneous terms.

\subsection{Contracting property}
\label{sec052}
In this section, we prove
\begin{align}
	\begin{split}
	\Vert \Lambda(v_1) - \Lambda(v_2) \Vert_{K^{s+1}}
	\leq
	\frac{1}{2} \Vert v_1 - v_2 \Vert_{K^{s+1}}
	\comma
	v_1, v_2 \in \mathcal{Z}
	,
	\label{EQ121}
	\end{split}
\end{align}
where $M\geq 1$ is fixed as in \eqref{EQ245} and $\TT >0$ is a sufficiently small constant as in the previous section, which is further restricted below.
We emphasize that the implicit constants below are allowed to depend on~$M$.

\begin{proof}[Proof of Theorem~\ref{T03}]
Let $v_1, v_2 \in \mathcal{Z}$. Let $(\RR_1, \xi_1, \xi_{1t}, \bar{v}_1)$ and $(\RR_2, \xi_2, \xi_{2t}, \bar{v}_2)$ be the corresponding solutions of \eqref{EQ146}--\eqref{EQ147}, \eqref{EQ72}--\eqref{EQ73}, and \eqref{EQ124}--\eqref{EQ125} with the same initial data $(\RR_0, w_0, w_1, v_0)$ and the same nonhomogeneous terms $(f,h)$.
We denote $\tilde{\VV} = \bar{v}_1 - \bar{v}_2$,
$\tilde{v}= v_1 - v_2$, $\tilde{\RR} = \RR_1-\RR_2$, and $\tilde{\xi} = \xi_1 - \xi_2$. 
The difference $\tilde{V}$ satisfies
\begin{align}
	\begin{split}
	\tilde{V}_t
	- 
	\lambda \RR_1 \dive 
	(\nabla \tilde{V}
	+ 
	(\nabla \tilde{V})^T)
	-
	\mu \RR_1 \nabla \dive \tilde{V}
	= 
	g
	\inon{in~$[0,1]\times \Omegaf$}
	,
	\end{split}
   \llabel{EQ155}
\end{align}
with the boundary conditions and the initial data
\begin{align*}
	&
	\lambda 
	(\partial_k \tilde{\VV}_j 
	+ 
	\partial_j \tilde{\VV}_k)
	\nu^k
	+
	\mu \partial_k \tilde{\VV}_k  \nu^j
	=
	\partial_k \tilde{\xi}_j \nu^k
	-
	\RR_1^{-1} \RR_2^{-1} \tilde{\RR} \nu^j
	\inon{on~$[0,1]\times\Gammac$}
	,
	\\&
	\tilde{\VV}
	=
	0
	\inon{on~$[0,1]\times \Gammaf$}
	,
	\\&
	\tilde{V}~~\text{periodic in the $y_1$~and~$y_2$ directions},
	\\&
	\tilde{\VV} (0) = 0
	\inon{in~$\Omegaf$}
	,
   \llabel{EQ156}
\end{align*}
for $j= 1,2,3$, where
\begin{align}
	\begin{split}
		g
		=	
		-
		\RR_1 \nabla \RR_1^{-1}
		+
		\RR_2 \nabla \RR_2^{-1}	
		+
		\lambda \tilde{R} \dive (\nabla \bar{v}_2 + (\nabla \bar{v}_2)^T)
		+
		\mu \tilde{R} \nabla \dive \bar{v}_2
		.
		\label{EQ128}
	\end{split}
\end{align}
Note that 
$
	g(0)=0
$.
We proceed as in \eqref{EQ111} to obtain
\begin{align}
	\begin{split}
	\Vert \tilde{\VV}\Vert_{K^{s+1}}	
	&
	\les
	\left\Vert \frac{\partial \tilde{\xi}}{\partial \nu} \right\Vert_{K^{s-1/2 }_{\Gammac}}
	+
	\Vert \tilde{\RR} \RR_1^{-1} \RR_2^{-1} \nabla \RR_2 \Vert_{K^{s-1} }
	+
	\Vert \RR_1^{-1} \nabla \tilde{R} \Vert_{K^{s-1}}
	+
	\Vert \tilde{R} D_x^2 \bar{v}_2 \Vert_{K^{s-1}}
	\\&\indeq
	+
	\Vert \RR_1^{-1} \RR_2^{-1} \tilde{\RR} \Vert_{K^{s-1/2}_{\Gammac}}
	,
	\label{EQ127}
	\end{split}
\end{align}
where the last inequality follows from~\eqref{EQ128}
and
$  -
		R_1 \nabla R_1^{-1}
		+
		R_2 \nabla R_2^{-1}
  =
  R_1^{-1} \nabla\tilde R
  - R_{1}^{-1} R_{2}^{-1} \tilde R\nabla R_2
$.
The difference $\tilde{\xi}$ satisfies the wave equation
\begin{align}
	\tilde{\xi}_{tt} 
	- 
	\Delta \tilde{\xi}
	= 
	0
	\inon{in~$[0,1] \times \Omegae$}
	,
   \llabel{EQ157}
\end{align}
with the boundary condition and the initial data
\begin{align}
	\begin{split}
		&
	\tilde{\xi}(t,x) 
	= 
	\int_0^t \phi_{\TT} 
	\tilde{v} (\tau, x) \,d\tau
	\inon{on~$[0, 1]\times \Gammac$}
	,
	\\
	&
	(\tilde{\xi}, \tilde{\xi}_t) (0,x) 
	= 
	(0, 0)
	\inon{in~$\Omegae$}
	.
	\llabel{EQ159}	
	\end{split}
\end{align}  
For the first term on the right side of \eqref{EQ127}, we proceed as in \eqref{EQ115}--\eqref{EQ973} to obtain
\begin{align}
	\begin{split}
	\left\Vert \frac{\partial \tilde{\xi}}{\partial \nu}\right\Vert_{K^{s-1/2 } (\Gammac)}
	\les
(\epsilon 
+ 
\tilde{\epsilon}C_\epsilon
+
C_{\epsilon, \tilde{\epsilon}}
\TT^{1/2})
	\Vert \tilde{v}\Vert_{K^{s+1}}
	,
	\label{EQ144}
	\end{split}
\end{align}
for any $\epsilon, \tilde{\epsilon} \in (0,1]$.

Since the difference $\tilde{R}$ satisfies the ODE system
\begin{align}
	&
	\tilde{\RR}_t	
	-
	\phi_{\TT}
	\tilde{\RR} \dive v_2
	=
	\RR_1 \phi_{\TT}
	\dive \tilde{v}
	\inon{in~$[0, 1]\times \Omegaf$}
	,
	\label{EQ174}
	\\
	&
	\tilde{\RR}(0)
	=
	0
	\inon{in~$\Omegaf$}
	,
	\label{EQ136}
\end{align}
the solution is given by
\begin{align}
	\begin{split}
	\tilde{\RR}(t,x)
	&=
	\int_0^t e^{\int_\tau^t \phi_{\TT} \dive v_2} \phi_{\TT} (\tau) \RR_1(\tau) \dive \tilde{v}(\tau) \,d\tau
	\inon{in~$[0,1]\times \Omegaf$ }.
	\label{EQ173}
	\end{split}
\end{align}
For the second term on the right side of \eqref{EQ127}, we obtain
\begin{align}
	\begin{split}
	\Vert\tilde{\RR}  \RR_1^{-1} \RR_2^{-1} \nabla \RR_2 \Vert_{L^2_t H^{s-1}_x}
	\les
	\Vert \tilde{\RR} \Vert_{L^2_t H^s_x} 
	\Vert  \RR_1^{-1} \Vert_{L^\infty_t H^s_x}
	\Vert  \RR_2^{-1} \Vert_{L^\infty_t H^s_x}
	\Vert \RR_2 \Vert_{L^\infty_t H^{s}_x}
	\les
	\Vert \tilde{\RR} \Vert_{L^\infty_t H^s_x} 
	,
	\label{EQ134}
	\end{split}
\end{align}
where we used H\"older's inequality and Lemma~\ref{L04}, and then
from \eqref{EQ173} it follows that
\begin{align}
	\begin{split}
	\Vert\tilde{\RR}  \RR_1^{-1} \RR_2^{-1} \nabla \RR_2 \Vert_{L^2_t H^{s-1}_x}
	\les	
	\Vert \tilde{\RR} \Vert_{L^\infty_t H^s_x} 
	\les
	\TT^{1/2} \Vert \tilde{v} \Vert_{L^2_t H^{s+1}_x},
	\end{split}
   \label{EQ80}
\end{align}
where we used the Cauchy-Schwarz inequality.
For the time component 
(note that $s/2-1/4\leq 1$),
we have
\begin{align}
	\begin{split}
	&
	\Vert (\tilde{\RR}  \RR_1^{-1} \RR_2^{-1} \nabla \RR_2)_t \Vert_{L^2_t L^2_x}	
	\\&\indeq
	\les
	\Vert \tilde{\RR}_t \nabla \RR_2 \Vert_{L^2_t L^2_x}
	+
	\Vert \tilde{\RR} \RR_{1t} \nabla \RR_2 \Vert_{L^2_t L^2_x}
	+
	\Vert \tilde{\RR} \RR_{2t} \nabla \RR_2 \Vert_{L^2_t L^2_x}
	+
	\Vert \tilde{\RR} \nabla \RR_{2t} \Vert_{L^2_t L^2_x}
	\\&\indeq
	\les
	\Vert \tilde{R} \Vert_{L^\infty_t L^\infty_x} 
	\Vert \dive v_2 \Vert_{L^2_t L^4_x} 
	\Vert \nabla \RR_2 \Vert_{L^\infty_t L^4_x}
	+
	\Vert \RR_1 \Vert_{L^\infty_t L^\infty_x}
	\Vert \dive \tilde{v}\Vert_{L^2  ((0, 2\TT), L^\infty (\Omegaf))}
	\Vert \nabla \RR_2 \Vert_{L^\infty_t L^2_x}
	\\&\indeq\indeq
	+
	\Vert \tilde{\RR} \Vert_{L^\infty_t L^\infty_x}
	\Vert \dive v_1 \Vert_{L^2_t L^4_x}
	\Vert \nabla \RR_2 \Vert_{L^\infty_t L^4_x}
		+
	\Vert \tilde{\RR} \Vert_{L^\infty_t L^\infty_x}
	\Vert \dive v_2 \Vert_{L^2_t L^4_x}
	\Vert \nabla \RR_2 \Vert_{L^\infty_t L^4_x}
	\\&\indeq\indeq
	+
	\Vert \tilde{\RR} \Vert_{L^\infty_t L^\infty_x}
	\Vert \nabla \RR_{2t}\Vert_{L^2_t L^2_x}
        \\&\indeq
	\les
	(
	\epsilon	+
	C_\epsilon \TT^{1/2}) \Vert \tilde{v} \Vert_{K^{s+1}} 
	,
	\label{EQ135}
	\end{split}
\end{align}
for any $\epsilon \in (0,1]$, where we used H\"older's inequality, Lemma~\ref{L04}, \eqref{EQ174}, and~\eqref{EQ173},
as well as
$	\Vert \tilde{\RR} \Vert_{L^\infty_t H^s_x} 
	\les
	\TT^{1/2} \Vert \tilde{v} \Vert_{L^2_t H^{s+1}_x}
$.
Note that
$\Vert \tilde{\RR}  \RR_1^{-1} \RR_2^{-1} \nabla \RR_2 \Vert_{L^2_t L^2_x}$
does not need to be estimated since it is dominated by~\eqref{EQ80}.
Similarly, the third term on the right side of \eqref{EQ127} is estimated as
\begin{align}
	\begin{split}
	\Vert \RR_1^{-1} \nabla \tilde{\RR} \Vert_{L^2_t H^{s-1}_x}
	\les
	\Vert \RR_1^{-1} \Vert_{L^\infty_t H^s_x}
	\Vert \tilde{\RR} \Vert_{L^\infty_t H^s_x}
	\les
	\TT^{1/2} \Vert \tilde{v}\Vert_{K^{s+1}}
	,
	\label{EQ138}
	\end{split}
\end{align}
and for the time component
\begin{align}
	\begin{split}
	\Vert (\RR_1 \nabla \tilde{\RR})_t \Vert_{L^2_t L^2_x}
	&\les
	\Vert \RR_{1t} \nabla \tilde{\RR} \Vert_{L^2_t L^2_x}	
	+
	\Vert \RR_1 \nabla \tilde{\RR}_t  \Vert_{L^2_t L^2_x}	
	\\&
	\les
	\Vert \RR_{1t} \Vert_{L^2_t L^\infty_x}
	\Vert \nabla \tilde{\RR} \Vert_{L^\infty_t L^2_x}
	+
		\Vert \RR_{1} \Vert_{L^\infty_t L^\infty_x}
	\Vert \nabla \tilde{\RR}_t \Vert_{L^2_t L^2_x}
	\les
	(\epsilon
	 +
	 C_\epsilon \TT^{1/2})
	 \Vert \tilde{v} \Vert_{K^{s+1}}
	 .
	 \label{EQ139}
	\end{split}
\end{align}
Again, the term $\Vert \RR_1 \nabla \tilde{\RR} \Vert_{L^2_t L^2_x}$ is dominated by~\eqref{EQ138}.
Regarding the fourth term on the right side of \eqref{EQ127}, we use Corollary~\ref{C02} to obtain
\begin{align}
	\begin{split}
	\Vert \tilde{\RR} D_x^2 \bar{v}_2\Vert_{L^2_t H^{s-1}_x}
	\les
	\Vert \tilde{\RR}\Vert_{L^\infty_t H^s_x}
	\Vert \bar{v}_2 \Vert_{L^2_t H^{s+1}_x}
	\les
	\TT^{1/2} 
	\Vert \tilde{v} \Vert_{K^{s+1}}
	.
	\label{EQ141}
	\end{split}
\end{align}
To treat
$\Vert \tilde{\RR} D_x^2 \bar{v}_2  \Vert_{H^{(s-1)/2}_t L^2_x}$, we 
claim that for any $\alpha>1/2$ and $\delta>0$, we have
  \begin{align}
  \begin{split}
   \Vert A B  \Vert_{H^{\alpha}_t L^2_x}	   
   &\les
       \Vert A   \Vert_{H^{\alpha}_t H^{3/2+\delta}_x}
   \Vert B   \Vert_{H^{\alpha}_t L_{x}^2}
  \end{split}
   \label{EQ100}
  \end{align}
on the domain $(0,1)\times\Omegaf$.
Using extensions, we may assume that the domain is actually $\mathbb{R}\times\mathbb{R}^{3}$.
Then
\begin{align}
  \begin{split}
   \Vert A B  \Vert_{H^{\alpha}_t L^2_x}
   &=    \Vert A B  \Vert_{L_{x}^{2}H^{\alpha}_t}
   \lec   \bigl\Vert
             \Vert A   \Vert_{H^{\alpha}_t}
             \Vert B   \Vert_{H^{\alpha}_t}
          \bigr\Vert_{L^2_{x}}
  \\&
   \lec
      \Vert A   \Vert_{L_{x}^{\infty}H^{\alpha}_t}
      \Vert B   \Vert_{L_{x}^2 H^{\alpha}_t}
   \lec
    \Vert A   \Vert_{H_x^{3/2+\delta} H^{\alpha}_t}
    \Vert B   \Vert_{L_{x}^2 H^{\alpha}_t}
   \\&
   =
    \Vert A   \Vert_{H^{\alpha}_t H_x^{3/2+\delta}}
    \Vert B   \Vert_{H^{\alpha}_t L_{x}^2}
   ,
  \end{split}
   \label{EQ98}
  \end{align}
since $\alpha >1/2$, and \eqref{EQ100} follows.
Using \eqref{EQ100}, we then write
  \begin{align}
  \begin{split}
  \Vert \tilde{\RR} D_x^2 \bar{v}_2  \Vert_{H^{(s-1)/2}_t L^2_x}
  &\lec
   \Vert \tilde R   \Vert_{H^{1}_t H_x^{2}}
   \Vert D_x^2 \bar v_2   \Vert_{H^{(s-1)/2}_t L_{x}^2}
   \lec
     \Vert \tilde R_t   \Vert_{L^2_t H_x^{2}}
     +       \Vert \tilde R   \Vert_{L^2_t H_x^{2}}
   ,
  \end{split}
   \label{EQ101}
  \end{align}
where we used Corollary~\ref{C02} in the last inequality.
From \eqref{EQ203} and \eqref{EQ174}, it follows that
\begin{align}
	\begin{split}
	\Vert \tilde{R}_t \Vert_{L^2_t H^2_x}
	&
	\les
	\Vert \tilde{R}\Vert_{L^\infty_t H^2_x}	
	\Vert v_2\Vert_{L^2_t H^{s+1}_x}
	+
	\Vert R_1\Vert_{L^\infty_t H^2_x}
	\Vert 	\phi_{\TT} \tilde{v}\Vert_{L^2_t H^3_x}
	\\&
	\les
	\TT^{1/2} \Vert
 \tilde{v}\Vert_{L^2_t H^{s+1}_x}
	+
	\epsilon 
	\Vert \tilde{v}\Vert_{L^2_t H^{s+1}_x}
	+
	C_\epsilon
	\TT^{1/2}
		\Vert \tilde{v}\Vert_{H^{(s+1)/2}_t L^2_x}
	\les
	(\epsilon + C_\epsilon \TT^{1/2})
	\Vert \tilde{v}\Vert_{K^{s+1}}
	,
	\label{EQ012}
	\end{split}
\end{align}
for any $\epsilon \in (0,1]$, since~$s>2$.
Combining \eqref{EQ101} and \eqref{EQ012}, we arrive at
\begin{align}
	\begin{split}
	\Vert \tilde{R} D^2_x \bar{v}_2\Vert_{H^{(s-1)/2}_t L^2_x}
	\les	
	(\epsilon + C_\epsilon \TT^{1/2})
	\Vert \tilde{v}\Vert_{K^{s+1}}.
	\end{split}
\end{align}
For the last term on the right side of \eqref{EQ127}, we use the trace inequality and arrive at
\begin{align}
	\begin{split}
	\Vert \RR_1^{-1} \RR_2^{-1} \tilde{\RR} \Vert_{L^2_t H_x^{s-1/2}(\Gammac)}
	\les	
	\Vert \tilde{\RR} \Vert_{L^2_t H^{s}_x}
	\les
	\TT^{1/2} \Vert \tilde{v}\Vert_{K^{s+1}}
	\label{EQ192}
	\end{split}
\end{align}
and
\begin{align}
	\begin{split}
	\Vert \RR_1^{-1} \RR_2^{-1} \tilde{\RR} \Vert_{H^{s/2-1/4}_t L_x^2(\Gammac)}
	\les
	\Vert \RR_1^{-1} \RR_2^{-1} \tilde{\RR} \Vert_{H^1_t H^1_x}
	\les
	\Vert (\RR_1^{-1} \RR_2^{-1} \tilde{\RR})_t \Vert_{L^2_t H^1_x}
	+
	\Vert \RR_1^{-1} \RR_2^{-1} \tilde{\RR} \Vert_{L^2_t H^1_x}
	,
	\label{EQ191}
	\end{split}
\end{align}
since $s\leq 5/2$.
For the first term on the right side of \eqref{EQ191}, we proceed as in \eqref{EQ135} to obtain
\begin{align}
	\begin{split}
	\Vert (\RR_1^{-1} \RR_2^{-1} \tilde{\RR})_t \Vert_{L^2_t H^1_x}
	&
	\les
	\Vert \RR_{1t} \tilde{\RR} \Vert_{L^2_t H^1_x}	
	+
	\Vert \RR_{2t} \tilde{\RR} \Vert_{L^2_t H^1_x}	
	+
	\Vert \tilde{\RR}_t \Vert_{L^2_t H^1_x}
	\les
	( \epsilon+
	C_{\epsilon} \TT^{1/2} )
	\Vert \tilde{v} \Vert_{K^{s+1}}
	.
	\end{split}
\end{align}
The second term on the right side of \eqref{EQ191} is estimated analogously to~\eqref{EQ192}, and we get
\begin{align}
	\begin{split}
		\Vert \RR_1^{-1} \RR_2^{-1} \tilde{\RR} \Vert_{L^2_t H^1_x}
		\les
		\TT		\Vert \tilde{v} \Vert_{K^{s+1}}
		.
	\label{EQ193}
	\end{split}
\end{align}

Applying the above estimates in \eqref{EQ127}, we obtain
\begin{align}
	\begin{split}
	\Vert \tilde{\VV} \Vert_{K^{s+1}}
	\les
	(\epsilon+ \tilde{\epsilon}
	+ 
	\TT^{1/2} C_{\tilde{\epsilon}, \epsilon}  ) \Vert \tilde{v}\Vert_{K^{s+1}}
	,
	\end{split}
   \llabel{EQ161}
\end{align}
for any $\epsilon, \tilde{\epsilon} \in (0,1]$.
Taking appropriate $\epsilon$, $\tilde{\epsilon}$, and $\TT>0$
(first $\epsilon$ sufficiently small, then $\tilde \epsilon$ sufficiently small depending on $\epsilon$, and then $\TT$ sufficiently small, depending on $\epsilon$ and~$\tilde \epsilon$), we conclude the proof of~\eqref{EQ121}.
Thus, the mapping $\Lambda$ is a contraction from $\mathcal{Z}$ to $\mathcal{Z}$. 
Using the Banach fixed point theorem, there exists a unique solution $v \in \mathcal{Z}$ such that $\Lambda (v) = v$ and which also satisfies \eqref{EQ245} for some~$M\geq 1$.

Now we fix the constant $\TT>0$ as above. 
Using Lemma~\ref{L03}, we have the interior regularity estimate 
\begin{align}
	\begin{split}
		&
		\Vert w \Vert_{C([0, 1], H^{s+1/4 -\epsilon_0} (\Omegae))}
		+
		\Vert w_t \Vert_{C([0, 1], H^{s-3/4-\epsilon_0}(\Omegae))}
		\\&\indeq	
		\les
		\Vert w_0\Vert_{H^{s+1/4-\epsilon_0} (\Omegae)}
		+
		\Vert w_1\Vert_{H^{s-3/4 -\epsilon_0} (\Omegae)}
		+
		\Vert w\Vert_{H^{s+1/4 -\epsilon_0, s+1/4 -\epsilon_0} ( \Gammac)}
		.
		\label{EQ180}
	\end{split}
\end{align}
For the last term on the right side, we appeal to \eqref{EQ007}, yielding
\begin{align}
	\begin{split}
		\Vert w\Vert_{H^{s+1/4 -\epsilon_0, s+1/4-\epsilon_0} (\Gammac)}
		&
		\les
		\Vert w_t \Vert_{H^{s-3/4-\epsilon_0}_t L^2_x ( \Gammac)}
		+
		\Vert w\Vert_{L^2_t H_x^{s+1/4 -
				\epsilon_0}( \Gammac)}
		\\&
		\les
		\Vert \phi_{\TT} v\Vert_{H^{s-3/4-\epsilon_0}_t L^2_x ( \Gammac)}
		+
		\Vert (1-\phi_{\TT}) v_0 \Vert_{H^{s-3/4-\epsilon_0}_t L^2_x ( \Gammac)}
		+
		\Vert w \Vert_{L^2_t H_x^{s+1/2}(\Gammac)}.
		\label{EQ181}
	\end{split}
\end{align}
For the first term on the far right side of \eqref{EQ181}, we appeal to Corollary~\ref{C01} and Sobolev inequality to get
\begin{align}
	\begin{split}
		\Vert \phi_{\TT} v\Vert_{H^{s-3/4-\epsilon_0}_t L^2_x ( \Gammac)}
		&
	\les
	\Vert v\Vert_{H^{s-3/4-\epsilon_0}_t L^2_x ( \Gammac)}
	\les
	\Vert  v\Vert_{H^{(s+1)/2}_t L^2_x}
	+
		\Vert  v\Vert_{L^2_t H^{s+1}_x}
		,
		\label{EQ010}
	\end{split}
\end{align}
since $s\leq 2+2\epsilon_0$.
From \eqref{EQ350} and \eqref{EQ180}--\eqref{EQ010}, it follows that
\begin{align}
	\begin{split}
		&
		\Vert w \Vert_{C([0, 1], H^{s+1/4 -\epsilon_0} (\Omegae))}
		+
		\Vert w_t \Vert_{C([0, 1], H^{s-3/4-\epsilon_0}(\Omegae))}
		\leq
		C,
		\label{EQ116}
	\end{split}
\end{align}
where $C>0$ is a constant.
By \eqref{EQ116} and Lemma~\ref{L04}, there exists a unique solution 
\begin{align}
	\begin{split}
		(v, \RR, w, w_t)
		&
		\in
		K^{s+1} ((0,\TT) \times \Omegaf)
		\times
		H^1 ((0, \TT), H^s(\Omegaf))
		\\&\indeqtimes
		C([0, \TT], H^{s+1/4-\epsilon_0}(\Omegae))
		\times
		C([0, \TT], H^{s-3/4-\epsilon_0}(\Omegae))
	\end{split}
   \llabel{EQ165}
\end{align}
to the system \eqref{EQ71}--\eqref{EQ76}, with the corresponding norms bounded by a function of the initial data and the nonhomogeneous terms.
\end{proof}

\begin{Remark}
\label{pressure}
{\rm
As pointed out at the end of Section~\ref{sec_setting},
the approach extends to more general pressure laws.
For general equation of state $q(r)$, we assume that $q(r)$ is smooth such that $q(0)=0$ and $q(r_1) - q(r_2) = (r_1 -r_2) \tilde{q}(r_1, r_2)$ for any $r_1$ and $r_2$, where $\tilde{q}$ is a smooth function.
We shall briefly outline the modifications needed for this general pressure law.
In Section~\ref{sec051}, we have $\Vert R \nabla (q(R^{-1})) \Vert_{K^{s-1}}$ instead of the second term on the right side of~\eqref{EQ111}.  For the space component, we use the H\"older and Sobolev inequalities to get
\begin{align*}
	\begin{split}
		\Vert R \nabla (q(R^{-1}))\Vert_{L^2_t H^{s-1}_x}
		&
		\les
		\Vert  q'(R^{-1}) R^{-1} \nabla R \Vert_{L^2_t H^{s-1}_x}
		\les
		\Vert q'(R^{-1}) \Vert_{L^\infty_t H^s_x}
		\Vert R^{-1} \nabla R \Vert_{L^2_t H^{s-1}_x}	
		\les 1,
	\end{split}
\end{align*}
where the last inequality follows from~\eqref{EQ214}.
For the time component, we appeal to \eqref{EQ98}, yielding
\begin{align*}
	\begin{split}
	\Vert R \nabla (q(R^{-1}))\Vert_{H^{(s-1)/2}_t L^2_x}
	&
	\les
	\Vert q'(R^{-1}) \Vert_{H^{1}_t H^{3/2+\delta}_x}
	\Vert R^{-1} \nabla R \Vert_{H^{(s-1)/2}_t L^2_x}
	\les 1
	\end{split}
\end{align*}
where we used Lemma~\ref{L04} and \eqref{EQ216} in the last inequality.
The third term on the right side of \eqref{EQ111} is replaced by $\Vert q(R^{-1})\Vert_{K^{s-1/2}_{\Gammac}}$, which can be estimated in a similar fashion.
In Section~\ref{sec052}, the first two terms on the right side of \eqref{EQ128} are replaced by $-R_1 \nabla (q( R_1^{-1} ))+ R_2 \nabla (q( R_2^{-1} ))$ and the $K^{s-1}$ norm can be estimated using the structural assumption on $q(r)$.
\square
}
\end{Remark}

\startnewsection{Solution to the Navier-Stokes-wave system}{sec06}
In this section, we provide the local existence for the coupled Navier-Stokes-wave system \eqref{EQ260}--\eqref{EQ23} with the boundary conditions \eqref{EQ262}--\eqref{EQ267} and the initial data~\eqref{EQ265}.
Let $v\in \mathcal{Z}$ where $\mathcal{Z}$ is as in \eqref{EQ188}, with constant $M\geq 1$ to be determined below. 
Let $\phi_{\TT} (t)$ be a smooth cutoff function as defined in Section~\ref{sec05};
here, $\TT \in (0, 1/4)$ is a constant to be determined below; it is assumed to be smaller
than the constant $\TT$ from the previous section, which we from here
on denote by~$\TT_0$.
We allow all constants to depend on $\TT_0$ (but not on $\TT$).

We again modify the system to be able to construct a solution on a unit time interval.
Let
\begin{align}
	\eta (t,x)
	=
	x
	+
	\int_0^t \phi_{\TT}(\tau) v(\tau, x) \,d\tau
	\inon{in~$[0,1]\times \Omegaf$}
	\label{EQ999}
\end{align}
be a modified Lagrangian flow map and
$a (t,x)= (\nabla \eta(t,x))^{-1}$ its inverse matrix, while
we denote by $\JJ(t,x) = \det (\nabla \eta(t,x))$ the corresponding Jacobian.
The density equations we consider is
\begin{align}
	&
	\RR_t - \RR \phi_{\TT} a_{kj} \partial_k v_j
	=
	0
	\inon{in~$[0,1]\times \Omegaf$}
	,
	\label{EQ940}
	\\
	&
	\RR(0) = \RR_0
	\inon{on~$\Omegaf$}
	,
	\label{EQ941}
\end{align}	
with the solution given by
\begin{align}
	\begin{split}
		\RR(t,x) 
		=
		\RR_0(x)  e^{\int_0^t \phi_{\TT} (\tau) a_{kj} (\tau, x) \partial_k v_j (\tau, x) \,d\tau}
		\inon{in~$[0, 1]\times \Omegaf$}
		.
	\end{split}
   \llabel{EQ212}
\end{align}
Next, we consider the solution $w$ to the wave equation \eqref{EQ72} with the boundary condition \eqref{EQ73}--\eqref{EQ76} and the initial data
$(w, w_t) (0) 
= 
(w_0, w_1)$
in~$\Omegae$.
With $(\eta, a, \JJ, \RR, w)$ constructed, we define
\begin{align}
	\Pi\colon v\in \mathcal{Z} \mapsto \bar{v}
	,
   \llabel{EQ166}
\end{align}
where $\bar{v}$ is the solution of the nonhomogeneous parabolic problem
\begin{align}
	&
	\partial_t \bar{v}_{j} 
	-
	\lambda \RR \partial_k 
	(\partial_j \bar{v}_k + \partial_k \bar{v}_j)
	-
	\mu \RR \partial_j \partial_k \bar{v}_k
	=
	f_j
	\inon{in~$[0,1]\times \Omegaf$}
	,
	\label{EQ175}
	\\
	&
	\lambda (\partial_k \bar{v}_j + \partial_j \bar{v}_k) \nu^k
	+
	\mu \partial_k \bar{v}_k \nu^j
	=
	\partial_k w_j \nu^k
	+
	h_j
	\inon{in~$[0,1] \times \Gammac$}
	,
	\label{EQ176}
		\\&
	\bar{v}~~\text{periodic in the $y_1$~and~$y_2$ directions},
	\label{EQ400}
	\\&
	\bar{v}(0)
	=
	v_0
	\inon{in~$\Omegaf$}
	;
	\label{EQ302}
\end{align}
in \eqref{EQ175}--\eqref{EQ176}, we set
\begin{align}
	\begin{split}
		f_j
		&
		=
		\lambda \RR \partial_k 
		(b_{mk} \partial_m \bar{v}_j
		+
		b_{mj} \partial_m \bar{v}_k)	
		+
		\lambda \RR b_{kl} 
		\partial_k
		(b_{ml} \partial_m \bar{v}_j
		+
		b_{mj} \partial_m \bar{v}_l)
		+
		\lambda \RR b_{kl} \partial_k 
		(\partial_l \bar{v}_j + \partial_j \bar{v}_l)
		\\&\indeq
		+
		\mu \RR \partial_j (b_{mi} \partial_m \bar{v}_i)
		+
		\mu \RR b_{kj} \partial_k (b_{mi} \partial_m \bar{v}_i)
		+
		\mu \RR b_{kj} \partial_k \partial_i \bar{v}_i
		-
		\RR b_{kj} \partial_k \RR^{-1}
		-
		\RR  \partial_j \RR^{-1}
		\\&
		=:
		I_1
		+
		I_2
		+
		I_3
		+
		I_4
		+
		I_5
		+
		I_6
		+
		I_7
		+
		I_8
		\label{EQ177}
	\end{split}
\end{align}
and 
\begin{align}
	\begin{split}	
		h_j
		&
		=
		\lambda (1-\JJ) (\partial_k \bar{v}_j + \partial_j \bar{v}_k) \nu^k
		+
		\mu (1-\JJ) \partial_k \bar{v}_k \nu^j
		-
		\lambda \JJ b_{kl} 
		(b_{ml} \partial_m \bar{v}_j 
		+ 
		b_{mj} \partial_m \bar{v}_l) \nu^k
		\\&\indeq
		+
		\JJ b_{kj} \RR^{-1} \nu^k
		+
		(\JJ-1) \RR^{-1} \nu^j
		-
		\lambda \JJ 
		(b_{mk} \partial_m \bar{v}_j + b_{mj} \partial_m \bar{v}_k) \nu^k
		-
		\lambda \JJ b_{kl} 
		(\partial_l \bar{v}_j 
		+
		\partial_j \bar{v}_l
		)\nu^k
		\\&\indeq
		-
		\mu J b_{kj} b_{mi} \partial_m \bar{v}_i \nu^k
		-
		\mu \JJ b_{mi} \partial_m \bar{v}_i \nu^j
		-
		\mu \JJ b_{kj} \partial_i \bar{v}_i \nu^k
		+
		\RR^{-1} \nu^j
		\\&
		=
		:
		K_1
		+
		K_2
		+
		K_3
		+
		K_4
		+
		K_5
		+
		K_6
		+
		K_7
		+
		K_8
		+
		K_9
		+
		K_{10}
		+
		K_{11}
		,
		\label{EQ178}
	\end{split}
\end{align}
for $j=1,2,3$, where
  \begin{equation}
  b = a - \mathbb{I}_3   
  ,
   \llabel{EQ106}
  \end{equation}
and $\mathbb{I}_3$ is the three-dimensional identity matrix.

Before we bound the terms in \eqref{EQ177}--\eqref{EQ178} and
prove the contracting property, as in Section~\ref{sec05}, we provide some necessary estimates on the variable coefficients.

\subsection{The Lagrangian flow map, Jacobian matrix, and density estimates}
We start with estimates on the Jacobian and the inverse matrix of the flow map.

\cole
\begin{Lemma}
	\label{L08}
Suppose that $\Vert v\Vert_{K^{s+1} ((0,1)\times \Omegaf)} \leq M$,
where $M\geq 1$, and let $\delta \in (0,1/5)$.
Then for $\TT>0$ sufficiently small depending on $M$ and $\delta$, the following statements hold:
\begin{enumerate}[label=(\roman*)]
	\item $\Vert b \Vert_{L^\infty_t H^s_x} +\Vert b\Vert_{H^1_t H^{3/2+\delta}_x}\les \TT^{1/20}$,
	\item $\Vert b\Vert_{H^1_t H^s_x} \les M$,
	\item $\Vert 1-\JJ \Vert_{L^\infty_t H^s_x} \les\TT^{1/20}$,
	\item $\Vert \JJ \Vert_{L^\infty_t L^\infty_x} +\Vert \JJ^{-1} \Vert_{L^\infty_t L^\infty_x}+ \Vert \JJ \Vert_{L^\infty_t H^s_x} 
	+ \Vert \JJ^{-1} \Vert_{L^\infty_t H^s_x} 
	\les 1$,
	\item $\Vert \JJ^{-1}\Vert_{H^1_t H^{3/2+\delta}_x} + \Vert \JJ\Vert_{H^1_t H^{3/2+\delta}_x}\les 1$,
	\item $ \Vert \JJ \Vert_{H^1_t H^s_x} \les M$,
\end{enumerate}
where the region of dependence is understood to be $(0,1)\times \Omegaf$.
\end{Lemma}
\colb
We emphasize that the implicit constants in the above inequalities (i)-(vi) are independent of~$M$ and~$\delta$.

\begin{proof}[Proof of Lemma~\ref{L08}]
(i) From \eqref{EQ204} and \eqref{EQ999} it follows that
\begin{align}
	b_t
	=
        -
	\phi_{\TT} (		
	b  \nabla v b
	+
	 b\nabla v
	+
	 \nabla v  b
	+
	\nabla v
	)
	\inon{in~$[0,1]\times \Omegaf$}
	,
	\label{EQ211}
\end{align}
while~$b(0) = 0$. 
By the Fundamental Theorem of Calculus, it follows that for $t\in (0,2\TT)$ we have
\begin{align}
	\begin{split}
	\Vert b(t) \Vert_{H^s} 
	&
	\les
	\int_0^{t} \Vert b \Vert_{H^s}^2 
	\Vert \nabla v \Vert_{H^s} \,d\tau
	+
	\int_0^{t} \Vert b \Vert_{H^s} 
	\Vert \nabla v\Vert_{H^s} \,d\tau
	+
	\int_0^{t} \Vert \nabla v\Vert_{H^s} \,d\tau
	\\&
	\les
	\int_0^{t} \Vert v\Vert_{H^{s+1}}
	(\Vert b\Vert_{H^s}^2 + \Vert b\Vert_{H^s})
	 \,d\tau
	+
	\TT^{1/2} M
	,
\end{split}
   \llabel{EQ167}
\end{align}
where we appealed to the Cauchy-Schwarz inequality in the last step.
Using Gronwall's inequality, we obtain
\begin{align}
	\Vert b\Vert_{L^\infty  ((0,2\TT), H^s (\Omegaf))}
	\les
	\TT^{1/2} M
	\les
	\TT^{1/20}
	,
   \llabel{EQ168}
\end{align}
for $\TT>0$ sufficiently small;
the choice of the power $1/20$ is apparent in \eqref{EQ303} below.
Since also $b_t = 0$ on $(2\TT, 1)$, we then infer that
\begin{align}
	\Vert b\Vert_{L^\infty_t H^s_x}
	\les
	\TT^{1/20}
	.
	\label{EQ280}
\end{align}
Applying \eqref{EQ280} in \eqref{EQ211} and using \eqref{EQ203}, we obtain
\begin{align}
	\begin{split}
	\Vert b_t\Vert_{L^2_t H^{3/2+\delta}_x}
	&
	\les
	\Vert v\Vert_{L^2  ((0,2\TT), H^{5/2+\delta} (\Omegaf))}
	\les
	\epsilon 	\Vert  v\Vert_{L^2_t H^{s+1}_x } 
	+ 
		\epsilon^{(5+2\delta)/(3+2\delta-2s)}
	\Vert v\Vert_{L^2  ((0,2\TT), L^2 (\Omegaf))}  
	\\&
		\les
	\epsilon 	M
	+ 
		\epsilon^{(5+2\delta)/(3+2\delta-2s)}
	\TT^{1/2}
	M
	,
	\end{split}
   \llabel{EQ169}
\end{align}
for any $\epsilon \in (0,1]$.
Letting $\epsilon = \TT^{1/20} M^{-1}$, we get
\begin{align}
	\Vert b_t\Vert_{L^2_t H^{3/2+\delta}_x}
	\les
	\TT^{1/20}
	+
	\TT^{1/2+(5+2\delta)/20(3+2\delta-2s)} M^{1+(5+2\delta)/(2s-3-2\delta)}
	\les
	\TT^{1/20}
	,
		\label{EQ303}
\end{align}
for $\TT>0$ sufficiently small.
Combining \eqref{EQ280}--\eqref{EQ303}, we conclude the proof of (i).

(ii) From \eqref{EQ211} and H\"older's inequality it follows that
\begin{align}
	\begin{split}
		\Vert b_t\Vert_{L^2_t H^s_x}
		&
		\les
		\Vert \nabla v\Vert_{L^2_t H^s_x }
		\Vert b \Vert_{L^\infty_t H^s_x}^2	
		+
		\Vert \nabla v\Vert_{L^2_t H^s_x}
		\Vert b \Vert_{L^\infty_t H^s_x}
		+
		\Vert \nabla v\Vert_{L^2_t H^s_x }
		\les
		M
		,
	\end{split}
   \llabel{EQ170}
\end{align}
which gives~(ii).

(iii) From \eqref{EQ210} and \eqref{EQ999} we infer that $\JJ$ satisfies the ODE system 
\begin{align}
  \begin{split}
	&
	J_t
	=
	\phi_{\TT} 
	\JJ  a_{kj}  \partial_k v_j 
	\inon{in~$[0,1] \times \Omegaf$}
	,
	\\&
	\JJ(0) = 1
	\inon{in~$\Omegaf$}
	.
  \end{split}
	\label{EQ910}
\end{align}
The solution is given by
\begin{align}
	\begin{split}
		\JJ (t,x)
		= 
		e^{\int_0^t \phi_{\TT} (\tau) a_{kj} (\tau, x)\partial_k v_j (\tau, x) \,d\tau}	
		\inon{in~$[0,1] \times \Omegaf$}.
		\llabel{EQ207}
	\end{split}
\end{align}
Using the nonlinear Sobolev estimate, we have
\begin{align}
	\begin{split}
	\Vert \JJ -1 \Vert_{L^\infty_t H^s_x}
	\les
	C^{\TT^{1/2} M} -1
	\les
	\TT^{1/20}
	,
	\end{split}
   \llabel{EQ171}
\end{align}
for $\TT>0$ sufficiently small.

The proofs of (iv), (v), and (vi) are  analogous to the proof of Lemma~\ref{L04}, and thus we omit the details.
\end{proof}

The following lemma provides the necessary a~priori density estimates.
\cole
\begin{Lemma}
	\label{L09}
Assume that
\begin{align}
	(\RR_0, \RR_0^{-1}, b)
	\in
	H^s(\Omegaf)
	\times
	H^s(\Omegaf)
	\times
	L^\infty ((0,T), H^s (\Omegaf))
   \llabel{EQ172}
\end{align}
and
$\Vert v\Vert_{K^{s+1} ((0,1)\times \Omegaf)} \leq M$, where~$M\geq 1$.
Let $\delta \in (0,1/2)$.
Then for $\TT >0$ sufficiently small depending on $M$ and $\delta$, the solution to the ODE system \eqref{EQ940}--\eqref{EQ941} satisfies
\begin{enumerate}[label=(\roman*)]
	\item $\Vert \RR \Vert_{L^\infty_t L^\infty_x} +\Vert \RR^{-1} \Vert_{L^\infty_t L^\infty_x} +\Vert \RR \Vert_{L^\infty_t H^s_x} +\Vert \RR^{-1} \Vert_{L^\infty_t H^s_x} \les 1$, 
	\item $\Vert \RR^{-1} \Vert_{H^1_t H^{3/2+\delta}_x}+
	\Vert \RR \Vert_{H^1_t H^{3/2+\delta}_x}
	\les 1$,
	\item $\Vert \RR\Vert_{H^1_t H^s_x} \les M$,
\end{enumerate}
where the norm of dependence is $(0, 1)\times \Omegaf$.
\end{Lemma}
\colb

We emphasize that the implicit constants in the above inequalities (i)--(iii) are independent of~$M$ and~$\delta$.
The proof of Lemma~\ref{L09} is analogous to the proof of Lemma~\ref{L04}. 
Thus we omit the details.

\subsection{Uniform boundedness of the iterative sequence}
In this section we shall prove that the mapping $\Pi$ is well-defined from $\mathcal{Z}$ to $\mathcal{Z}$, for a sufficiently large constant $M \geq 1$ and a sufficiently small constant~$\TT>0$.
From Lemmas~\ref{L02} and~\ref{L09}, it follows that
\begin{align}
	\begin{split}
	\Vert \bar{v}\Vert_{K^{s+1} }
	&
	\les	
	\left
	\Vert \frac{\partial w}{\partial \nu} \right\Vert_{K^{s-1/2  }_{\Gammac}}
		+
	\Vert f\Vert_{K^{s-1} }
	+
	\Vert  h\Vert_{K^{s-1/2}_{\Gammac}}
	+
	\Vert f(0) \Vert_{H^{s-2}}
	+
	\Vert v_0\Vert_{H^s}
		\\&\indeq
		+
	\Vert v_0|_{\Gammac}\Vert_{H^{s+1/2} (\Gammac)}
	+
	\Vert \partial_3 v_0|_{\Gammaf}\Vert_{H^{s-1/2} (\Gammaf)}
	,
	\label{EQ201}
	\end{split}
\end{align}
where $f$ and $h$ are as in~\eqref{EQ177}--\eqref{EQ178}. 
Here and below, the time and space domains in the norms are understood to be $(0,1)$ and $\Omegaf$, respectively, unless indicated otherwise.
We emphasize that the implicit constants in this section are independent of~$M$.

For the first term on the right side of \eqref{EQ201}, we proceed as in \eqref{EQ115}--\eqref{EQ973} to obtain
\begin{align}
	\begin{split}
		\left\Vert \frac{\partial w}{\partial \nu} 
		\right\Vert_{K^{s-1/2} (\Gammac)}
		\les 
		(\epsilon+\tilde{\epsilon} C_\epsilon
		+ 
		\TT^{1/2} C_{\tilde{\epsilon}, \epsilon} ) M
		+
		C_\epsilon
		,
		\label{EQ931}
		\end{split}
\end{align}
for any $\epsilon, \tilde{\epsilon} \in (0,1]$.
Next, we estimate the $K^{s-1}$ norm of the terms on the right side of \eqref{EQ177} for $j=1,2,3$.
For the space component of the term $I_1$ in \eqref{EQ177}, we use H\"older's inequality and Lemmas~\ref{L08}--\ref{L09} to get
\begin{align}
	\begin{split}
		\Vert I_1\Vert_{L^2_t H^{s-1}_x}	
		&
		\les
		\Vert \RR \nabla b \nabla \bar{v}\Vert_{L^2_t H^{s-1}_x}		
		+
		\Vert \RR b D_x^2 \bar{v}\Vert_{L^2_t H^{s-1}_x}	
		\les
		\Vert b\Vert_{L^\infty_t H^{s}_x}
		\Vert  \bar{v} \Vert_{L^2_t H^{s+1}_x}
		\les
		\TT^{1/20} \Vert \bar{v} \Vert_{K^{s+1}}
		.
		\label{EQ228}
	\end{split}
\end{align}
For the time component, we have
\begin{align}
	\begin{split}
		\Vert I_1\Vert_{H^{(s-1)/2}_t L^2_x}	
		&
		\les
		\Vert \RR \nabla b \nabla \bar{v}\Vert_{H^{(s-1)/2}_t L^2_x}	
		+
		\Vert \RR b D_x^2 \bar{v}\Vert_{H^{(s-1)/2}_t L^2_x}.	
		\label{EQ217}
	\end{split}
\end{align}
To treat the first term on the right side, we claim that for any $\alpha>1/2$ we have
\begin{align}
	\begin{split}
		\Vert A B  \Vert_{H^{\alpha}_t L^2_x}	   
		&\les
		\Vert A   \Vert_{H^\alpha_t H^1_x}
		\Vert B   \Vert_{H^{\alpha}_t H^{1/2}_x}
	\end{split}
	\label{EQ102}
\end{align}
on the domain $(0,1)\times\Omegaf$.
Using extensions, we may assume that the domain is actually $\mathbb{R}\times\mathbb{R}^{3}$.
We proceed as in \eqref{EQ98} and estimate
\begin{align}
	\begin{split}
		\Vert A B  \Vert_{H^{\alpha}_t L^2_x}
		&=    \Vert A B  \Vert_{L_{x}^{2}H^{\alpha}_t}
		\lec   \bigl\Vert
		\Vert A   \Vert_{H^{\alpha}_t}
		\Vert B   \Vert_{H^{\alpha}_t}
		\bigr\Vert_{L^2_{x}}
		\\&
		\les
		\Vert A   \Vert_{L_{x}^{6} H^{\alpha}_t}
		\Vert B   \Vert_{L_{x}^3 H^{\alpha}_t}
		\lec
		\Vert A   \Vert_{H_x^{1} H^{\alpha}_t}
		\Vert B   \Vert_{H^{1/2}_x H^{\alpha}_t}
		\\&
		=
		\Vert A   \Vert_{H^{\alpha}_t H_x^{1}}
		\Vert B   \Vert_{H^{\alpha}_t H_{x}^{1/2}}
		,
	\end{split}
   \llabel{EQ179}
\end{align}
since $\alpha >1/2$, and \eqref{EQ102} follows.
For the first term on the right side of \eqref{EQ217}, we now use \eqref{EQ102} and write
\begin{align}
	\begin{split}
	\Vert \RR \nabla b \nabla \bar{v} \Vert_{H^{(s-1)/2}_t L^2_x}	
	&\les
		\Vert \nabla \bar{v} \Vert_{H^{(s-1)/2}_t H^{1}_x}
		\Vert R\nabla b\Vert_{H^1_t H^{1/2}_x}	
	\\&
	\les		
	\Vert \bar{v}\Vert_{K^{s+1}}
	(\Vert R_t \nabla b\Vert_{L^2_t H^{1/2}_x}
	+
	\Vert R \nabla b_t \Vert_{L^2_t H^{1/2}_x}
	+
	\Vert R \nabla b\Vert_{L^2_t H^{1/2}_x}
	)
	\\&
	\les
	\TT^{1/20} 
	\Vert \bar{v}\Vert_{K^{s+1}}
	,
	\label{EQ219}
	\end{split}
\end{align}
for any $\epsilon \in (0,1]$,
where we used Corollary~\ref{C02} and Lemmas~\ref{L08}--\ref{L09}.
For the second term on the right side of \eqref{EQ217}, we appeal to \eqref{EQ100} to get
\begin{align}
	\begin{split}
		\Vert \RR b D_x^2 \bar{v}\Vert_{H^{(s-1)/2}_t L^2_x}
	&
	\les
	\Vert \RR b \Vert_{H^{(s-1)/2}_t H^{3/2+\delta}_x}	
	\Vert D^2_x \bar{v}\Vert_{H^{(s-1)/2}_t L^2_x}		
	\les
\Vert \RR b \Vert_{H^{1}_t H^{3/2+\delta}_x}	
	\Vert  \bar{v}\Vert_{K^{s+1}}	
	\\&
	\les
	\TT^{1/20} \Vert  \bar{v}\Vert_{K^{s+1}},
	\label{EQ215}
	\end{split}
\end{align}
for $\delta \in (0,1/2)$, where we used Corollary~\ref{C02} and Lemmas~\ref{L08}--\ref{L09}.
Combining \eqref{EQ228}--\eqref{EQ217} and \eqref{EQ219}--\eqref{EQ215}, we obtain
\begin{align}
	\Vert I_1 \Vert_{K^{s-1}}
	\les
	\TT^{1/20}
	\Vert \bar{v} \Vert_{K^{s+1}}
	.
	\llabel{EQ240}
\end{align}
The terms $I_2$, $I_3$, $I_4$, $I_5$, and $I_6$ are estimated analogously to $I_1$, and we get
\begin{align}
	\Vert I_2\Vert_{K^{s-1}} 
	+
	\Vert I_3 \Vert_{K^{s-1}} 
	+
	\Vert I_4 \Vert_{K^{s-1}} 
	+
	\Vert I_5 \Vert_{K^{s-1}} 
	+
	\Vert I_6 \Vert_{K^{s-1}} 
	\les
	\TT^{1/20}
	 \Vert \bar{v} \Vert_{K^{s+1}}
	.
	\llabel{EQ229}
\end{align}
For the term $I_7$, we use H\"older's inequality and obtain
\begin{align}
	\begin{split}
	\Vert I_7\Vert_{L^2_t H^{s-1}_x}
	&
	\les
	\Vert \RR^{-1} b \nabla \RR\Vert_{L^2_t H^{s-1}_x}
	\les
	\Vert \RR^{-1} \Vert_{L^\infty_t H^s_x}
	\Vert b \Vert_{L^\infty_t H^s_x}
	\Vert \RR \Vert_{L^\infty_t H^{s}_x}
	\les
	1
	\llabel{EQ243}
	\end{split}
\end{align}
and
\begin{align}
	\begin{split}
		\Vert I_7\Vert_{H^{(s-1)/2}_t L^2_x}
		&
		\les
		\Vert \RR^{-1} b \nabla \RR\Vert_{H^{(s-1)/2}_t L^2_x}
		\les
		\Vert \RR^{-1} b \nabla \RR\Vert_{H^{1}_t L^2_x}
		\\&
		\les
		\Vert \RR^{-1} \Vert_{H^1_t L^\infty_x}
		\Vert b \Vert_{L^\infty_t L^\infty_x}
		\Vert \nabla \RR \Vert_{L^\infty_t L^2_x}
		+
		\Vert \RR^{-1} \Vert_{L^\infty_t L^\infty_x}
		\Vert b \Vert_{L^\infty_t L^\infty_x}
		\Vert \nabla \RR \Vert_{H^1_t L^2_x}
		\\&\indeq
		+
		\Vert \RR^{-1} \Vert_{L^\infty_t L^\infty_x}
		\Vert b \Vert_{H^1_t L^\infty_x}
		\Vert \nabla \RR \Vert_{L^\infty_t L^2_x}
		\les 
	1,
		\llabel{EQ244}
	\end{split}
\end{align}
where we appealed to Lemmas~\ref{L08}--\ref{L09}.
The term $I_8$ is estimated analogously to $I_7$, leading to
\begin{align}
	\begin{split}
	\Vert I_8\Vert_{K^{s-1}}
	\les
	1
	.
   \llabel{EQ213}
	\end{split}
\end{align}
Using the estimates on $I_1$--$I_8$ in \eqref{EQ177}, we conclude that
\begin{align}
	\Vert f\Vert_{K^{s-1}}
	\les
	\TT^{1/20}
	\Vert \bar{v}\Vert_{K^{s+1}}
	+
	 1
	.
	\label{EQ218}
\end{align}

Next, we bound the $K^{s-1/2} (\Gammac)$ norm of the terms on the right side of \eqref{EQ178}, for
every fixed $j=1,2,3$.
For $K_1$, we use H\"older's and trace inequalities
along with Lemma~\ref{L08}
to obtain
\begin{align}
	\begin{split}
		\Vert K_1\Vert_{L^2_t H_x^{s-1/2} (\Gammac)}
		&
		\les
		\Vert (1-\JJ) \nabla \bar{v} \Vert_{L^2_t H^s_x}
		\les
		\Vert 1-\JJ \Vert_{L^\infty_t H^s_x}
		\Vert \bar{v} \Vert_{L^2_t H^{s+1}_x}
		\les
		\TT^{1/20}
		\Vert \bar{v} \Vert_{K^{s+1}}
		.
	\end{split}
		\label{EQ220}
\end{align}
For the time component, we appeal to Corollary~\ref{C01}, obtaining
\begin{align}
	\begin{split}
		\Vert K_1\Vert_{H^{s/2-1/4}_t L_x^2 (\Gammac)}
		&
		\les
		\epsilon_1
		\Vert (1-\JJ) \nabla \bar{v} \Vert_{H^{s/2}_t L^2_x}
		+
		\epsilon_1^{1-2s}
		\Vert (1-\JJ) \nabla \bar{v} \Vert_{L^2_t H^s_x}
		\\&
\les
		\epsilon_1
		\Vert \nabla \bar{v} \Vert_{H^{s/2}_t L^2_x}
		\Vert 1-\JJ \Vert_{H^{1}_t H^{3/2+\delta}_x}
		+
		\epsilon_1
		\Vert \nabla \bar{v} \Vert_{H^{(s-1)/2}_t H^1_x}
		\Vert 1-\JJ \Vert_{H^{s/2}_t H^{1/2}_x}
		\\&\indeq
		+
		\epsilon_1^{1-2s}
		\Vert 1-\JJ \Vert_{L^\infty_t H^s_x}
		\Vert\nabla \bar{v} \Vert_{L^2_t H^s_x}
		=:
		K_{11}
		+
		K_{12}
		+
		K_{13}
		,
		\label{EQ231}
	\end{split}
\end{align}
for any $\epsilon_1 \in (0,1]$, where $\delta \in (0,1/5)$.
For the term $K_{11}$, we use Corollary~\ref{C02} and Lemma~\ref{L08}
and obtain
\begin{align}
	\begin{split}
		K_{11}
		\les	
		\epsilon_1
		\Vert \bar{v}\Vert_{K^{s+1}}
		.
		\label{EQ922}
	\end{split}
\end{align}
Similarly, the term $K_{12}$ is estimated as
\begin{align}
	\begin{split}
		K_{12}
		\les
		\epsilon_1
		\Vert \bar{v}\Vert_{K^{s+1}}
		\Vert 1-J \Vert_{H^{s/2}_t H^{1/2}_x}.	
		\label{EQ925}
	\end{split}
\end{align}
From \eqref{EQ910}, Corollary~\ref{C02}, and Lemma~\ref{L08}, it follows that
\begin{align}
	\begin{split}
		\Vert \JJ_t \Vert_{H^{(s-2)/2}_t H^{1/2}_x}	
		&
		\les
		\Vert \phi_{\TT} \Vert_{H^{(s-2)/2}_t}
		\Vert J a \nabla v\Vert_{H^{(s-1)/2}_t H^{1/2}_x}
		\\&
		\les
		\Vert \phi_{\TT} \Vert_{ H^{(s-2)/2}_t}
	\Vert  J \Vert_{H^{(s-1)/2}_t H^{3/2+\delta}_x}
	\Vert a \Vert_{H^{(s-1)/2}_t H^{3/2+\delta}_x}
		\Vert \nabla v \Vert_{H^{(s-1)/2}_t H^{1}_x}
		\\&
		\les
			\Vert \phi_{\TT} \Vert_{ H^{(s-2)/2}_t}
		M
		,
		\label{EQ233}
	\end{split}
\end{align}
since $1/2<(s-1)/2<1$ and $\delta \in (0,1/5)$.
Using Lemma~\ref{L14}
in \eqref{EQ233} and applying the resulting inequality in
\eqref{EQ925}, we get
\begin{align}
	\begin{split}
		K_{12}
		&\les
		\epsilon_1
		 \Vert \bar{v} \Vert_{K^{s+1}}
		( \Vert J_t\Vert_{H^{(s-2)/2}_t H^{1/2}_x}	
		 +
		  \Vert J\Vert_{H^{(s-2)/2}_t H^{1/2}_x}	
		 )
		 \les
		 \epsilon_1
		 \Vert \bar{v} \Vert_{K^{s+1}}
                 (
		M
		+
		1)
		\les
		\epsilon
		\Vert \bar{v} \Vert_{K^{s+1}}
		,
		\label{EQ923}
	\end{split}
\end{align}
where $\epsilon \in (0,1]$,
by taking
$\epsilon_1= \epsilon M^{-1}$.
For the term $K_{13}$, we have
\begin{align}
	\begin{split}
		K_{13}
		\les	
		C_\epsilon
		M^{2s-1}
		\TT^{1/20}
		\Vert \bar{v}\Vert_{K^{s+1}}
		,	
		\label{EQ924}
	\end{split}
\end{align}
for $\TT>0$ sufficiently small. 
Combining \eqref{EQ220}--\eqref{EQ922} and \eqref{EQ923}--\eqref{EQ924}, we arrive at
\begin{align}
	\begin{split}
	\Vert K_1\Vert_{K^{s-1/2}_{\Gammac}}
	&
	\les
	\epsilon
	\Vert \bar{v}\Vert_{K^{s+1}}
	,
	\label{EQ234}
	\end{split}
\end{align}
for any $\epsilon \in (0,1]$, by taking $\TT>0$ sufficiently small.
The term $K_2$ is estimated analogously to $K_1$, and we obtain
\begin{align}
	\begin{split}
		\Vert K_2 \Vert_{K^{s-1/2}_{\Gammac} }
		\les
	\epsilon
	\Vert \bar{v}\Vert_{K^{s+1}}
		.
		\label{EQ928}
	\end{split}
\end{align}
For the space component of the term $K_3$, we use H\"older's and trace inequalities to obtain
\begin{align}
	\begin{split}
	\Vert K_3\Vert_{L^2_t H_x^{s-1/2}}
	&
	\les
	\Vert \JJ \Vert_{L^\infty_t H^s_x}
	\Vert b\Vert_{L^\infty_t H^s_x}^2 
	\Vert \bar{v} \Vert_{L^2_t H^{s+1}_x}
	\les
	\TT^{1/20}
	 \Vert \bar{v} \Vert_{K^{s+1}}
	 ,
	\label{EQ238}
\end{split}
\end{align}
where we appealed to Lemma~\ref{L08}.
For the time component, using Corollary~\ref{C01}, we have
\begin{align}
	\begin{split}
	\Vert K_3\Vert_{H^{s/2-1/4}_t L_x^2 ( \Gammac)}
	&
	\les
	\epsilon_1
	\Vert \JJ b b \nabla \bar{v} \Vert_{H^{s/2}_t L^2_x}
	+
	\epsilon_1^{1-2s}
	\Vert \JJ b b \nabla \bar{v} \Vert_{L^2_t H^s_x}
	\\&
	\les
	\epsilon_1
	\Vert \JJ \Vert_{H^{s/2}_t H^{1/2}_x}
	\Vert \nabla \bar{v} \Vert_{H^{(s-1)/2}_t H^1_x}
	+
	\epsilon_1
	\Vert b\Vert_{H^{s/2}_t H^{1/2}_x}
	\Vert \nabla \bar{v} \Vert_{H^{(s-1)/2}_t H^1_x}
			\\&\indeq
	+
	\epsilon_1
	\Vert \nabla \bar{v} \Vert_{H^{s/2}_t L^2_x}
	+
	\epsilon_1^{1-2s}
	\Vert b\Vert_{L^\infty_t H^s_x}
	\Vert \bar{v} \Vert_{L^2_t H^{s+1}_x}
	=:
	K_{31}+K_{32}+K_{33}+K_{34}
	,
	\label{EQ237}
	\end{split}
\end{align}
for any $\epsilon_1 \in (0,1]$.
The term $K_{31}$ is estimated analogously to \eqref{EQ925}--\eqref{EQ923}, and we obtain
\begin{align}
	\begin{split}
	K_{31}
	\les
	\epsilon
	\Vert \bar{v}\Vert_{K^{s+1}}
	,	
	\end{split}
   \llabel{EQ184}
\end{align}
by taking $\epsilon_1= \epsilon  M^{-1}$ in \eqref{EQ237}, where $\epsilon \in (0,1]$ is a constant. 
The term $\Vert b\Vert_{H^{s/2}_t H^{1/2}_x}$ is estimated analogously to \eqref{EQ925}--\eqref{EQ923}, and we get
\begin{align}
	\begin{split}
	\Vert b \Vert_{H^{s/2}_t L^3_x}
	&
	\les
	M+1
	.
	\end{split}
   \llabel{EQ185}
\end{align}
Therefore, we infer that
\begin{align}
	\begin{split}
	K_{32}
	\les
	\epsilon
	\Vert \bar{v}\Vert_{K^{s+1}}
	.
	\end{split}
   \llabel{EQ194}
\end{align}
The term $K_{33}$ is estimated using Corollary~\ref{C02} as
\begin{align}
	\begin{split}
	K_{33}
	\les
	\epsilon
	\Vert \bar{v}\Vert_{K^{s+1}}
	,
	\end{split}
   \llabel{EQ195}
\end{align}
while the term $K_{34}$ is estimated analogously to \eqref{EQ924} as
\begin{align}
	\begin{split}
	K_{34}
	\les
	\epsilon
	\Vert \bar{v}\Vert_{K^{s+1}},
	\llabel{EQ927}
	\end{split}
\end{align}
by taking $\TT>0$ sufficiently small.
Combining \eqref{EQ238}--\eqref{EQ237} and the estimates on $K_{31}$--$K_{34}$, we conclude that
\begin{align}
	\begin{split}
	\Vert K_3\Vert_{K^{s-1/2}_{\Gammac}}
	\les
	\epsilon
	\Vert \bar{v}\Vert_{K^{s+1}}
	,
   \label{EQ52}
	\end{split}
\end{align}
for any $\epsilon \in (0,1]$.
Regarding the term $K_4$, we proceed as in \eqref{EQ189}--\eqref{EQ190} to obtain
\begin{align}
	\begin{split}
	\Vert K_4\Vert_{K^{s-1/2}_{\Gammac}}
	& 
	\les
	\Vert \JJ b \RR^{-1} \Vert_{L^2_t H^s_x}
	+
	\Vert \JJ b \RR^{-1} \Vert_{H^1_t H^1_x}
	\les
	1
	+
	\Vert \JJ_t \Vert_{L^2_t H^1_x}
	+
	\Vert b_t \Vert_{L^2_t H^1_x}
	+
	\Vert \RR_t \Vert_{L^2_t H^1_x}
	\les
	1
	,
	\end{split}
   \label{EQ53}
\end{align}
where we used Lemmas~\ref{L08}--\ref{L09}.
The term $K_5$ is estimated in a similar fashion as $K_4$, and we arrive at
\begin{align}
	\begin{split}
	\Vert K_5 \Vert_{K^{s-1/2}_{\Gammac}}
	\les
	1
	.
	\end{split}
   \label{EQ54}
\end{align}
The terms $K_6$, $K_7$, $K_8$, $K_9$, and $K_{10}$ are estimated analogously to $K_3$, and we have
\begin{align}
	\begin{split}
		&
	\Vert K_6\Vert_{K^{s-1/2}_{\Gammac}}
	+
	\Vert K_7 \Vert_{K^{s-1/2}_{\Gammac}}
	+
	\Vert K_7 \Vert_{K^{s-1/2}_{\Gammac}}
	+
	\Vert K_8 \Vert_{K^{s-1/2}_{\Gammac}}
	+
	\Vert K_9 \Vert_{K^{s-1/2}_{\Gammac}}
	+
	\Vert K_{10} \Vert_{K^{s-1/2}_{\Gammac}}
	\les
	\epsilon
	\Vert \bar{v}\Vert_{K^{s+1}}
	,
	\label{EQ239}
	\end{split}
\end{align}
for any $\epsilon \in (0,1]$.
For the term $K_{11}$, we proceed as in \eqref{EQ189}--\eqref{EQ190} using Lemma~\ref{L09} to obtain
\begin{align}
	K_{11}
	\les
	1
	.
   \label{EQ55}
\end{align}
Collecting the estimates \eqref{EQ234}--\eqref{EQ928} and \eqref{EQ52}--\eqref{EQ55}, we conclude
\begin{align}
	\begin{split}
	\Vert h\Vert_{K^{s-1/2}_{\Gammac}}
	\les
	\epsilon
	\Vert \bar{v}\Vert_{K^{s+1}}
	+
	1
	,
	\label{EQ930}
	\end{split}
\end{align}
for any $\epsilon \in (0,1]$.
For the fourth term on the right side of \eqref{EQ201}, we have
\begin{align}
	\begin{split}
	\Vert f(0)\Vert_{H^{s-2}}
	\les	
	\Vert R_0^{-1} \nabla R_0 \Vert_{H^{s-2}}
	\les
	1.
	\label{EQ1111}
	\end{split}
\end{align}
From \eqref{EQ201}--\eqref{EQ931}, \eqref{EQ218}, and \eqref{EQ930}--\eqref{EQ1111} it follows that
\begin{align}
	\begin{split}
		\Vert \bar{v}\Vert_{K^{s+1} }
		&
		\les	
	(\epsilon
	+
	\TT^{1/20})
	\Vert \bar{v}\Vert_{K^{s+1}}
	+
	(\epsilon+\tilde{\epsilon} C_\epsilon
	+ 
	\TT^{1/2} C_{\tilde{\epsilon}, \epsilon} ) M
	+
	C_\epsilon
	,
	\end{split}
   \llabel{EQ196}
\end{align}
for any $\epsilon, \tilde{\epsilon} \in (0,1]$.
We first take $\epsilon$ sufficiency small, and $\tilde{\epsilon}$ sufficiently small depending on $\epsilon$, and then $\TT$ sufficiently small depending on $\epsilon, \tilde{\epsilon}$, yielding
\begin{align}
	\begin{split}
	\Vert \bar{v}\Vert_{K^{s+1}}
	\leq
	M,	
	\label{EQ932}
	\end{split}
\end{align}
by allowing $M\geq 1$ sufficiently large.
Thus, the mapping $\Pi \colon v\mapsto \bar{v}$ is well-defined from $\mathcal{Z}$ to $\mathcal{Z}$, for some constant $M\geq 1$, which depends on the size of the initial data.

\subsection{Contracting property}
In this section, we prove
\begin{align}
	\Vert \Pi (v_1) - \Pi (v_2) \Vert_{K^{s+1}}
	\leq
	\frac{1}{2} \Vert v_1 - v_2 \Vert_{K^{s+1}}
	\comma
	v_1, v_2 \in \mathcal{Z},
	\label{EQ298}
\end{align}
where $M\geq 1$ is fixed as in~\eqref{EQ932}
and $\TT $ sufficiently small.
Note that the implicit constants below are allowed to depend on~$M$.
Let $\TT>0$ be a sufficiently small constant such that Lemmas~\ref{L08}--\ref{L09} hold.

Let $v_1, v_2 \in \mathcal{Z}$ and $(\eta_1, \eta_2)$ be the corresponding Lagrangian flow maps as in~\eqref{EQ999}.
Denote by
$(\JJ_1, a_1)$ and $(\JJ_2, a_2)$ the Jacobians and the inverse matrices of the corresponding flow map.
First we solve for $(\RR_1, \RR_2)$ from \eqref{EQ940}--\eqref{EQ941} with the same initial data~$\RR_0$.
Then we solve for $(\xi_1, \xi_{1t})$ and $(\xi_2, \xi_{2t})$ from \eqref{EQ72} with the boundary conditions \eqref{EQ73}--\eqref{EQ76} and the same initial data $(w_0, w_1)$.
To obtain the next iterate $(\bar{v}_1, \bar{v}_2)$, we solve \eqref{EQ124} with the boundary conditions and the initial data~\eqref{EQ125}.
Denote $b_1 = a_1 - \mathbb{I}_3$, $b_2 = a_2 - \mathbb{I}_3$, 
 $\tilde{b} = b_1 - b_2$, $\tilde{V} = \bar{v}_1 - \bar{v}_2$, $\tilde{v} = v_1 - v_2$, $\tilde{\RR} = \RR_1 - \RR_2$, $\tilde{\xi} = \xi_1 - \xi_2$, $\tilde{\eta}= \eta_1 -\eta_2$, and $\tilde{\JJ} = \JJ_1 - \JJ_2$.
The difference $\tilde{\VV}$ satisfies
\begin{align}
	\begin{split}
	&
	\partial_t \tilde{\VV}_j 
	-
	\lambda \tilde{\RR} \partial_k 
	(\partial_j \tilde{\VV}_k 
	+ 
	\partial_k \tilde{\VV}_j)
	-
	\mu \tilde{\RR} \partial_j \partial_k \tilde{\VV}_k
	=
	\tilde{f}_j
	\inon{in~$[0,1]\times \Omegaf$}
	,
	\\
	&
	\lambda (\partial_k \tilde{\VV}_j 
	+ 
	\partial_j \tilde{\VV}_k) \nu^k
	+
	\mu \partial_k \tilde{\VV}_k \nu^j
	=
	\partial_k \tilde{\xi}_j \nu^k
	+
	\tilde{h}_j
	\inon{in~$[0,1] \times \Gammac$}
	,
		\\&
	\tilde{\VV} ~\text{periodic in the $y_1$ and $y_2$ directions}
	,
	\\&
	\tilde{\VV}(0) = 0
	\inon{in~$\Omegaf$}
	,
	\label{EQ304}
\end{split}
\end{align}
where
\begin{align}
	\begin{split}
	\tilde{f}_j
	&
	=	
	\lambda \tilde{\RR} \partial_k
	(\partial_j \bar{v}_{2k} 
	+ 
	\partial_k \bar{v}_{2j})
	+
	\mu \tilde{\RR} \partial_j \partial_k \bar{v}_{2k}
	+
	\lambda \tilde{\RR} \partial_k 
	(b_{1mk} \partial_m \bar{v}_{1j} 
	+ 
	b_{1mj} \partial_m \bar{v}_{1k})
	\\&\indeq
	+
	\lambda \RR_2 \partial_k
	(b_{1mk} \partial_m \tilde{\VV}_j
	+
	b_{1mj} \partial_m \tilde{\VV}_k)
	+
	\lambda \RR_2 \partial_k 
	(\tilde{b}_{mk} \partial_m \bar{v}_{2j}
	+
	\tilde{b}_{mj} \partial_m \bar{v}_{2k})
	\\&\indeq
	+
	\lambda \tilde{\RR} b_{1kl} \partial_k
	(b_{1ml} \partial_m \bar{v}_{1j} 
	+ 
	b_{1mj} \partial_m \bar{v}_{1l})
	+
	\lambda \RR_2 \tilde{b}_{kl} 
	\partial_k (b_{1ml} \partial_m \bar{v}_{1j}
	+
	b_{1mj} \partial_m \bar{v}_{1l}
	)
	\\&\indeq
	+
	\lambda \RR_2 b_{2kl} \partial_k 
	(\tilde{b}_{ml} \partial_m \bar{v}_{1j}
	+ \tilde{b}_{mj} \partial_m \bar{v}_{1l}
	)
	+
	\lambda \RR_2 b_{2kl} \partial_k
	( b_{2ml} \partial_m \tilde{\VV}_j
	+
	b_{2mj} \partial_m \tilde{\VV}_l
	)
	\\&\indeq
	+
	\lambda \tilde{\RR} b_{1kl} \partial_k
	(\partial_l \bar{v}_{1j}
	+
	\partial_j \bar{v}_{1l}
	)
	+
	\lambda \RR_2 \tilde{b}_{kl} \partial_k 
	(\partial_l \bar{v}_{1j}
	+
	\partial_j \bar{v}_{1l})
	+
	\lambda \RR_2
	b_{2kl} \partial_k 
	(\partial_l \tilde{\VV}_j
	+
	\partial_j \tilde{\VV}_l
	)
	\\&\indeq
	+
	\mu \RR \partial_j
	(b_{1mi} \partial_m \bar{v}_{1i})
	+
	\mu \RR_2 \partial_j (\tilde{b}_{mi} \partial_m \bar{v}_{1i})
	+
	\mu \RR_2 \partial_j (b_{2mi} \partial_m \tilde{\VV})
	\\&\indeq
	+
	\mu \tilde{\RR} b_{1kj} \partial_k
	(b_{1mi} \partial_m \bar{v}_{1i})
	+
	\mu \RR_2 \tilde{b}_{kj} \partial_k
	(b_{1mi} \partial_m \bar{v}_{1i})
	+
	\mu \RR_2 b_{2kj} \partial_k
	( \tilde{b}_{mi} \partial_m \bar{v}_{1i})
	\\&\indeq
	+
	\mu \RR_2 b_{2kj} \partial_k
	(b_{2mi} \partial_m \tilde{\VV}_{i})
	+
	\mu \tilde{\RR} b_{1kj} \partial_k \partial_i \bar{v}_{1i}
	+
	\mu \RR_2 \tilde{b}_{kj} \partial_k \partial_i \bar{v}_{1i}
	+
	\mu \RR_2 b_{2kj} \partial_k \partial_i \tilde{\VV}_{i}
	\\&\indeq
	-
	\RR_1^{-1} \RR_2^{-1} \tilde{\RR} b_{1kj} \partial_k \RR_1
	+
	\RR_2^{-1} \tilde{b}_{kj} \partial_k \RR_1
	+
	\RR_2^{-1} b_{2kj} \partial_k \tilde{\RR}
	-
	\RR_1^{-1} \RR_2^{-1} \tilde{\RR} \partial_j \RR_1
	\\&\indeq
	+
	\RR_2^{-1} \partial_j \tilde{\RR}
	\label{EQ248}
	\end{split}
\end{align}
and
\begin{align}
	\begin{split}
	\tilde{h}_j
	&
	=	
	-
	\lambda \tilde{\JJ} 
	(\partial_k \bar{v}_{1j} 
	+ 
	\partial_j \bar{v}_{1k}
	) \nu^k
	+
	\lambda (1- \JJ_2) 
	(\partial_k \tilde{\VV}_j + \partial_j \tilde{\VV}_k) \nu^k
	-
	\mu \tilde{\JJ} \partial_k \bar{v}_{1k} \nu^j
	+
	\mu (1- \JJ_2) \partial_k \tilde{\VV}_k \nu^j
	\\&\indeq
	+
	\lambda \tilde{\JJ} b_{1kl}
	(b_{1ml} \partial_m \bar{v}_{1j} 
	+
	b_{1mj} \partial_m \bar{v}_{1l}
	) \nu^k
	+
	\lambda \JJ_2 \tilde{b}_{kl}
	(b_{1ml} \partial_m \bar{v}_{1j} 
	+
	b_{1mj} \partial_m \bar{v}_{1l}
	) \nu^k
	\\&\indeq
	+
	\lambda \JJ_2 b_{2kl}
	(\tilde{b}_{ml} \partial_m \bar{v}_{1j} 
	+
	\tilde{b}_{mj} \partial_m \bar{v}_{1l}
	) \nu^k
	+
	\lambda \JJ_2 b_{2kl}
	(b_{2ml} \partial_m \tilde{\VV}_{j} 
	+
	b_{2mj} \partial_m \tilde{\VV}_{l}
	) \nu^k
	\\&\indeq
	+
	\tilde{\JJ} b_{1kj} \RR_1^{-1} \nu^k
	+
	\JJ_2 \tilde{b}_{kj} \RR_1^{-1} \nu^k
	-
	\JJ_2 b_{2kj} \RR_1^{-1} \RR_2^{-1} \tilde{\RR} \nu^k
	-
	\tilde{\JJ} \RR_1^{-1} \nu^j
	-
	(\JJ_2-1) \RR_1^{-1} \RR_2^{-1} \tilde{\RR} \nu^j
	\\&\indeq
	-
	\lambda \tilde{\JJ} 
	(b_{1mk} \partial_m \bar{v}_{1j} + b_{1mj} \partial_m \bar{v}_{1k}) \nu^k
	-
	\lambda \JJ_2
	(\tilde{b}_{mk} \partial_m \bar{v}_{1j} 
	+ 
	\tilde{b}_{mj} \partial_m \bar{v}_{1k}) \nu^k
	\\&\indeq
	-
	\lambda \JJ_2
	(b_{2mk} \partial_m \tilde{\VV}_{j} 
	+ 
	b_{2mj} \partial_m \tilde{\VV}_{1k}) \nu^k
	-
	\lambda \tilde{\JJ} b_{1kl} 
	(\partial_l \bar{v}_{1j} + \partial_j \bar{v}_{1l} )\nu^k
	-
	\lambda \JJ_2 \tilde{b}_{kl}
	(\partial_l \bar{v}_{1j} + \partial_j \bar{v}_{1l} )\nu^k
	\\&
	\indeq
	-
	\lambda \JJ_2 b_{2kl}
	(\partial_l \tilde{\VV}_{j} 
	+ 
	\partial_j \tilde{\VV}_{l} )\nu^k
	-
	\mu \tilde{\JJ} b_{1kj} b_{1mi} \partial_m \bar{v}_{1i} \nu^k
	-
	\mu \JJ_2 \tilde{b}_{kj} b_{1mi} \partial_m \bar{v}_{1i} \nu^k
	\\&\indeq
	-
	\mu \JJ_2 b_{2kj} \tilde{b}_{mi} \partial_m \bar{v}_{1i} \nu^k
	-
	\mu \JJ_2 b_{2kj} b_{2mi} \partial_m \tilde{\VV}_{i} \nu^k
	-
	\mu \tilde{\JJ} b_{1mi} \partial_m \bar{v}_{1i} \nu^j
	-
	\mu \JJ_2 \tilde{b}_{mi} \partial_m \bar{v}_{1i} \nu^j
	\\&
	\indeq
	-
	\mu \JJ_2 b_{2mi} \partial_m \tilde{\VV}_{i} \nu^j
	-
	\mu \tilde{\JJ} b_{1kj} \partial_i \bar{v}_{1i} \nu^k
	-
	\mu \JJ_2 \tilde{b}_{kj} \partial_i \bar{v}_{1i} \nu^k
	-
	\mu \JJ_2 b_{2kj} \partial_i \tilde{\VV}_{i} \nu^k
		-
	\RR_1^{-1} \RR_2^{-1} \tilde{\RR} \nu^j
	,
	\label{EQ249}
	\end{split}
\end{align}
for  $j=1,2,3$.

Before we bound the terms on the right sides of \eqref{EQ248} and~\eqref{EQ249}, we provide necessary a~priori estimates for the differences of
densities, Jacobians, and inverse matrices of the flow map.

\cole
\begin{Lemma}
\label{L10}
Let $v_1, v_2 \in \mathcal{Z}$.
Suppose $\Vert v_1\Vert_{K^{s+1} ((0,1)\times \Omegaf)} \leq M$ and $\Vert v_2\Vert_{K^{s+1} ((0,1)\times \Omegaf)} \leq M$, where $M\geq 1$ is fixed as in~\eqref{EQ245}.
Let $\delta \in (0,1/5)$.
Then, for $\TT>0$ sufficiently small depending on $\delta$, we have
\begin{enumerate}[label=(\roman*)]
	\item $\Vert \tilde{b} \Vert_{L^\infty_t H^s_x} + \Vert \tilde{b} \Vert_{H^1_t H^{3/2+\delta}_x} \les \TT^{1/20} \Vert \tilde{v} \Vert_{K^{s+1}}$,
	\item $\Vert \tilde{\RR} \Vert_{L^\infty_t H^s_x} + \Vert \tilde{\RR} \Vert_{H^1_t H^{3/2+\delta}_x}
	\les \TT^{1/20} \Vert \tilde{v} \Vert_{K^{s+1}}$,
	\item $\Vert \tilde{\JJ} \Vert_{L^\infty_t H^s_x} +
	\Vert \tilde{\JJ} \Vert_{H^1_t H^{3/2+\delta}_x}
	 \les \TT^{1/20} \Vert \tilde{v} \Vert_{K^{s+1}}$,
	\item $\Vert \tilde{\RR} \Vert_{H^1_t H^s_x} + \Vert \tilde{b} \Vert_{H^1_t H^s_x} + \Vert \tilde{\JJ} \Vert_{H^1_t H^s_x} \les \Vert \tilde{v} \Vert_{K^{s+1}}$,
\end{enumerate}
for any $\delta \in (0,1)$, where the norm of dependence is $(0,1)\times \Omegaf$.
\end{Lemma}
\colb

\begin{proof}[Proof of Lemma~\ref{L10}]
(i) From \eqref{EQ211} it follows that the difference $\tilde{b}$ satisfies
\begin{align}
	\begin{split}
	-\tilde{b}_t
	&
	=
	\phi_{\TT}
	\Bigl(
    	\tilde{b} ( \nabla v_1 ) b_1
  	+
  	b_2 ( \nabla \tilde{v} ) b_1
  	+
  	b_2 ( \nabla v_2 ) \tilde{b}
  	+
  	(\nabla \tilde{v} ) b_1
	+
	(\nabla v_2 ) \tilde{b}
	+
	\tilde{b} ( \nabla v_1 )
	\\&\indeq\indeq\indeq\indeq\indeq\indeq
	+
	b_2 ( \nabla \tilde{v} )
	+
	\nabla \tilde{v}
	\Bigr)
	\inon{in~$[0, 1]\times \Omegaf$}
	,
	\label{EQ251}
	\end{split}
\end{align}	
with the initial data $\tilde{b}(0) = 0$.
Using the fundamental theorem of calculus, we obtain that for $t\in (0, 2\TT)$
\begin{align}
	\begin{split}
	\Vert \tilde{b} (t) \Vert_{H^s}	
	&
	\les
	\int_0^t \Vert \tilde{b} \Vert_{H^s}
	\Vert \nabla v_1 \Vert_{H^s}
	\Vert b_1 \Vert_{H^s}
	+
	\int_0^t \Vert b_2 \Vert_{H^s}
	\Vert \nabla \tilde{v} \Vert_{H^s}
	\Vert b_1 \Vert_{H^s}
	\\&\indeq
	+
	\int_0^t \Vert b_2 \Vert_{H^s}
	\Vert \nabla v_2 \Vert_{H^s}
	\Vert \tilde{b} \Vert_{H^s}
	+
	\int_0^t \Vert \nabla \tilde{v} \Vert_{H^s}
	\Vert b_1 \Vert_{H^s}
	+
	\int_0^t \Vert \nabla \tilde{v} \Vert_{H^s}
	\\&\indeq
	+
	\int_0^t 
	\Vert \nabla v_2 \Vert_{H^s}
	\Vert \tilde{b} \Vert_{H^s}
	+
		\int_0^t \Vert \tilde{b} \Vert_{H^s}
	\Vert \nabla v_1 \Vert_{H^s}
	+
	\int_0^t \Vert b_2 \Vert_{H^s}
	\Vert \nabla \tilde{v} \Vert_{H^s}
	\\&
	\les
	\int_0^t \Vert \tilde{v} \Vert_{H^{s+1}}
		+
	\int_0^t \Vert \tilde{b} \Vert_{H^s}
	(\Vert  v_1\Vert_{H^{s+1}}
	+ 
	\Vert v_2\Vert_{H^{s+1}})
	,
	\end{split}
   \llabel{EQ197}
\end{align}
where the last inequality follows from Lemma~\ref{L08}.
Using Gronwall's inequality, we arrive at
\begin{align}
	\begin{split}
	\Vert \tilde{b} \Vert_{L^\infty  ((0,2\TT), H^s (\Omegaf))}
	\les
	\TT^{1/2}
	\Vert \tilde{v} \Vert_{K^{s+1}}
	.
	\label{EQ255}
	\end{split}
\end{align}
Therefore, we have
\begin{align}
	\Vert \tilde{b} \Vert_{L^\infty_t H^s_x}
	\les
	\TT^{1/2}
	\Vert \tilde{v} \Vert_{K^{s+1}}
	,
	\label{EQ003}
\end{align}
since $\tilde{b}_t = 0$ on~$(2\TT,1)$.
From \eqref{EQ251} and H\"older's and the Sobolev inequalities it follows that
\begin{align}
	\begin{split}
	\Vert \tilde{b}_t \Vert_{L^2_t H^{3/2+\delta}_x}
	&
	\les
 	\Vert \tilde{b}\Vert_{L^\infty_t H^{3/2+\delta}_x}
	+
	\Vert \nabla \tilde{v} \Vert_{L^2 ((0, 2\TT), H^{3/2+\delta}  (\Omegaf))}
	\\&
	\les
	\TT^{1/2}
	\Vert \tilde{v} \Vert_{K^{s+1} }
	+
	\epsilon_1 \Vert \tilde{v} \Vert_{L^2_t H^{s+1}_x }
	+
	\epsilon_1^{(3+2\delta)/(3+2\delta-2s)} \Vert \tilde{v} \Vert_{L^2 ((0, 2\TT), L^2 (\Omegaf))}
	,
	\label{EQ295}
	\end{split}
\end{align}
for any $\epsilon_1 \in (0,1]$,
where we used Corollary~\ref{C02} and Lemma~\ref{L08}.
Letting $\epsilon_1 = \TT^{1/20}$, we obtain
\begin{align}
	\Vert \tilde{b}_t \Vert_{L^2_t H^{3/2+\delta}_x}
	\les
	\TT^{1/20}
	\Vert \tilde{v}\Vert_{K^{s+1}}
	+
		\TT^{1/2+(3+2\delta)/20(3+2\delta-2s)}
	\Vert \tilde{v}\Vert_{K^{s+1}}
	\les
		\TT^{1/20}
	\Vert \tilde{v}\Vert_{K^{s+1}}
	.
	\label{EQ004}
\end{align}
Combining \eqref{EQ003} and \eqref{EQ004}, we conclude the proof of (i).

(ii) Since the difference $\tilde{R}$ satisfies the ODE system
\begin{align} 
		&
	\tilde{\RR}_t	
	-
	\phi_{\TT}
	\tilde{\RR} 
	(\dive v_2
	+
	b_{1kj} \partial_k v_{1j}
	)
	=
	\phi_{\TT} 
	(
	\RR_1 \dive \tilde{v}
	+
	\RR_2 \tilde{b}_{kj} \partial_k v_{1j}
	+
	\RR_2 b_{2kj} \partial_k \tilde{v}_j
	)
	\inon{in~$[0,1] \times \Omegaf$}
	,
	\label{EQ256}
	\\
	&
	\tilde{\RR} (0) = 0
	\inon{in~$\Omegaf$}
	,
	\label{EQ257}
\end{align}
the solution is given by
\begin{align}
	\begin{split}
		\tilde{\RR}(t,x)
		&
		=
		e^{\int_0^t \phi_{\TT} (\dive v_2 + b_{1kj} \partial_k v_{1j} ) \,d\tau}
		\int_0^t e^{-\int_0^\tau \phi_{\TT} (\dive v_2 + b_{1kj} \partial_k v_{1j} )}
		\\&	\indeqtimes 
		\phi_{\TT}
		(\RR_1  \dive \tilde{v}  
		+
		\RR_2 \tilde{b}_{kj} \partial_k v_{1j}
		+
		\RR_2 b_{2kj} \partial_k \tilde{v}_j
		)
		\,d\tau
		\inon{in~$[0, 1] \times \Omegaf$}
		.
		\llabel{EQ252}
	\end{split}
\end{align}
From H\"older's inequality, it follows that
\begin{align}
	\begin{split}
	\Vert \tilde{R} \Vert_{L^\infty_t H^s_x}
	\les	
	\int_0^{2\TT} 
	(\Vert \tilde{v} \Vert_{H^{s+1}}
	+
	\Vert  \tilde{b} \nabla v_{1} \Vert_{H^s}) \,d\tau
	\les
	\TT^{1/2} \Vert \tilde{v} \Vert_{K^{s+1}}
	,
	\label{EQ402}
	\end{split}
\end{align}
where we used~\eqref{EQ003} in the last inequality.
Using \eqref{EQ256} and H\"older's and Sobolev inequalities, we obtain
\begin{align}
	\begin{split}
	\Vert \tilde{\RR}_t \Vert_{L^2_t H^{3/2+\delta}_x}	
	&
	\les
	\Vert \tilde{\RR} \Vert_{L^\infty_t H^{3/2+\delta}_x}
	\Vert \nabla v_2 \Vert_{L^2_t H^{3/2+\delta}_x}
	+
	\Vert \tilde{\RR} \Vert_{L^\infty_t H^{3/2+\delta}_x}
	\Vert b_1 \Vert_{L^\infty_t H^{3/2+\delta}_x}
	\Vert \nabla v_1 \Vert_{L^2_t H^{3/2+\delta}_x}
			\\&\indeq
	+
	\Vert  \tilde{v} \Vert_{L^2  ((0,2\TT), H^{5/2+\delta} (\Omegaf))}
	+
	\Vert \tilde{b} \Vert_{L^\infty_t H^{3/2+\delta}_x}
	\Vert \nabla v_1 \Vert_{L^2_t H^{3/2+\delta}_x}
	+
	\Vert  \tilde{v} \Vert_{L^2  ((0,2\TT), H^{5/2+\delta} (\Omegaf))}
	.
	\label{EQ290}
	\end{split}
\end{align}
We proceed analogously to Lemma~\ref{L09} to get
\begin{align}
	\Vert \tilde{\RR}_t \Vert_{L^2_t H^{3/2+\delta}_x}	
	\les
	\TT^{1/20}
	\Vert \tilde{v} \Vert_{K^{s+1}}
	.
   \llabel{EQ198}
\end{align}
By combining \eqref{EQ402}--\eqref{EQ290},
we conclude the proof of (ii)

The proofs of (iii) and (iv) are analogous to the proofs of (i)--(iii),
and thus we omit the details.
\end{proof}

\begin{proof}[Proof of Theorem~\ref{T01}]
From Lemmas~\ref{L02} and~\ref{L09}, it follows that the solution $\tilde{V}$ of \eqref{EQ304} satisfies
\begin{align}
	\begin{split}
		\Vert \tilde{\VV}\Vert_{K^{s+1} ((0,1) \times \Omegaf)}	
		&
		\les
		\left\Vert \frac{\partial \tilde{\xi}}{\partial \nu} \right\Vert_{K^{s-1/2} ((0,1)\times \Gammac)}
		+
		\Vert \tilde{h} \Vert_{K^{s-1/2} ((0,1)\times \Gammac)}
		+
		\Vert \tilde{f} \Vert_{K^{s-1} ((0,1)\times \Omegaf)}
		,
		\label{EQ293}
	\end{split}
\end{align}
where $\tilde{f}_j$ and $\tilde{h}_j$ are as in \eqref{EQ248}--\eqref{EQ249}, for $j=1,2,3$.

For the first term on the right side of \eqref{EQ293}, we proceed as in \eqref{EQ115}--\eqref{EQ973} to obtain
\begin{align}
	\left\Vert \frac{\partial \tilde{\xi}}{\partial \nu} 
	\right\Vert_{K^{s-1/2} (\Gammac)}
	\les
			(\epsilon+\tilde{\epsilon} C_\epsilon
	+ 
	\TT^{1/2} C_{\tilde{\epsilon}, \epsilon} ) \Vert \tilde{v}\Vert_{K^{s+1}}
	,
	\label{EQ297}
\end{align}
for any $\epsilon, \tilde{\epsilon} \in (0,1]$.

Next we estimate the $K^{s-1}$ norm of terms on the right side of~\eqref{EQ248} for $j=1,2,3$.
The space component of the term $\tilde{\RR} b_{1kj} \partial_k (b_{1mi} \partial_m \bar{v}_{1i})$ is bounded as
\begin{align}
	\begin{split}
	\Vert \tilde{\RR} b_{1kj} \partial_k (b_{1mi} \partial_m \bar{v}_{1i}) \Vert_{L^2_t H^{s-1}_x}
	&
	\les
	\Vert \tilde{\RR} \Vert_{L^\infty_t H^s_x} 
	\Vert b_1 \Vert_{L^\infty_t H^s_x}^2
	\Vert \bar{v}_1 \Vert_{L^2_t H^{s+1}_x}
	\les
	\TT^{1/20} 
	\Vert \tilde{v} \Vert_{K^{s+1}}
	,
	\end{split}
   \llabel{EQ108}
\end{align}
while for the time component, we have
\begin{align}
	\begin{split}
		&
		\Vert \tilde{\RR} b_{1kj} \partial_k (b_{1mi} \partial_m \bar{v}_{1i}) \Vert_{H^{(s-1)/2}_t L^2_x}
		\les
		\Vert \tilde{\RR} b_{1} \partial_k b_{1} \nabla \bar{v}_{1} \Vert_{H^{(s-1)/2}_t L^2_x}
		+
		\Vert \tilde{\RR} b_{1} b_{1} D_x^2  \bar{v}_{1} \Vert_{H^{(s-1)/2}_t L^2_x}
		\\
		&
		\les
		\Vert \tilde{\RR} \Vert_{H^{(s-1)/2}_t H^{3/2+\delta}_x} 
		\Vert \nabla \bar{v}_1 \Vert_{H^{(s-1)/2}_t H^1_x}
		+
		\Vert \tilde{\RR} \Vert_{H^{(s-1)/2}_t H^{3/2+\delta}_x} 
		\Vert D_x^2 \bar{v}_1 \Vert_{H^{(s-1)/2}_t L^2_x}
		\\&
		\les
		\TT^{1/20}
		 \Vert \tilde{v} \Vert_{K^{s+1}}
		,
	\end{split}
   \llabel{EQ112}
\end{align}
where we used Corollary~\ref{C02} and Lemmas~\ref{L08}--\ref{L10}.
Similarly, the term $\mu \RR_2 b_{2kj} \partial_k \partial_i \tilde{\VV}_{i}$ is estimated as
\begin{align}
	\begin{split}
	\Vert \mu \RR_2 b_{2kj} \partial_k \partial_i \tilde{\VV}_{i} \Vert_{L^2_t H^{s-1}_x}
	&
	\les
	\Vert \RR_2 \Vert_{L^\infty_t H^s_x}
	\Vert b_2 \Vert_{L^\infty_t H^s_x}
	\Vert D_x^2  \tilde{\VV} \Vert_{L^2_t H^{s-1}_x}
	\les
	\TT^{1/20} 
	\Vert \tilde{\VV} \Vert_{K^{s+1}}
	\end{split}
   \llabel{EQ110}
\end{align}
and
\begin{align}
	\begin{split}
		\Vert \mu \RR_2 b_{2kj} \partial_k \partial_i \tilde{\VV}_{i} \Vert_{H^{(s-1)/2}_t L^2_x}
		&
		\les
				\Vert b_2 \Vert_{H^{(s-1)/2}_t H^{3/2+\delta}_x}
		\Vert D_x^2 \tilde{\VV} \Vert_{H^{(s-1)/2}_t L^2_x}
		\les
		\TT^{1/20}
		\Vert \tilde{\VV} \Vert_{K^{s+1}}
		.
	\end{split}
   \llabel{EQ113}
\end{align}
Other terms on the right side of \eqref{EQ248} are treated analogously
as in the proof of
Theorem~\ref{T03} using Lemmas~\ref{L08}--\ref{L10}, and we arrive at
\begin{align}
	\Vert \tilde{f} \Vert_{K^{s-1}}
	\les
	\TT^{1/20}
	\Vert \tilde{v} \Vert_{K^{s+1}}
	+
	\TT^{1/20}
	\Vert \tilde{\VV} \Vert_{K^{s+1}}
	.
	\label{EQ291}
\end{align}

Next we estimate the $K^{s-1/2}_{\Gammac}$ norm of the terms on the right side of~\eqref{EQ249}, for $j=1,2,3$. 
The term $\lambda (1- \JJ_2) \partial_k \tilde{\VV}_j \nu^k$ is estimated using the trace inequality and Lemma~\ref{L08} as
\begin{align}
	\begin{split}
        &
	\Vert \lambda (1- \JJ_2) 
	(\partial_k \tilde{\VV}_j + \partial_j \tilde{\VV}_k) \nu^k	\Vert_{L^2_t H_x^{s-1/2} (\Gammac)} 
\\&\indeq	
	\les
	\Vert (1- \JJ_2)
	\nabla \tilde{V}	\Vert_{L^2_t H^{s}_x} 
	\les
	\Vert 1- \JJ_2 \Vert_{L^\infty_t H^s_x}
	\Vert \tilde{\VV} \Vert_{L^2_t H^{s+1}_x}
	\les
	\TT^{1/20} 
	\Vert \tilde{\VV} \Vert_{K^{s+1}}
	.
	\end{split}
   \llabel{EQ222}
   \end{align}
For the time component, we proceed analogously to \eqref{EQ231}--\eqref{EQ924}, obtaining
\begin{align}
	\begin{split}
	\Vert \lambda (1- \JJ_2) 
	(\partial_k \tilde{\VV}_j + \partial_j \tilde{\VV}_k) \nu^k	\Vert_{H^{s/2-1/4}_t L^2_x (\Gammac)} 	
	&\les
	(\epsilon_1 + C_{\epsilon_1} \TT^{1/30} ) \Vert \tilde{\VV} \Vert_{K^{s+1}}
	,
	\end{split}
   \llabel{EQ199}
\end{align}
for any $\epsilon_1 \in (0,1]$.
Other terms on the right side of \eqref{EQ249} are treated analogously to Theorem~\ref{T03} using Lemmas~\ref{L08}--\ref{L10}, and we arrive at
\begin{align}
	\Vert \tilde{h} \Vert_{K^{s-1/2}_{\Gammac}}
	\les
	\TT^{1/20}
	 \Vert \tilde{v} \Vert_{K^{s+1}}
	+
	(\epsilon_1 + C_{\epsilon_1} \TT^{1/30})
	\Vert \tilde{\VV} \Vert_{K^{s+1}}
	,
	\label{EQ296}
\end{align}
for any $\epsilon_1 \in (0,1]$.

Since the terms involving $\Vert \tilde{\VV} \Vert_{K^{s+1}}$ on the right side of \eqref{EQ291}--\eqref{EQ296} are absorbed to the left side \eqref{EQ293} by taking $\epsilon_1$ and $\TT>0$ sufficiently small, we obtain from \eqref{EQ293}--\eqref{EQ296} that
\begin{align}
	\begin{split}
		\Vert \tilde{\VV}\Vert_{K^{s+1}}	
		&
		\leq
		\frac{1}{2}
		\Vert \tilde{v} \Vert_{K^{s+1}}	
		,
	\end{split}
   \llabel{EQ200}
\end{align}
by taking suitable $\epsilon$, $\tilde{\epsilon}$, and~$\TT>0$.
Thus the mapping $\Pi$ is contracting from $\mathcal{Z}$ to~$\mathcal{Z}$. Using the Banach fix point theorem, there exists a unique solution $v \in \mathcal{Z}$ such that $\Pi (v) = v$.

Fix the constant~$\TT>0$.
We proceed as in \eqref{EQ180}--\eqref{EQ116} to obtain the interior regularity estimate
\begin{align}
	\begin{split}
		&
		\Vert w \Vert_{C([0, 1], H^{s+1/4 -\epsilon_0} (\Omegae))}
		+
		\Vert w_t \Vert_{C([0, 1], H^{s-3/4-\epsilon_0}(\Omegae))}
		\leq
		C
		,
		\label{EQ202}
	\end{split}
\end{align}
where $C>0$ is a constant.
From \eqref{EQ202} and Lemma~\ref{L09} it follows that the system \eqref{EQ260}--\eqref{EQ267} admits a unique solution
\begin{align}
	\begin{split}
	(v, \RR, w, w_t)
	&
	\in
	K^{s+1} ((0,\TT)\times \Omegaf) \times
	H^1 ((0, \TT), H^s (\Omegaf))
	\\&
	\indeqtimes
	C([0, \TT], H^{s+1/4 -\epsilon_0} (\Omegae))	
	\times
	C([0, \TT], H^{s-3/4 -\epsilon_0}(\Omegae))	,
	\end{split}
   \llabel{EQ206}
\end{align}
for some constant $\TT>0$, with the corresponding norms bounded by a function of the initial data.
\end{proof}

  \colb
\section*{Acknowledgments}
IK was supported in part by the
NSF grants DMS-1907992 and DMS-2205493, while
LL was supported in part by the NSF grants  DMS-1907992, DMS-2009458, and DMS-2205493.
The work was undertaken while the authors were members of the MSRI
program ``Mathematical problems in fluid dynamics'' during the Spring~2021 semester (NSF~DMS-1928930).

%Data sharing not applicable to this article as no datasets were generated or analyzed during the current study.

\end{document}